%% file: main.tex
\documentclass{amsart}

\title[Quadratically enriched tropical intersections]{Quadratically enriched tropical intersections}
\author{Jaramillo Puentes, Andr\'es \and Pauli, Sabrina}

\usepackage[english]{babel}
\usepackage{amsmath}
\usepackage{graphicx}
\usepackage[dvipsnames]{xcolor}
\PassOptionsToPackage{unicode}{hyperref}
\PassOptionsToPackage{naturalnames}{hyperref}
\usepackage[colorlinks, allcolors=MidnightBlue]{hyperref}

\input{commands}

\begin{document}

\input{0abstract}
\maketitle
\setcounter{tocdepth}{2}
\tableofcontents
\input{1Introduction}
\input{2settingquadratic}
\input{3settingtrop}

%\input{4etic}
\input{5higherdimensions}

\input{6applications}

\bibliographystyle{plain}
\bibliography{enriched}
\end{document}

%% file: commands.tex
\usepackage{amssymb,amsmath,amsthm,amscd,mathtools,mathdots}
\usepackage{multirow}
\usepackage{enumitem}
\usepackage{caption}
\usepackage{subcaption}
\usepackage{tikz}
\usepackage{tikz-cd}
\usepackage{xcolor}
\usepackage[all,cmtip]{xy}

\newtheorem{thm}{Theorem}[section]
\newtheorem{lm}[thm]{Lemma}
\newtheorem{coro}[thm]{Corollary}
\newtheorem{prop}[thm]{Proposition}
\newtheorem{conj}[thm]{Conjecture}
\theoremstyle{definition}
\newtheorem{df}[thm]{Definition}
\newtheorem{ex}[thm]{Example}
\theoremstyle{plain}
\newtheorem{rmk}[thm]{Remark}

\newcommand{\R}{\mathbb{R}}
\newcommand{\C}{\mathbb{C}}
\newcommand{\Z}{\mathbb{Z}}
\newcommand{\A}{\mathbb{A}}

\newcommand{\Q}{\mathbb{Q}}
\newcommand{\GW}{\operatorname{GW}}
\newcommand{\W}{\operatorname{W}}
\newcommand{\Tr}{\operatorname{Tr}}
\newcommand{\Jac}{\operatorname{Jac}}
\newcommand{\glob}{\mathcal{O}}
\newcommand{\sign}{\operatorname{sign}}

\newcommand{\lra}{\longrightarrow}
\newcommand{\lmt}{\longmapsto}

\newcommand{\val}{\operatorname{val}}

\newcommand{\tC}{\widetilde{C}}
\newcommand{\mult}{\widetilde{\operatorname{mult}}}

\newcommand{\Card}{\operatorname{Card}}
\newcommand{\MV}{\operatorname{MVol}}
\newcommand{\Lop}{\Lambda^{\mathrm{odd}}}
\newcommand{\Ini}{\operatorname{In}}
\newcommand{\GWk}{\operatorname{GW}\left(k\right)}

\newcommand{\qinv}[1]{\left\langle#1\right\rangle}
\newcommand{\qqinv}[1]{\left(\left\langle#1\right\rangle\right)}
\newcommand{\Tra}[2][E/k]{\operatorname{Tr}_{#1}\left(#2\right)}

\newcommand{\tV}{\widetilde{V}}
\newcommand{\Puiseux}{\{\!\{t\}\!\}}
\newcommand{\NP}{\operatorname{NP}}
\newcommand{\cL}{\mathcal{L}}

\makeatletter
\let\save@mathaccent\mathaccent
\newcommand*\if@single[3]{%
    \setbox0\hbox{${\mathaccent"0362{#1}}^H$}%
    \setbox2\hbox{${\mathaccent"0362{\kern0pt#1}}^H$}%
    \ifdim\ht0=\ht2 #3\else #2\fi
}
\newcommand*\rel@kern[1]{\kern#1\dimexpr\macc@kerna}
\newcommand*\widebar[1]{\@ifnextchar^{{\wide@bar{#1}{0}}}{\wide@bar{#1}{1}}}
\newcommand*\wide@bar[2]{\if@single{#1}{\wide@bar@{#1}{#2}{1}}{\wide@bar@{#1}{#2}{2}}}
\newcommand*\wide@bar@[3]{%
    \begingroup
    \def\mathaccent##1##2{%
        \let\mathaccent\save@mathaccent
        \if#32 \let\macc@nucleus\first@char \fi
        \setbox\z@\hbox{$\macc@style{\macc@nucleus}_{}$}%
        \setbox\tw@\hbox{$\macc@style{\macc@nucleus}{}_{}$}%
        \dimen@\wd\tw@
        \advance\dimen@-\wd\z@
        \divide\dimen@ 3
        \@tempdima\wd\tw@
        \advance\@tempdima-\scriptspace
        \divide\@tempdima 10
        \advance\dimen@-\@tempdima
        \ifdim\dimen@>\z@ \dimen@0pt\fi
        \rel@kern{0.6}\kern-\dimen@
        \if#31
        \overline{\rel@kern{-0.6}\kern\dimen@\macc@nucleus\rel@kern{0.4}\kern\dimen@}%
        \advance\dimen@0.4\dimexpr\macc@kerna
        \let\final@kern#2%
        \ifdim\dimen@<\z@ \let\final@kern1\fi
        \if\final@kern1 \kern-\dimen@\fi
        \else
        \overline{\rel@kern{-0.6}\kern\dimen@#1}%
        \fi
    }%
    \macc@depth\@ne
    \let\math@bgroup\@empty \let\math@egroup\macc@set@skewchar
    \mathsurround\z@ \frozen@everymath{\mathgroup\macc@group\relax}%
    \macc@set@skewchar\relax
    \let\mathaccentV\macc@nested@a
    \if#31
    \macc@nested@a\relax111{#1}%
    \else
    \def\gobble@till@marker##1\endmarker{}%
    \futurelet\first@char\gobble@till@marker#1\endmarker
    \ifcat\noexpand\first@char A\else
    \def\first@char{}%
    \fi
    \macc@nested@a\relax111{\first@char}%
    \fi
    \endgroup
}
\makeatother

\makeatletter
\def\mathcenterto#1#2{\mathclap{\phantom{#1}\mathclap{#2}}\phantom{#1}}
\let\old@widetilde\widetilde
\def\widetildeto#1#2{\mathcenterto{#2}{\old@widetilde{\mathcenterto{#1}{#2\,}}}}
\let\old@widehat\widehat
\def\widehatto#1#2{\mathcenterto{#2}{\old@widehat{\mathcenterto{#1}{#2\,}}}}
\makeatother

%% file: 0abstract.tex
\begin{abstract}
Using tropical geometry one can translate problems in enumerative geometry to combinatorial problems. Thus tropical geometry is a powerful tool in enumerative geometry over the complex and real numbers. Results from $\A^1$-homotopy theory allow to enrich classical enumerative geometry questions and get answers over an arbitrary field. In the resulting area, $\A^1$-enumerative geometry, the answer to these questions lives in the Grothendieck-Witt ring of the base field $k$.
In this paper, we use tropical methods in this enriched set up by showing B\'ezout's theorem and a generalization, namely the Bernstein-Kushnirenko theorem, for tropical hypersurfaces enriched in $\GW(k)$.
\end{abstract}

%% file: 1Introduction.tex
\section{Introduction}
%Intro A1-enumerative geometry
%Enriched Bezout in GW(k)
%Intro tropical geometry
%Enriched patchworking
%Proof of tropical Bezout (classical result)
%Summary of results

Classically, B\'ezout's theorem states that any $n$ hypersurfaces in $\mathbb{P}^n_{\mathbb{C}}$ of degrees $d_1,\ldots,d_n$ that intersect transversally, intersect in $d_1\cdots d_n$ points. The count is invariant of the choice of hypersurfaces. This invariance breaks down if we replace the base field~$\C$ by a non-algebraically closed field $k$.
For example for $k=\R$ some of the intersections might only be defined over the complex numbers.
% add a picture with example
Motivated by results from $\A^1$-homotopy theory there is a new way of counting geometric objects when the base field $k$ is not algebraically closed. The resulting count is valued in the Grothendieck-Witt ring $\GW(k)$ of $k$. This way of counting restores the invariance in the ``relatively oriented" case.

Recall that
the Grothendieck-Witt ring $\GW(k)$ is generated by $\qinv{ a}$ with $a\in k^{\times}/(k^{\times})^2$ and that~$h\coloneqq \qinv{1}+\qinv{-1}$ denotes the \emph{hyperbolic form}.
In \cite{McKean} McKean proves B\'ezout's theorem enriched in $\GW(k)$: Let $V_1,\ldots,V_n$ be hypersurfaces in $\mathbb{P}^n_k$ defined by homogeneous polynomials $F_1,\ldots,F_n$ of degrees $d_1,\ldots,d_n$, respectively, such that $\sum_{i=1}^nd_i\equiv n+1\mod 2$. 
Assume that all the common zeros of $F_1,\ldots,F_n$ lie in $U_0=\{x_0\neq0\}\subset\mathbb{P}^n_k$ and set polynomials on the affine chart $f_i(x_1,\ldots,x_n)\coloneqq F_i(1,x_1,\ldots,x_n)$. Furthermore, assume that the hypersurfaces intersect transversally at every intersection point, then 
\begin{equation}
\label{eq:McKean}
\sum_{p\in V_1\cap\ldots\cap V_n}\Tr_{k(p)/k}\qinv{\det \operatorname{Jac}(f_1,\ldots,f_n)(p) }=\frac{d_1\cdot d_2\cdots d_n}{2}\cdot h\in \GW(k).
\end{equation}

In this paper we reprove \eqref{eq:McKean} and generalize the result using tropical geometry.
To do this, we define \emph{enriched tropical hypersurfaces}. The definition is motivated by Viro's patchworking.
\begin{df}
An \emph{enriched tropical hypersurface} $\tV=(V,(\alpha_I))$ in $\R^n$ is a tropical hypersurface~$V$ in $\R^n$ together with an element $\alpha_I\in k^{\times}/(k^{\times})^2$ assigned to each connected component of $\R^n\setminus V$.
\end{df}

Recall that connected components of $\R^n\setminus V$ correspond to vertices in the dual subdivision of~$V$. So we can alternatively assign coefficients $\alpha_v$ to each vertex $v$ in the dual subdivision of $V$.
An enriched tropical hypersurface can be seen as a homotopy equivalence class of \emph{Viro polynomials} associated to the dual subdivision of an embedded tropical hypersurface in $\R^n$ and a choice of coefficients in $k^{\times}$ for each vertex in the dual subdivision. These Viro polynomials are of the form
\[\sum_{I\in A}\alpha_Ix^It^{\varphi(I)}\]
where $A$ is a finite subset of $\Z^n$, $\alpha_I\in k^\times$ for any $I\in A$, $x=(x_1,\ldots,x_n)$ and $\varphi:A\longrightarrow \mathbb{Q}$ is a convex function that assigns a rational number to each exponent $I\in A$.
These Viro polynomials can be viewed as polynomials over the field of \emph{Puiseux series} over~$k$ 
\[k\Puiseux=\left\{\sum_{i=i_0}^\infty c_it^{\frac{i}{n}}\mid c_i\in k,i_0\in \Z,n\in \mathbb{N}\right\}.\]

Let $\tV_1,\ldots,\tV_n$ be $n$ enriched tropical hypersurfaces in $\R^n$ with Viro polynomials \linebreak $f_1,\ldots,f_n$ in~$k\Puiseux[x_1,\ldots,x_n]$. Assume that $k$ is a field of characteristic $0$ or characteristic bigger than the diameter of the Newton polygons of the $f_i$.
Inspired by \eqref{eq:McKean} we define the \emph{enriched intersection multiplicity} of $\tV_1,\ldots,\tV_n$ at an intersection point $p$ to be
\begin{multline}
\label{eq:enrichedIMDef}
\mult_p(\tV_1,\ldots,\tV_n)\coloneqq\\
\sum_{z} \Tr_{\kappa(z)/k\Puiseux}\qinv{ \det \operatorname{Jac}(f_1,\ldots,f_n)(z)}\in\GW(k\Puiseux)\cong \GW(k)
\end{multline}
%{\color{blue} are we sure it's $E_t$ and not $E$. Why the $t$?}
where $\kappa(z)$ is the residue field of $z$ and $z$ ranges over all the common zeros of $f_1,\ldots,f_n$ that tropicalize to $p$. Note that taking the rank of \eqref{eq:enrichedIMDef} recovers the classical (tropical) intersection multiplicity.

It is rather tedious to compute this enriched intersection multiplicity. The main theorem of this paper finds a purely combinatorial rule to determine it.
\begin{df}
Let $\Lambda^{\text{odd}}$ be the subset of $\Z^n$ consisting of tuples $(a_1,\ldots,a_n)$ with $a_i\equiv1\mod2$ for $i=1,\ldots n$. We call the elements of $\Lop$ \emph{odd} lattice points.
\end{df}
Our Main Theorem states that the enriched intersection multiplicity is determined by the coefficients of the \emph{odd} vertices of the dual subdivision.
To avoid confusion, we say that the vertices of a polytope that belong to a minimal generating set with respect to the convex hull are its \emph{corner vertices}.
\begin{thm}[Main Theorem]
\label{thm:Main}
Let $P$ be the parallelepiped in the dual subdivision of the  union $\tV_1\cup\tV_2\cup\ldots\cup\tV_n$ corresponding to the intersection point $p$.
Assume that the volume of $P$ equals~$m$. 
Let $r$ be the number of corner vertices of $P$ that are elements of $\Lambda^{\text{odd}}$.
Then
\[\mult_p(V_1,\ldots,V_n)=\sum_{v\text{ an odd corner of }P} \qinv{ \epsilon(v) \alpha_v} +\frac{m-r}{2}\cdot h\in\GW(k).\]
Here, $\alpha_v$ denotes the coefficient of the vertex $v$ and $\epsilon(v)$ is a sign determined by the intersection.
\end{thm}

There is an easy proof for B\'ezout's theorem for tropical curves (not enriched in $\GW(k)$) in~\cite{Sturmfels}. Namely, the intersection points of two tropical curves $C_1$ and $C_2$ with Newton polygons $\Delta_{d_1}$ and $\Delta_{d_2}$ (see \eqref{eq:deltad} for the definition of $\Delta_d$), respectively, correspond to parallelograms in the dual subdivision of $C_1\cup C_2$. The area of such a parallelogram is equal to the intersection multiplicity at the corresponding intersection point and the rest of the dual subdivision of $C_1\cup C_2$ consists of the dual subdivisions of $C_1$ and~$C_2$. Thus the number of intersections counted with multiplicities is
\[
\operatorname{Area}(\Delta_{d_1+d_2})-\operatorname{Area}(\Delta_{d_1})-\operatorname{Area}(\Delta_{d_2})=\frac{(d_1+d_2)^2}{2}-\frac{d_1^2}{2}-\frac{d_2^2}{2}=d_1\cdot d_2.
\]

This is a particular instance of a more general statement. 
Two tropical curves $C_1$ and $C_2$ in~$\R^2$ of with Newton polygons~$\Delta_1$ and~$\Delta_2$, respectively, intersect in a number of points counted with multiplicities equal to the \emph{mixed volume}
\begin{equation}
\label{eq:MixedArea}
\MV(\Delta_1,\Delta_2)=\operatorname{Area}(\Delta_{1}+\Delta_{2})-\operatorname{Area}(\Delta_{1})-\operatorname{Area}(\Delta_{2}),
\end{equation}
where $\Delta_1+\Delta_2$ is the Minkowski sum of the polygons $\Delta_1$ and $\Delta_2$.

%to the enriched setting by assigning intersection multiplicities valued in $\GW(k)$ to each new parallelogram in the dual subdivision of $C_1\cup C_2$.
Two tropical curves embedded in $\R^2$ intersect \emph{tropically transversely} if they intersect in a finite number of points and every point of the intersection belongs to an edge in each of the curves. 
A direct consequence of our Main Theorem \ref{thm:Main} is that we can quadratically enrich the above argument for enriched tropical curves that intersect tropically transversally at every intersection point.
Theorem \ref{thm:Main} implies that the only non-hyperbolic contribution to \eqref{eq:McKean} comes from odd points on the boundary of $\Delta_1+\Delta_2$ (see Proposition \ref{prop:number of parallelepipeds} and Corollary \ref{coro:EnrichedTropicalBezout}). Hence, we get the following for which we need the characteristic of $k$ to be $0$ or bigger than the maximum of the diameters of $\Delta_1$ and $\Delta_2$.
\begin{thm}
Let $\tC_1$ and $\tC_2$ be two enriched tropical curves with Newton polygons $\Delta_1$ and $\Delta_2$, respectively. Assume that they intersect tropically transversally at every intersection point and  that~$\partial(\Delta_1+\Delta_2)\cap \Lop=\emptyset$.
Then 
\[\sum_{p\in \tC_1\cap \tC_2}\mult_p(\tC_1,\tC_2)=\frac{\MV(\Delta_1,\Delta_2)}{2}\cdot h\in \GW(k). \]
\end{thm}
In particular, we recover \eqref{eq:McKean} for curves, since the condition $d_1+d_2$ being odd holds if and only if  $\partial(\Delta_{d_1}+\Delta_{d_2})\cap \Lop=\partial(\Delta_{d_1+d_2})\cap \Lop=\emptyset$.
%\begin{thm}
%\label{thm:dim2}
%Let $\tC_1$ and $\tC_2$ be enriched tropical curves over $k$ of degree $\Delta_{1}$ and $\Delta_{2}$, respectively. If $\tC_1$ and $\tC_2$ intersect tropically transversely,  then 
%\[\sum_{p\in C_1\cap C_2}\widetilde{\operatorname{mult}}_p(\tC_1,\tC_2)=
%\left(\left\lfloor\frac{\MV(\Delta_1,\Delta_2)}{2}\right\rfloor-\lambda\right)\left(\qinv{1}+\qinv{-1}\right)+
%\left\qinv{ a_1,\dots,a_{\lambda}\right} \in\W(k),\]
%where $\lambda=\min\left(\Card\left(\partial(\Delta_1+\Delta_2)\cap\Lop\right),\MV(\Delta_1,\Delta_2)\right)$.
%\end{thm}

\begin{coro}
Let $\tC_1$ and $\tC_2$ be enriched tropical curves over $k$ of with Newton polygons $\Delta_{d_1}$ and $\Delta_{d_2}$, respectively. If $d_1+d_2\equiv 1\mod 2$, then 
\[\sum_{p\in C_1\cap C_2}\widetilde{\operatorname{mult}}_p(\tC_1,\tC_2)=\frac{d_1\cdot d_2}{2}\;h\in \GW(k). \]
\end{coro}

When $d_1+d_2\equiv0\mod2$, we are dealing with the \emph{non-relatively orientable case}, in which the left hand side in \eqref{eq:McKean} depends on choice of coefficients.
However, the left hand side of \eqref{eq:McKean} can not equal any arbitrary element of $\GW(k)$. More precisely, we have the following possibilities for the intersection of two enriched tropical curves. 
\begin{coro}
Let $\tC_1$ and $\tC_2$ be enriched tropical curves over $k$ with Newton polygons $\Delta_{d_1}$ and~$\Delta_{d_2}$, respectively. If $\tC_1$ and $\tC_2$ intersect tropically transversally and $d_1+d_2\equiv 0\mod 2$, then 
\[\sum_{p\in C_1\cap C_2}\widetilde{\operatorname{mult}}_p(\tC_1,\tC_2)=
\frac{d_1\cdot d_2-\min(d_1,d_2)}{2}\;h+
\qinv{ a_1,\dots,a_{\min(d_1,d_2)}}\in \GW(k),\]
where $a_1,\ldots,a_{\min(d_1,d_2)}$ can be any element in $k^{\times}/(k^{\times})^2$.
\end{coro}

This also works in higher dimensions (see \cite{HuberSturmfels}): Let $V_1,\ldots,V_n$ be tropical hypersurfaces in $\R^n$ with Newton polytopes $\Delta_1,\ldots,\Delta_n$. Then the number of intersection points counted with multiplicities equals the \emph{mixed volume}
\begin{equation}
\label{eq:MixedVolume}
\MV(\Delta_1,\ldots,\Delta_n)\coloneqq \text{coefficient of $\lambda_1\cdots\lambda_n$ in $%R(\lambda_1,\ldots,\lambda_n)\coloneqq
\operatorname{vol}(\lambda_1\Delta_1+\ldots+\lambda_n\Delta_n)$.} \end{equation}
This agrees with \eqref{eq:MixedArea} in dimension $2$.

We replace the relative orientation condition by a purely combinatorial condition, namely we assume that 
\[\partial(\Delta_1+\ldots+\Delta_n)\cap \Lop=\emptyset.\]
In this case our Main Theorem \ref{thm:Main} implies the following corollary. For this we assume the characteristic of $k$ to be $0$ or bigger than the maximum of the diameters of the $\Delta_i$.
\begin{coro}
Let $\tV_1,\ldots,\tV_n$ be $n$ tropical hypersurfaces with Newton polytopes \linebreak$\Delta_1,\ldots,\Delta_n$, respectively, such that $\partial(\Delta_1+\ldots+\Delta_n)\cap \Lop=\emptyset$.
Assume that $\tV_1,\ldots,\tV_n$ intersect tropically transversally at every point.
Then 
\[\sum_{p\in \tV_1\cap\ldots \cap \tV_n}\mult_p(\tV_1,\ldots,\tV_n)=\frac{\MV(\Delta_1,\ldots,\Delta_n)}{2}\;h\in \GW(k).\]
\end{coro}

From this we can derive a quadratic enrichment of a Theorem of Bernstein and Kushnirenko. McKean's B\'ezout theorem \eqref{eq:McKean} is a special case of this. %{\color{red} do we need to add a characteristic assumption here?}
\begin{coro}[Enriched Bernstein-Kushnirenko Theorem]
Let $f_1,\ldots,f_n$ be Laurent polynomials in $n$ variables with Newton polytopes $\Delta_1,\ldots,\Delta_n$, respectively. %Assume $k$ is a field of characteristic $0$ or bigger than $\max_i \{\operatorname{diam}\Delta_{i}\}$.

If 
${\partial(\Delta_1+\ldots+\Delta_n)\cap \Lop=\emptyset}$, then
\[\sum_z \Tr_{k(z)/k}\qinv{ \det \operatorname{Jac}(f_1,\ldots,f_n)(z)}=\frac{\MV(\Delta_1,\ldots,\Delta_n)}{2}\;h\in \GW(k).\]
Here, the sum runs over all solutions $z$ to $f_1=\ldots=f_n=0$ in $\operatorname{Spec}k[x_1^{\pm1},\ldots,x_n^{\pm1}]$.
\end{coro}

%In case 
%\[\partial(\Delta_1+\ldots+\Delta_n)\cap \Lop\not=\emptyset\]
%we can still say something.

%\begin{coro}
%Let $\tV_1,\ldots,\tV_n$ be $n$ tropical hypersurfaces with Newton polytopes $\Delta_1,\ldots,\Delta_n$ that intersect tropically transversally at every point.
%Then 
%\[\sum_{p\in \tV_1\cap\ldots \cap \tV_n}=\qinv{ a_1,\ldots,a_{\lambda}}\in \W(k)\]
%where 
%\[\lambda = \min(\Card(\partial(\Delta_1+\ldots+\Delta_n)),\MV(\Delta_1,\ldots,\Delta_n)).\]
%\end{coro}
%{\color{red} I think we have to reformulate theorem 1.2 and Corollary 1.6: Example if $\lambda$ is even, but the mixed volume is odd, this doesn't work, we need the difference to be even!}

\subsection{Related work and the a general strategy for quadratically enriched tropical counts}

In \cite{MPS} Markwig, Payne and Shaw use tropical methods to redo the quadratically enriched count of bitangents to a quartic curve by Larson-Vogt \cite{LV}. Their strategy is similar to ours, they also use \emph{enriched tropical curves}: In many cases, questions in enumerative geometry can be solved by counting zeros of a general section of some vector bundle. This is the case for B\'ezout's theorem and the count of bitangents to a quartic curve. Hence, these counts equal the degree of the Euler class of the respective vector bundles and their quadratic enrichments are the degree of the ``$\mathbb{A}^1$-Euler class" of the respective vector bundles. There is a ``Poincar\'e-Hopf Theorem" (see Theorem \ref{thm:PoincareHopf} for the classical Poincar\'e-Hopf Theorem) for the degree of the $\A^1$-Euler class of a vector bundle, that is, the degree of the $\A^1$-Euler class equals the sum of ``local indices" at the zeros of a chosen section. One can write down explicit polynomials to compute these local indices and then interpret them tropically. In our case, the local index is what we call the \emph{enriched tropical multiplicity} defined in \eqref{eq:enrichedIMDef} and its tropical interpretation is our Main Theorem \ref{thm:Main}. 
Note that there are now also expository lecture notes on this approach and the main result of this paper for the case of curves \cite{pauli2024motivic}.

A major breakthrough in the use of tropical geometry in enumerative geometry was Mikhalkin's correspondence theorem \cite{MikhalkinCorrespondence}, which translates the question of counting plane algebraic curves to counting their tropical counterparts.
The authors have proved a quadratic enrichment of Mikhalkin's correspondence theorem \cite{JPPcorrespondence} for the counting of plane algebraic curves passing through a configuration of $k$ points, yielding quadratically enriched tropical curve counting invariants, and studied its arithmetic implications in \cite{JPMPRarithmetic} together with Markwig and R\"ohrle. There is work in progress \cite{JPMPRworkinprogress} to extend the quadratically enriched tropical correspondence theorem \cite{JPPcorrespondence} to curve counts for point conditions consisting not only of $k$-points but also of points defined over quadratic field extensions of $k$. A tropical wall crossing formula was derived for these counts in \cite{puentes2024wall}.

There is another refinement in tropical geometry that specializes to complex and real curve counting invariants. Block-Goettsche invariants are also curves counting invariants that interpolate between the complex and real curve counts \cite{BlockGoettsche}.
Nicaise, Payne and Schroeter  proposed a geometric interpretation of Block and G\"ottsche's refined tropical curve counting invariants in terms of virtual $\chi_{-y}$-specializations of motivic measures of semialgebraic sets in relative Hilbert schemes \cite{NicaisePayneSchroeter}.

%Here is the general recipe for using tropical methods in $\A^1$-enumerative geometry:
%\begin{enumerate}
   % \item Find the vector bundle such that counting the zeros of a general section solves your enumerative problem.
%    \item Check whether this vector bundle is relatively orientable (see Definition \ref{df:relativeorientation}).
 %   \item Find a polynomial formula for the local indices at the zeros by choosing local coordinates and a trivialization of the vector bundle around your zeros compatible with a relative orientation (see Definition \ref{df:local index}).
  %  \item Take a toric deformation (see subsection \ref{subsection:ToricDeformations}) of these polynomials and tropicalize and remember the coefficients to get a purely combinatorial rule and solve your problem using combinatorics.
%\end{enumerate}
%\begin{rmk}
%One can skip step $2$ in the recipe above.
%For example the vector bundle for the count of bitangents to a quartic curve is not relatively orientable. We also deal with the non-relatively orientable case in this paper. The answer to the enumerative problem will not necessarily be invariant in this case but our tropical methods still gives a list of possible answers which was not available earlier.
%\end{rmk}
\subsection{Outline} We start by recalling the prerequisites from $\A^1$-enumerative geometry in section \ref{section:A1enumerativegeometry} and from tropical geometry in section \ref{section:IntroTropicalGeometry}.
In section \ref{section:Higher dimensions} we prove our main theorem on the enriched tropical intersection multiplicity \ref{thm:mainthm}. In section \ref{section:applications} we use our main theorem to prove the enriched tropical B\'ezout theorem and the enriched Bernstein-Kushnirenko theorem.

\subsection{Acknowledgements}
Both authors have been supported by the ERC programme QUADAG.  This paper is part of a project that has received funding from the European Research Council (ERC) under the European Union's Horizon 2020 research and innovation programme (grant agreement No. 832833).\\ 
\includegraphics[scale=0.08]{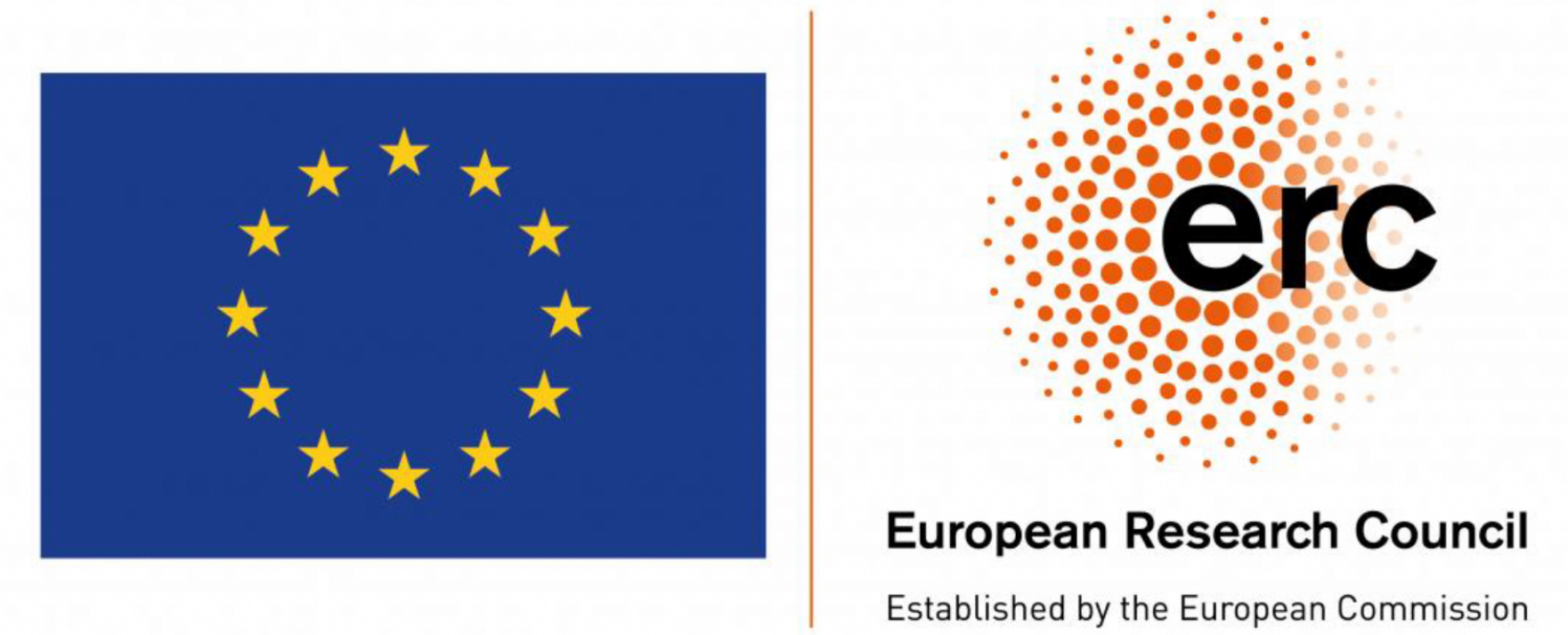}

The second author acknowledges the support of the Centre for Advanced Study at the Norwegian Academy of Science and Letters in Oslo, Norway, which funded and hosted the research project ``Motivic Geometry" during the 2020/21 academic year.

The second author would like to thank the Isaac Newton Institute for Mathematical Sciences, Cambridge, for support and hospitality during the programme ``K-theory, algebraic cycles and motivic homotopy theory" where work on this paper was undertaken. This work was supported by EPSRC grant no EP/R014604/1.

We would like to thank Marc Levine, Stephen McKean, Kris Shaw and Helge Ruddat for fruitful discussions. 

%% file: 2settingquadratic.tex
\section{Introduction to $\A^1$-enumerative geometry}
\label{section:A1enumerativegeometry}
In $\A^1$-enumerative geometry one uses machinery from $\A^1$-homotopy theory to enrich classical results from enumerative geometry yielding results over an arbitrary field $k$.
In this section, we introduce one way of doing this following the work of Jesse Kass and Kirsten Wickelgren in  \cite{KW-Euler,KW-EKL}. 

\subsection{The Grothendieck-Witt ring}
The enriched enumerative results will be valued in the Grothendieck-Witt ring. We recall the definition of $\GW(R)$ where $R$ is a commutative ring with unity.
For this, consider the set 
\[S\coloneqq \{\textrm{isometry classes of non-degenerate symmetric bilinear forms over }R\}.\]
On this set we have two binary operations. If 
\[b_1\colon V_1\times V_1\longrightarrow R\text{ and }b_2\colon V_2\times V_2\longrightarrow R\] are two non-degenerate symmetric bilinear forms over $R$, then the direct sum, as well as the tensor product,
\[b_1\oplus b_2\colon (V_1\oplus V_2)\times (V_1\oplus V_2)\longrightarrow R, \quad b_1\otimes b_2\colon (V_1\otimes V_2)\times (V_1\otimes V_2)\longrightarrow R,\]
are non-degenerate symmetric bilinear forms over $R$. One can check that this respects isometry classes, and thus, the set $S$ together with the operations $\oplus$ and $\otimes$ is a semi-ring $(S,\oplus,\otimes)$.
Group completion with respect to the direct sum yields a ring.
\begin{df}
\label{df:GW}
The Grothendieck-Witt ring $\GW(R)$ of a ring $R$ is the group completion of the semi-ring $(S,\oplus,\otimes)$ of isometry classes of non-degenerate symmetric bilinear forms over~$R$ with respect to taking the direct sum $\oplus$. 
\end{df}

We are mainly interested in the case when $R$ is a field $k$ of characteristic not equal to~$2$ in which case we have a nice presentation of $\GW(k)$:
Let $k^{\times}=k\setminus\{0\}$ be the set of units in $k$. If $\operatorname{char}k\neq 2$, any form can be diagonalized, i.e., for any symmetric bilinear form $\beta:V\times V\longrightarrow k$. Hence, we can find a basis for the $k$-vector space $V$, such that \begin{equation}
\label{eq:diagonal}
\beta((x_1,\ldots,x_n),(y_1,\ldots,y_n))=a_1x_1y_1+\ldots+a_nx_ny_n\end{equation}
for some $a_1,\ldots,a_n\in k^\times$ in this basis. Also note that if we replace one of the $a_i$ by $a_ib^2$ for some~$b\in k^{\times}$, the resulting form is in the same isometry class as $\beta$. 
Thus the form $\beta$ above \eqref{eq:diagonal} can be expressed as the direct sum of $n$ symmetric bilinear forms on a one-dimensional $k$-vector space. Indeed 
$\GW(k)$ is generated by the classes of bilinear forms 
\[\qinv{ a} \colon k\times k\longrightarrow k\text{, } (x,y)\mapsto axy\]
for $a\in k^{\times }/(k^{\times})^2$ (since classes in $\GW(k)$ are non-degenerate, we need $a\neq 0$) subject to the following two relations
\begin{enumerate}
    \item $\qinv{ a} \qinv{ b}=\qinv{ ab} $ for $a,b\in k^{\times}$
    \item $\qinv{ a}+\qinv{ b}=\qinv{ a+b}+\qinv{ ab(a+b)}$ for $a,b,a+b\in k^{\times}$.
\end{enumerate}
We use the notation $\qinv{ a_1,\ldots,a_s}\coloneqq \qinv{ a_1}+\ldots +\qinv{ a_s} \in \GW(k)$.
\begin{df}
We write $h$ for the \emph{hyperbolic form}, that is the form on a $2$-dimensional $k$-vector space (or free rank two R-module over $R$ when $R$ is not a field), with Gram matrix \[\begin{pmatrix}0&1\\1&0\end{pmatrix}.\]
\end{df}
\begin{rmk}
\label{rmk:hyperbolic}
\begin{enumerate}
\item Assume $a$ is a unit in $R$, then the class of the symmetric bilinear form on a rank $2$ $R$-module defined by $\begin{pmatrix}0&a\\a&0\end{pmatrix}$ equals $h$.
\item If $k$ is a field of characteristic not equal to $2$, then after diagonalizing, we get that the hyperbolic form equals  $$h=\qinv{1} +\qinv{-1}.$$
Furthermore, one can deduce from relation $2$ above that for $a\in k^{\times}$ the equality $$\qinv{ a}+\qinv{ -a}=\qinv{ 1}+\qinv{-1} =h$$ holds in $\GW(k)$.
\end{enumerate}
\end{rmk}
\begin{df}
We say that the \emph{rank} of a symmetric bilinear form
\[\beta:V\times V\longrightarrow R\]
is the rank of the $R$-module $V$.
\end{df}
Taking the rank extends to a homomorphism  
\[\operatorname{rank}:\GW(R)\longrightarrow \Z.\]
\begin{ex}
\label{ex:GW(C)}
Let $k=\C$. Since $\C$ is algebraically closed, there is only one generator $\qinv{ 1} \in \GW(k)$ and thus $\GW(k)\cong \Z$ where the isomorphism is the rank homomorphism. In particular, results in classical enumerative geometry coincide with the counts enriched in $\GW(k)$ for an algebraically closed field $k$ by taking the rank.
\end{ex}

\begin{ex}
For $k=\mathbb{R}$, $\GW(k)$ has two generators, namely $\qinv{1}$ and $\qinv{ -1}$. In fact, an element in $\GW(\mathbb{R})$ is completely determined by its rank and its signature.
\end{ex}

\begin{ex}
\label{ex:GW of Puiseux}
The Grothendieck-Witt ring of the field $k\Puiseux$ of Puiseux series over $k$ is isomorphic to $\GW(k)$. This is because $k\Puiseux^{\times}/(k\Puiseux^{\times})^2\cong k^{\times}/(k^{\times})^2$.
More precisely when $a_m\in k^{\times}$,
\begin{align*}
\sum_{i=m}^{\infty}a_it^{\frac{i}{n}}=a_m\cdot\left(t^{\frac{m}{n}}+\sum_{i=m+1}^{\infty}\frac{a_i}{a_m}t^{\frac{i}{n}}\right)
\end{align*}
and $\left(t^{\frac{m}{n}}+\sum_{i=m+1}^{\infty}\frac{a_i}{a_m}t^{\frac{i}{n}}\right)$ is a square in $k\Puiseux$. So
$\sum_{i=m}^{\infty}a_it^{\frac{i}{n}}\mapsto a_{m}$ defines an isomorphism
$k\Puiseux^{\times}/(k\Puiseux^{\times})^2\cong k^\times/(k^\times)^2$.

It follows that $\GW(k)$ and $\GW(k\Puiseux)$ have the same generators and
Markwig-Payne-Shaw show that also the relations in the Grothendieck-Witt ring are respected by this isomorphism \cite[Theorem 4.7]{MPS}.
Hence,
we get the following isomorphism 
\[ \Ini\colon
\begin{array}{ccc}
\GW\left(k\Puiseux\right) &   \lra   &  \GWk\\
\qinv{\sum_{i=m}^\infty a_it^\frac{1}{n}     }      &   \lmt   &    \qinv{ a_m}
\end{array}
\]
sending a generator to its initial.
\end{ex}

\subsubsection{The Witt ring}
A non-degenerate symmetric bilinear form $\beta:V\times V\longrightarrow R$ over a ring $R$ is \emph{split} if there exists a submodule $N\subset V$ such that $N$ is a direct summand of $V$ and $N$ is equal to its orthogonal complement $N^{\bot}$.
We say that two non-degenerate symmetric bilinear forms $\beta\colon V\times V\longrightarrow R$ and~${\beta'\colon V'\times V'\longrightarrow R}$ are \emph{stably equivalent} if there exist split symmetric bilinear forms $s$ and $s'$ such $\beta\oplus s\cong\beta'\oplus s'$.
\begin{df}
The \emph{Witt ring} of $R$ is the set of classes of stably equivalent non-degener\-ate symmetric bilinear forms with addition given by the direct sum $\oplus$ and multiplication by the tensor product $\otimes$.
\end{df}

\begin{rmk}
\label{rmk:GWandW}
If $R$ is a field $k$ of characteristic not equal to $2$, then 
the split non-degenerate symmetric bilinear forms are exactly the multiples of the hyperbolic form $h$.
Recall that in this case we have that $\langle a\rangle +\langle-a\rangle=h$ for any unit $a$ and hence 
\[\operatorname{W}(k)=\frac{\GW(k)}{\Z\cdot h}.\]
More generally, if $R$ is local and $2$ is invertible then an element of $\GW(R)$ is completely determined by its rank and the associated element of $\W(R)$.
\end{rmk}

%Assume $k$ is a field of characteristic not equal to $2$.
%By Remark \ref{rmk:hyperbolic} the ideal in $\GW(k)$ generated by $h$ is equal to $\Z\cdot h$. 

\begin{ex}
The Witt ring of $\C$ is isomorphic to $\Z/2\Z$.
\end{ex}

%\begin{rmk}
%\label{rmk:GWandW}
%An element $\beta\in \GW(k)$ is completely determined by its rank and its image in $\W(k)$.
%\end{rmk}

\subsubsection{Trace \label{ss:trace}}
Assume that $R$ is a commutative ring. We are particularly interested in the case that $R$ is a finite \'etale $k$-algebra.
For a finite projective $R$-algebra $L$ one can define the trace $\Tr_{L/R}:L\longrightarrow R$ that sends $b\in L$ to the trace of the multiplication map $m_b(x)=b\cdot x$.
If $L$ is \'etale over $R$ this induces the \emph{trace map}
$\Tr_{L/R}\colon \GW(L)\longrightarrow \GW(R)$ which sends the class of a bilinear form $\beta\colon V\times V\longrightarrow L$ over $L$ to the form
\[V\times V\xrightarrow{\beta}L\xrightarrow{\operatorname{Tr}_{L/R}}R\]
over $R$.

We will compute several trace forms in our main result. So we already collect some facts and computations about the trace form here.
Let $E$ be a finite \'{e}tale $R$-algebra.  
\begin{enumerate}
    \item If $R=k$ is a field, then $E=L_1\times \ldots\times L_s$ for some finite separable field extensions $L_1,\ldots,L_s$ of $k$ and the trace map $\Tr_{E/k}:E\longrightarrow k$ equal to the sum of  field traces $\Tr_{E/k}=\sum_{i=1}^s\Tr_{L_i/k}$.
    \item $\Tr_{E/R}$ is $R$-linear.
    \item Let $F$ be a finite \'{e}tale $E$-algebra. Then 
    \[\Tr_{F/R}=\Tr_{E/R}\circ\Tr_{F/E}.\]
\end{enumerate}

%Let $\beta:V\times V\longrightarrow E$ a non-degenerate symmetric bilinear form on $E$. This defines a non-degenerate symmetric bilinear form on $R$, the \emph{trace form}
%\[V\times V\xrightarrow{\beta}E\xrightarrow{\Tr_{E/R}}R.\]
%One can define $\GW(E)$ in the same way as in Definition \ref{df:GW} and we get a map
%\[\Tr_{E/R}:\GW(E)\longrightarrow \GW(R).\]

From now on let $k$ be a field.
\begin{lm}
\label{lm:trace}
Let $L$ be a finite \'{e}tale $k$-algebra.
Let $E=\frac{L[x]}{(x^m-D)}$ for some $D\in L^\times$, and assume that $E$ is \'{e}tale over $L$. Then $\Tr_{E/L}(1)=m$ and $\Tr_{E/L}(x^s)=0$ for $s=1,\ldots,m-1$.
\end{lm}
\begin{proof}
We have the following $L$-basis for $E$: $1,x,x^2,\ldots,x^{m-1}$. Recall that $\Tr_{E/L}(a)$ for $a\in E$ is the trace of the $L$-linear map $m_a:E\longrightarrow E$ defined by $m_a(y)=a\cdot y$. So we are looking for the matrix of $m_a$ with respect to the basis $1,x,\ldots,x^{m-1}$. %Write $y=y_1\cdot 1+y_2\cdot x+\ldots+y_m\cdot x^{d-1}$.
\begin{enumerate}
    \item If $a=1$, then this matrix is the identity matrix and its trace equals the $L$-dimension of $E$, namely $m$.
    \item If $a=x^{s}$ for some $s\in \{1,\ldots,m-1\}$, then every entry of the diagonal of this matrix equals $0$.
\end{enumerate}\vspace{-20pt}
\end{proof}
Since the trace $\Tr_{E/L}$ is $L$-linear, Lemma \ref{lm:trace} tells us what $\Tr_{E/L}(a)$ is for any element $a\in E$ in case $E=\frac{L[x]}{(x^m-D)}$.

\begin{lm}
Let $E$ be a finite \'{e}tale $R$-algebra of rank $m$.
Then $\Tr_{E/R}(h)=m\cdot h$.
\end{lm}
\begin{proof}
This follows directly from the fact that the hyperbolic form $h$ is split.
%\[\Tr_{E/L}(\qinv{ 1}+\qinv{-1})=\Tr_{E/k}(\qinv{1})+\Tr_{E/k}(\qinv{-1}).\]
%Assume $\Tr_{E/L}(\qinv{ 1})$ is represented by a form $\beta$, then $\Tr_{E/L}(\qinv{-1})$ is represented by $-\beta$.
\end{proof}

%\begin{lm}
%\label{lm:hyperbolic form}
%Let $L$ be a finite \'{e}tale $k$-algebra.
%The non-degenerate symmetric bilinear form  which has Gram matrix $\begin{pmatrix}0 &a\\a&0\end{pmatrix}$ where $a\in L^{\times}$, is equal to the hyperbolic form $h$ in $\GW(L)$. {\color{red} characteristic 2!!!}
%\end{lm}
%\begin{proof}
%We perform the following basis change on quadratic forms
%\[\begin{pmatrix}\frac{1}{2}& \frac{1}{2}\\-\frac{1}{a}&\frac{1}{2}\end{pmatrix}\cdot \begin{pmatrix}0& a\\a&0\end{pmatrix}\cdot \begin{pmatrix}\frac{1}{2}& -\frac{1}{a}\\\frac{1}{2}&\frac{1}{2}\end{pmatrix}=\begin{pmatrix}1& 0\\0&-1\end{pmatrix}\]
%\[\begin{pmatrix}1&1\\1&-1\end{pmatrix}\cdot \begin{pmatrix}0& a\\a&0\end{pmatrix}\cdot \begin{pmatrix}1&1\\1&-1\end{pmatrix}=\begin{pmatrix}2a& 0\\0&-2a\end{pmatrix}\]
%and get the Gram matrix of a quadratic form equivalent to $h$ in $\GW(k)$. 
%\end{proof}

\begin{prop}
\label{prop:traces}
Let $L$ be a finite \'{e}tale $k$-algebra and let $E=\frac{L[x]}{(x^m-D)}$, for some $D\in L^\times$. Further, assume that $\operatorname{char}k$ does not divide $m$. Then for $a\in L$ we get
\begin{enumerate}
    \item $\Tr_{E/L}(\qinv{ m\cdot a})=\begin{cases}
\qinv{ a} +\frac{m-1}{2}\;h & m\text{ odd}, \\
\qinv{ a}+\qinv{  a\cdot D}+\frac{m-2}{2}\;h &  m \text{ even.}
\end{cases}$
    \item $\Tr_{E/L}(\qinv{ m\cdot a\cdot x})=\begin{cases}
\qinv{ a\cdot D} +\frac{m-1}{2}\;h & m\text{ odd}, \\
\frac{m}{2}\;h &  m \text{ even.}
\end{cases}$
\end{enumerate}

\end{prop}
\begin{proof}

We have the following $L$-basis for $E$: $1,x,\ldots,x^{m-1}$. 
Let $M=(M_{ij})$ be the Gram matrix of $\Tr_{E/L}(\qinv{ m\cdot a})$.
Then the $(i,j)$th entry $M_{ij}$ of $M$ equals
\[\Tr_{E/L}(m\cdot a\cdot b_i\cdot b_j)=\Tr_{E/L}(m\cdot a\cdot x^{i-1}\cdot x^{j-1}),\] where $b_i=x^{i-1}$ is the $i$th basis element of the chosen $L$-basis of $E$.
In particalur, we have 
\[M_{ij}=\begin{cases}m^2\cdot a & \text{if } i=j=1\\
m^2\cdot a\cdot D & \text{if } i+j=m+1 \\
0 & \text{otherwise}\end{cases}\]
by Lemma \ref{lm:trace}
and thus the Gram matrix of $\Tr_{E/L}(\qinv{ m\cdot a})$ looks like
\begin{align*}%&\begin{pmatrix}\Tr(m\cdot a)& \Tr(m\cdot a\cdot x)&\dots&\dots& \dots & \Tr(m\cdot a\cdot x^{m-1})\\
%\Tr(m\cdot a\cdot x)& \ddots&\ddots&\ddots&  \Tr(m\cdot a\cdot x^{m-1}) & \Tr(m\cdot a\cdot D)\\ 
%\vdots & \iddots &\iddots& \iddots& \Tr(m\cdot a\cdot D) & \vdots\\
%\vdots & \iddots &\iddots& \iddots& \vdots & \vdots\\
%\vdots & \Tr(m\cdot a\cdot x^{m-1}) &\Tr(m\cdot a\cdot D) & \ldots& \ldots  & \vdots\\
%\Tr(m\cdot a\cdot x^{m-1})& \Tr(m\cdot a\cdot D) & \ldots&\ldots &\ldots & \Tr(a\cdot m\cdot a\cdot x^{2m-2})\end{pmatrix}\\
%\overset{\ref{lm:trace}}{=}&
M=\begin{pmatrix}m^2\cdot a& 0&\ldots&\ldots& \ldots & 0\\
0& \ldots&\ldots&\ldots&  0 & m^2\cdot a\cdot D\\ \ldots & \ldots &\ldots& \ldots& m^2\cdot a\cdot D & 0\\
\ldots & \ldots &\ldots& \iddots& \ldots & \ldots\\
\ldots & 0 &m^2\cdot a\cdot D & \ldots& \ldots  & \ldots\\
0& m^2\cdot a\cdot D & 0&\ldots &\ldots & 0\end{pmatrix}.
\end{align*}
By Remark \ref{rmk:hyperbolic}, this is equivalent to
%$\begin{cases}\qinv{ a} +\frac{m-1}{2}\;h & m\text{ odd} \\ \qinv{  a}+\qinv{  a\cdot D}+\frac{m-2}{2}\;h &  m \text{ even} \end{cases}$ in $\GW(L)$.
$\qinv{ a} +\frac{m-1}{2}\;h$ if $m$ is odd, or to $\qinv{  a}+\qinv{  a\cdot D}+\frac{m-2}{2}\;h$ if~$m$ is even, in $\GW(L)$.
Now let $M$ be the Gram matrix of $\Tr_{E/L}(\qinv{ m\cdot a\cdot x})$.
Then the $(i,j)$-th entry %~$M_{ij}$ 
equals
\[
%\Tr_{E/L}(\qinv{m\cdot a\cdot x \cdot x^{i-1}\cdot x^{j-1}})
M_{ij}=\begin{cases} m^2\cdot a\cdot D &\text{if } i+j=n\\
 0 & \text{otherwise.} 
\end{cases}\]
by Lemma \ref{lm:trace} and we get that 
\begin{align*}%&\begin{pmatrix}\Tr(m\cdot a\cdot x)& \Tr(m\cdot a\cdot x^2)&\ldots&\ldots& \ldots & \Tr(m\cdot a\cdot D)\\
%\Tr(m\cdot a\cdot x^2)& \ldots&\ldots&\ldots&  \Tr(m\cdot a\cdot D) & \Tr(m\cdot a\cdot D\cdot x)\\ \ldots & \ldots &\ldots& \ldots& \Tr(m\cdot a\cdot d) & \ldots\\
%\ldots & \ldots &\ldots& \ldots& \ldots & \ldots\\
%\ldots & \Tr(m\cdot a\cdot D) &\Tr(m\cdot a\cdot D\cdot x) & \ldots& \ldots  & \ldots\\
%\Tr(m\cdot a\cdot D)& \Tr(m\cdot a\cdot D\cdot x) & \ldots&\ldots &\ldots & \Tr( m\cdot a\cdot D\cdot x^{m-1})\end{pmatrix}\\
%\overset{\ref{lm:trace}}{=}&
M=\begin{pmatrix}0& 0&\ldots&\ldots& \ldots & m^2\cdot a\cdot D\\
0& \ldots&\ldots&\ldots&   m^2\cdot a\cdot D& 0\\ \ldots & \ldots &\ldots& m^2\cdot a\cdot D & 0&0\\
\ldots & \ldots &\ldots& \ldots& \ldots & \ldots\\
 0 &m^2\cdot a\cdot D &0& \ldots& \ldots  & \ldots\\
 m^2\cdot a\cdot D & 0&\ldots &\ldots& \ldots & 0\end{pmatrix}
\end{align*}
which is the Gram matrix of a quadratic form with class in $\GW(L)$ equal to 
%$\begin{cases}\qinv{ a\cdot D} +\frac{m-1}{2}\;h & m\text{ odd} \\\frac{m}{2}\;h &  m \text{ even}\end{cases}$
$\qinv{ a\cdot D} +\frac{m-1}{2}\;h$ if $m$ is odd, or to $\frac{m}{2}\;h$ if $m$ is even,
by Remark \ref{rmk:hyperbolic}.
\end{proof}

\subsection{The $\A^1$-degree}
$\A^1$-homotopy theory is a new branch of mathematics in which one aims to apply techniques from homotopy theory to the category of smooth algebraic varieties over a field $k$. 
Most constructions from classical homotopy theory work in this set up. In particular, we have an analog of the Brouwer degree. Recall (for example from \cite{hatcher}) that the Brouwer degree from classical topology is an isomorphism from the homotopy classes of the endomorphisms of the $n$-sphere to the integers
\[\deg:[S^n,S^n]\overset{\cong}{\longrightarrow} \Z\]
for $n\ge 1$.
Morel defines the $\A^1$-analog in \cite{morel}. His \emph{$\A^1$-degree} assigns an element of $\GW(k)$ to an $\A^1$-homotopy class of an endomorphism of the motivic sphere $\mathbb{P}^n_k/\mathbb{P}^{n-1}_k$
\[\deg^{\A^1}:[\mathbb{P}^n_k/\mathbb{P}^{n-1}_k, \mathbb{P}^n_k/\mathbb{P}^{n-1}_k]_{\A^1}\longrightarrow \GW(k).\]

Just like for the classical Brouwer degree, the $\A^1$-degree splits up as a sum of \emph{local $\A^1$-degrees}. We refer to \cite{KW-EKL} for the definition of the local $\A^1$-degree and merely recall some of their formulas to compute the local $\A^1$-degree $\deg_xf$ of a map $f:\A^n_k\longrightarrow \A^n_k$ at an isolated zero $x$.

%\begin{df}
%Assume $x$ is closed point of $f:\A^n_k\longrightarrow \A^n_k$ such that $f(x)$ is $k$-rational. Choose a neighborhood $U$ of $x$ such that $f$ maps $U\setminus \{x\}$ into $\A^n_k\setminus \{f(x)\}$. The \emph{local $\A^1$-degree} of $f$ at $x$ is the $\A^1$-degree of
%\[\mathbb{P}^n_k/\mathbb{P}^{n-1}_k\longrightarrow \mathbb{P}^n_k/(\mathbb{P}^{n}_k\setminus \{x\})\xleftarrow{\cong} U/U\setminus \{x\}\xrightarrow{\bar{f}}\A^n_k/(\A^n_k\setminus\{f(x)\})\xleftarrow{\cong}\mathbb{P}^n_k/\mathbb{P}^{n-1}_k.\]
%\end{df}

\subsubsection{Formulas for the local $\A^1$-degree}
Assume $f:\A^n_k\longrightarrow \A^n_k$ has an isolated zero $x$ with residue field $\kappa(x)$ separable over $k$. Furthermore, assume that the determinant of the Jacobian $\Jac(f)$ of $f$ at $x$ does not vanish.
In this case the local $\A^1$-degree at $x$ equals 
%\begin{equation}
%\label{eq:rational}
%\deg^{\A^1}_xf=\qinv{\det \Jac f(x)}\in \GW(k).\end{equation}
%\begin{ex}
%Let $f:\A^1_k\longrightarrow \A^1_k$ be defined by $f(x)=ax$ for some $a\in k^{\times}$. Then the local $\A^1$-degree at $x=0$ equals 
%\[\qinv{ \det \Jac f(0)}=\qinv{ a}\in \GW(k).\]
%\end{ex}
%In the case when $k(x)$ is separable over $k$, the local $\A^1$-degree of $f$ at $x$ is given by
\begin{equation}
\label{eq:nonrational}
\deg^{\A^1}_xf=\Tr_{\kappa(x)/k}\qqinv{\det \Jac f(x)}\in \GW(k).\end{equation}
There are also formulas for the local $\A^1$-degree in case $\det \operatorname{Jac} f(x)=0$ or $k(x)$ is not separable over~$k$ \cite{KW-EKL,BBMMO,BMP}, but in this paper we restrict to the case of zeros with a non-vanishing Jacobican determinant with residue field separable over our base field $k$.

\newpage
\subsection{The Poincar\'e-Hopf theorem and the $\A^1$-Euler number}
\subsubsection{Motivation from classical topology}
Let $V\longrightarrow X$ be an oriented vector bundle of rank $r$ on a smooth, closed, connected, oriented manifold $X$ of dimension $r$. The \emph{Euler number} $n(V)$ is the Poincar\'e dual of the Euler class $e(V)$
\[n(V):=e(V)\cap[X]\in H_0(X,\Z)\cong \Z.\]
Assume $\sigma:X\longrightarrow V$ is a section and $x\in X$ an isolated zero of $\sigma$. Choose oriented coordinates around $x$ and a trivialization of~$V$ in a neighborhood around $x$ compatible with the orientation of $V$. In these coordinates and trivialization, the section $\sigma$ is a map $\sigma:\R^r\longrightarrow \R^r$.
The \emph{local index}~$\operatorname{ind}_x\sigma$ of $\sigma$ at $x$ is the local Brouwer degree of $\sigma$ at $x$.

\begin{thm}[Poincar\'e-Hopf Theorem]
\label{thm:PoincareHopf}
Let ${\sigma:X\longrightarrow V}$ be a section 
of the bundle ${V\longrightarrow X}$, having only isolated zeros. Then
\[n(V)=\sum_{x:\sigma(x)=0}\operatorname{ind}_x\sigma.\]

\end{thm}

\subsubsection{$\A^1$-Euler number}

Kass and Wickelgren define the \emph{$\A^1$-Euler number} of a \emph{relatively oriented} vector bundle $V\longrightarrow X$ of rank $r$ on a $r$-dimensional smooth, proper variety $X$ over $k$ as the sum of local indices defined using the local $\A^1$-degree analogous to the Poincar\'e-Hopf theorem. We recall the following definitions from \cite{KW-Euler}.
\begin{df}
\label{df:relativeorientation}
Let $V\longrightarrow X$ be a vector bundle. A \emph{relative orientation} of ${V\longrightarrow X}$ consists of a line bundle $\mathcal{L}\longrightarrow X$ and an isomorphism \[\rho:\det V\otimes \omega_{X/k}\longrightarrow \mathcal{L}^{\otimes 2}.\] Here $\omega_{X/k}$ is the canonical line bundle on $X$.
\end{df}
\begin{ex}
Since $\omega_{\mathbb{P}^1_k/k}=\mathcal{O}_{\mathbb{P}^1_k}(-2)$, the line bundle $\mathcal{O}_{\mathbb{P}^1_k}(n)\longrightarrow \mathbb{P}^1_k$ is relatively oriented if and only if $n$ is even.
\end{ex}

We need to restrict to relatively oriented bundles since otherwise we would not have a well-defined local index: In order to define the local indices at all the zeros in a consistent way, we need to choose coordinates (called Nisnevich coordinates) and a trivialization of the vector bundle \emph{compatible} with the coordinates and the relative orientation. This means that the section \allowbreak of ${\det V\otimes \omega_{X/k}}$ defined by the chosen coordinates and the chosen trivialization is sent to a square in~$\mathcal{L}^{\otimes 2}$ by $\rho$. Hence, different choices of coordinates and trivializations compatible with the relative orientation only differ by a square, so they do not differ in $\GW(k)$.

\begin{df}
Let $X$ be a smooth and proper $k$-scheme of dimension $r$ and let $x\in X$ be a closed point. An \'{e}tale map $\psi:U\longrightarrow \A^r_k$ from a Zariski neighborhood $U$ of $x$ which induces an isomorphism of residue fields of $x$ is called \emph{Nisnevich coordinates} around $x$.
\end{df}

\begin{rmk}
\label{rmk:Nisnevichcoordinates}
Nisnevich coordinates always exist given that $r\ge 1$ by \cite[Lemma 19]{KW-Euler}.
Since~$\psi$ in the definition of Nisnevich coordinates is \'{e}tale,  the standard basis for the
tangent space of $\A^r_k$ defines a trivialization of $TX\vert_U$ where $TX$ is the tangent bundle of $X$.
\end{rmk}

\begin{df}
Let $V\longrightarrow X$ be a vector bundle of rank $r$ over an $r$-dimensional scheme $X$ over~$k$ equipped with a relative orientation $\rho:\det V\otimes \omega_{X/k}\longrightarrow \mathcal{L}^{\otimes 2}$ and let $\psi:U\longrightarrow \A^r_k$ be Nisnevich coordinates around a closed point $x\in X$. By Remark \ref{rmk:Nisnevichcoordinates}, a choice of Nisnevich coordinates defines a section of $\det TX\vert_U$. A trivialization $V\vert_U\cong U\times \A^r$ defines a section of $\det V\vert_U$. We say a trivialization of $V\vert_U$
is \emph{compatible} with the relative orientation $\rho$ and the Nisnevich coordinates if the section of $\det V\vert_U\otimes (\det TX\vert_U)^{\vee}$, equivalently, $\det V\vert_U\otimes \omega_{X/k}\vert_U$ defined by the trivialization and the Nisnevich coordinates is sent to a square by $\rho$, that is to a section of $\mathcal{L}^{\otimes 2}$ of the form $\ell\otimes \ell$.
\end{df}

We are now ready to define the local index at an isolated zero $x$ of a section $\sigma$ valued in $\GW(k)$.
Let $\sigma$ be a section of a relatively oriented vector bundle $V\longrightarrow X$ and let $x$ be an isolated zero of~$\sigma$. Choose Nisnevich coordinates $\psi:U\longrightarrow \A^r_k$ and a trivialization $\phi$ of $V\vert_U$ compatible with the relative orientation of $V\longrightarrow X$ and the Nisnevich coordinates around $x$. 
\begin{df}
\label{df:local index}
The \emph{local index} $\operatorname{ind}_x\sigma$ of $\sigma$ at $x$ is the local $\A^1$-degree of 
\[ U\xrightarrow{\sigma\vert_U} V\vert_U\cong U\times \A^r\xrightarrow{\operatorname{pr}_2}\A^r \]
at $x$.
\end{df}

Now assume that $V\longrightarrow X$ is a relatively oriented vector bundle of rank $r$ over a smooth proper $r$-dimensional $k$-scheme $X$ and $\sigma$ is a section with only isolated zeros.

\begin{df}[Kass-Wickelgren]
The $\A^1$-Euler number of $V\longrightarrow X$ is the sum of local indices at the zeros of $\sigma$
\[n^{\A^1}(V)\coloneqq\sum_{x:\sigma(x)=0}\operatorname{ind}_x\sigma \in \GW(k).\]
\end{df}
By \cite[Theorem 1.1]{BW}, this is independent of the choice of section.

\subsection{B\'ezout's Theorem enriched in $\GW(k)$}
The classical B\'{e}zout theorem for algebraically closed fields counts the intersections points of $n$ hypersurfaces in $\mathbb{P}^n$ defined by homogeneous polynomials $f_1,\ldots, f_n$ of degrees $d_1,\ldots,d_n$, respectively.
A homogeneous polyonomial $f$ in $n+1$ variables of degree $d$ defines a section of $\glob_{\mathbb{P}^n}(d)\longrightarrow \mathbb{P}^n$. So $f_1,\ldots,f_n$ define a section of 
\[V:=\glob_{\mathbb{P}^n}(d_1)\oplus \ldots\oplus \glob_{\mathbb{P}^n}(d_n)\longrightarrow \mathbb{P}^n\]
and the zeros of this section are exactly the intersection points of $f_1,\ldots,f_n$.
Since we have that~$\omega_{\mathbb{P}^n/k}\cong\glob(-n-1)$, the bundle $V$ is relatively oriented if and only if $\sum_{i=1}^nd_i-n-1$ is even.
In this case, McKean computes the $\A^1$-Euler number yielding an enrichment of B\'{e}zout's theorem in $\GW(k)$.
\begin{thm}[McKean]
\label{thm:McKean}
\[n^{\A^1}(V)=\sum \operatorname{ind}_x(f_1,\ldots,f_n)=\frac{d_1\cdot\ldots \cdot d_n}{2}\cdot h\in \GW(k)\]
where the sum runs over the intersection points of $f_1,\ldots,f_n$.
\end{thm}
McKean uses the standard open affine subsets $U_i=\{x_i\neq 0\}$ of $\mathbb{P}^n_k$ as Nisnevich coordinates and the usual trivialization of $V\vert_{U_i}$. In particular, the section $(f_1,\ldots,f_n)$ in these coordinates and trivialization becomes
\[(f_1,\ldots,f_n)\colon U_i\cong \A^n\longrightarrow \A^n\]
setting $x_i=1$.

\subsubsection{Non-orientable case and representability of the $\A^1$-degree}
When $V$ is not relatively orientable, McKean shows that one can still orient $V$ relative to the divisor $D=\{x_0=0\}$ in the sense of Larson and Vogt \cite{LV}. 
Geometrically, this counts the intersection points in $\A^n\cong U_0=\{x_0\neq0\}\subset\mathbb{P}^n$. %which equals the sum of local $\A^1$-degree the zeros of 
%\[f:\A^n\longrightarrow \A^n\]
%where $f=(f_1,\ldots,f_n)$ and $\A^n$ is $U_0$.
In section \ref{section:non orientable} we explain how to get all possible counts for B\'ezout in this case using (enriched) tropical methods. In particular, we will see that we cannot get any element of $\GW(k)$. More precisely, we show find a lower bound for the number of hyperbolic summands in Corollary \ref{cor:non-orientable}.
%In \cite{QSW} it is asked which quadratic forms are representable as the local $\A^1$-degree. In the classical setting, each integer can be realized as the local degree. This is not the case for the local $\A^1$-degree. For example, the only quadratic form of rank $2$ which is representable, is the hyperbolic form $h$. Eisenbud and Levine find lower bounds for the number of hyperbolic summands of a quadratic form which is representable by the local $\A^1$-degree \cite{EL}.

%Similarly, one can ask which quadratic forms represent the sum of local $\A^1$-degrees at the zeros of a map $f:\A^n_k\longrightarrow \A^n_k$ with only isolated zeros. 
%We answer this question in the tropical set up...

%% file: 3settingtrop.tex
% !TeX root=main.tex

\section{Introduction to Tropical Geometry}
\label{section:IntroTropicalGeometry}
%{\color{red} should we somewhere in this section say something about Puiseux series over $k$ with positive characteristic?}

%{\color{red} In the whole section or even paper: should we use different notation for the tropical polynomials vs the polynomials over $k\Puiseux$? For example kapital letters for polynomials over $k\Puiseux$?}

In this section we introduce the basic notions of tropical geometry we use in the subsequent sections. For more details in tropical geometry we refer the reader to \cite{BIMS}, \cite{IMS} and \cite{Mikh}.

\subsection{Toric deformations}
\label{subsection:ToricDeformations}
Given a Laurent polynomial $f=\sum_{I\in A}\alpha_I x^{i_1}y^{i_2}$ in $ k[x,y]$, where $A\subset \Z^2$ is a finite set of tuples $I=(i_1,i_2)$, we consider a \emph{toric deformation} of $f$, that is a family of polynomials given by
\[f_t(x,y)=\sum_{I\in A}\alpha_I x^{i_1}y^{i_2}t^{\varphi(I)},\]
where $\varphi:A\longrightarrow \Q$ is the restriction of a convex rational function to the set of indices $A$.
We can think of the variable $t$ as the variable of $\mathbb{G}_m$ whose specialization to $t=1$ is our initial polynomial.
The family $f_t$ can be seen as an element of the polynomial ring $k\Puiseux[x,y]$ with coefficients in the field of Puiseux series $k\Puiseux$.
The field of Puiseux series has a valuation given by
\[
\operatorname{val}\colon \begin{array}{ccc}
k\Puiseux  & \longrightarrow & \Q\cup \{\infty\} \\
 \sum_{i=i_0}^{\infty} a_i t^{i/N}   & \longmapsto & i_0/N\\
 0&\longmapsto \infty
\end{array}
\]
given that $a_{i_0}\neq 0$.
Let $\nu\coloneqq -\val\colon k\Puiseux\longrightarrow\mathbb{Q}\cup \{-\infty\}$.
Then $\nu$ satisfies
\[
\begin{array}{rcl}
\nu(x+y) & = & \max\{\nu(x),\nu(y)\} \quad\text{ if } \nu(x)\neq\nu(y),\\
\nu(xy)& = &\nu(x)+\nu(y).
\end{array}
\]
Note that these are exactly the operations in the tropical semifield introduced in the next subsection. 
Given a toric deformation $f_t\in k\Puiseux[x,y]$ as above, assume that $\operatorname{char}k=0$ or $\operatorname{char}k>\max\{\deg(f_t),\deg (g_t)\}$. Then we have that
\[f_t(x,y)=\sum_{I\in A}\alpha_Ix^{i_1}y^{i_2}t^{\varphi(I)}=0\]% \text{ and, }
%g_t(x,y)=\sum_{J\in B}\beta_Jx^{j_1}y^{j_2}t^{\psi(J)}=0\]
has a solution in $\overline{k}\Puiseux^2$ given by %{\color{red} do we need a characteristic assumption here!?}
\[
x(t)  =  x_0 t^{i_0}+\text{ higher order terms in } t,\quad y(t) =  y_0 t^{j_0}+\text{ higher order terms in } t,
\]
% \[
% \begin{array}{rcl}
% x(t) & = & x_0 t^{i_0}+\text{ higher order terms in } t,\\
% y(t) & = & y_0 t^{j_0}+\text{ higher order terms in } t,
% \end{array}
% \]
that is,
\begin{align*}
    0=f_t(x(t),y(t))&=\sum_{I\in A}\left(\alpha_Ix_0^{i_1}y_0^{i_2}t^{\varphi(I)-i_1\nu (x)-i_2\nu(y)}+ \text{ higher order terms in }t\right)%\\
   % 0=g_t(x(t),y(t))&=\sum_{J\in B}\left(\beta_Jx_0^{j_1}y_0^{j_2}t^{\psi(J)-j_1\nu (x)-j_2\nu(y)}+ \text{ higher order terms in }t\right),
\end{align*}
if and only if
the term of lowest power in $t$, that is where $t$ has the exponent
\[\{\varphi(I)-i_1\nu (x)-i_2\nu(y): I\in A\}\]%\textrm{ and }\{\psi(J)-j_1\nu (x)-j_2\nu(y): J\in B\}\]
appears at least twice in $f_t(x(t),y(t))$ %and $g_t(x(t),y(t))$. 
Equivalently, the maximum of
\[\{-\varphi(I)+i_1\nu (x)+i_2\nu(y): I\in A\}\]%\textrm{ and }\{-\psi(J)+j_1\nu (x)+j_2\nu(y): J\in B\}\]
has to be attained at least twice.
This is exactly the definition of the tropical vanishing locus in the next section. For example, if $f_t$ is a polynomial of degree $1$, then the locus where the maximum is attained at least twice is a tropical line and looks like Figure \ref{figure:exampletroplineandconic} on the left.

These notions extend naturally to more variables.

\subsection{Tropical hypersurfaces and tropicalization maps}

\subsubsection{Tropical hypersurfaces}

The tropical semifield is the set $\mathbb{T}=\mathbb{R}\cup\{-\infty\}$ endowed with the operations (denoted by ``$+$" and~``$\cdot$")
\[
\begin{array}{rcl}
`` x+y" & = & \max\{x,y\},\\
`` x\cdot y" & = &x+y.
\end{array}
\]
Then $\mathbb{T}$ with these two operations forms a semifield, i.e., it satisfies all axioms of a field but the existence of additive inverse.
We write $\mathbb{T}^*$ for $\mathbb{T}\setminus \{-\infty\}=\R$.
A \emph{tropical polynomial} in $n$ variables is a polynomial %$p\in \mathbb{T}[x_1,\ldots,x_n]$ 
given by
\begin{equation}
\label{eq:tropPolynvariables}
f(x)=``\sum_{I\in A}a_Ix^I"\coloneqq ``\sum_{I\in A}a_Ix_1^{I_1}\cdots x_n^{I_n}"\in \mathbb{T}[x_1,\ldots,x_n],\end{equation}
where $x=(x_1,\ldots,x_n)$, $A$ is a finite set of tuples $I=(I_1,\ldots,I_n)\in \Z^n$ and~$a_I\neq -\infty$ for $I\in A$. 

% \begin{ex}
% A polynomial $p\in\mathbb{T}[x]$ in one variable has the form
% \[
% p(x)=``\sum_{i=0}^{d} a_i x^i"=\max_{i=0}^{d}\{a_i+i x\}.
% \]
% \end{ex}

%A polynomial $p\in \mathbb{T}[x_1,\ldots,x_n]$ 
Such polynomial defines a function in~$(\mathbb{T}^*)^n=\R^n$ that is piecewise linear.  
Its \emph{tropical vanishing locus} is defined as the locus of non-differentiability, i.e., the points in $\R^n$ such that the maximum is obtained at least twice. We denote this locus by~$V_{\operatorname{Trop}}$,
and it is expressed by
\[V_{\operatorname{Trop}}(f)=\{x\in \R^n\mid \exists I,I'\in A: I\neq I', p(x)=I\cdot x+a_I=I'\cdot x+a_{I'}\},\]
where $\cdot$ denotes the scalar product.

\begin{ex}
Let \[f(x,y)=``(-3)x+(-3)y+0"=\max\{x-3,y-3,0\}.\]
Then $V_{\operatorname{Trop}}(f)$ is the \emph{tropical line} in Figure \ref{figure:exampletroplineandconic} (a). In the left lower component of $\R^2\setminus V_{\operatorname{Trop}}(f)$, $0$ is maximal, in the component on the right, $x-3$ is maximal and in the upper left component, $y-3$ is maximal.
The right picture in Figure \ref{figure:exampletroplineandconic} shows a \emph{tropical conic}, i.e., it is the tropical vanishing locus of a tropical polynomial of degree $2$.
\begin{figure}[t]
     \begin{subfigure}[b]{0.45\textwidth}
         \centering
         \begin{tikzpicture}[scale=0.9]
\draw[blue,thick] (-3,0)--(0,0)--(0,-3);
\draw[blue,thick] (0,0)--(2,2);
\node at (0.55,-0.15) {$(3,3)$};
\end{tikzpicture}
         \caption{tropical line}
     \end{subfigure}
     \hfill
     \begin{subfigure}[b]{0.45\textwidth}
         \centering
         \begin{tikzpicture}[scale=0.6]
\draw[red,thick] (-1,2)--(2,2)--(2,-1);
\draw[red,thick] (2,2)--(3,3)--(4,3)--(4,-1);
\draw[red,thick] (4,3)--(6.5,5.5);
\draw[red,thick] (3,3)--(3,4)--(5.5,6.5);
\draw[red,thick] (3,4)--(-1,4);
\end{tikzpicture}
         \caption{tropical conic}
     \end{subfigure}
        \caption{Examples of tropical curves}
        \label{figure:exampletroplineandconic}
\end{figure}
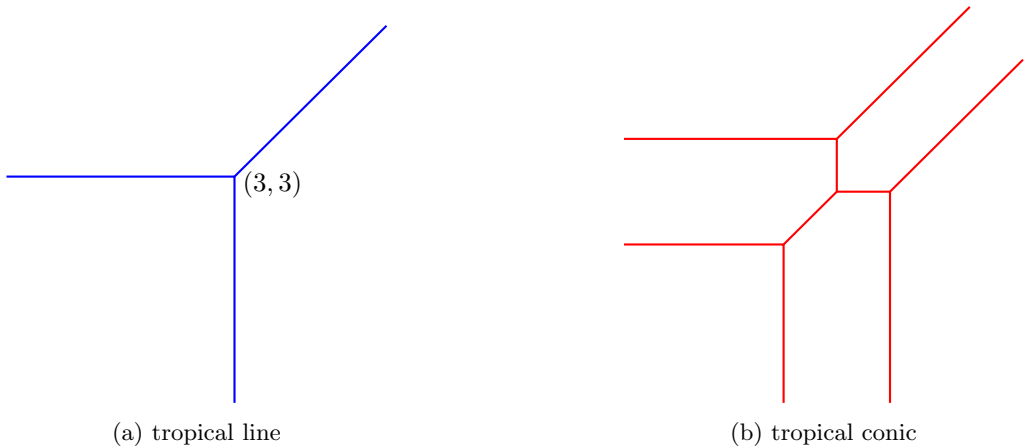
\end{ex}

The \emph{Newton polytope} of a tropical polynomial $f$ is given by the convex hull of its monomials with non-trivial coefficient in $\R^n$ 
\[
\operatorname{NP}(f)=
\operatorname{Conv}\left(\{
I\in \R^n \mid a_I\neq-\infty
\}\right)=\operatorname{Conv}(A).
\]
We say that a polynomial $f$ has \emph{degree} $d$ if its Newton polytope is 
\begin{equation}
\label{eq:deltad}
\Delta_d=\operatorname{Conv}\{(0,0,\ldots,0),(d,0,\ldots,0),(0,d,0,\ldots,0),\ldots,(0,\ldots,0,d)\}\subset\R^n.
\end{equation}

\subsubsection{Tropicalization maps}

Using the aforementioned map
\[
\nu\colon k\Puiseux \longrightarrow \Q\cup \{-\infty\}
\]
we can \emph{tropicalize} a polynomial over $k\Puiseux$ by taking $\nu=-\operatorname{val}$ of its coefficients and reinterpreting the addition and multiplication
\[
``\cdot"\colon \begin{array}{cccl}
k\Puiseux[x_1,\dots,x_n]     & \longrightarrow & \mathbb{T}[x_1,\dots,x_n] &\\
 \sum_{I\in A} \alpha_I(t) x^I   & \longmapsto & ``\displaystyle\sum_{I\in A} \nu(\alpha_I(t)) x^I"&=\displaystyle\max_{I\in A}\{\nu(\alpha_I(t)) +I\cdot x\}.
\end{array}
\]

Now let $f\in k\Puiseux[x_1,\ldots,x_n]$ and let $k$ be a field with $\operatorname{char}k=0$ or $\operatorname{char}k>\deg f$. Then we can also \emph{tropicalize} the zeros of a polynomial $f$, %(x)=\sum_{I\in A}\alpha_I(t)x^I$, 
at the level of sets, by taking the closure  of the image of the valuation taken point-wise in $\R^n$. More precisely, let $X$ be the vanishing locus  $V(f)\subset\mathbb{A}^n_{k\Puiseux}$.
Then 
\[
\operatorname{Trop}(X)\coloneqq
\overline{\left\{
(\nu(x))\mid x\in X
\text{ geometric point}\right\}}\subset\R^n,
\]
where $\nu(x)\in\R^n$ is given by $\nu$ coordinatewise.

If $X$ is an algebraic hypersurface, its tropicalization is a tropical hypersurface, defined by the tropicalization of a defining polynomial for $X$. %given that $k$, and thus also $k\Puiseux$, is algebraically closed and of characteristic $0$. 
%In other words, we have the following theorem. %{\color{red} I added the characteristic $0$ assumption here, otherwise the sentence was wrong. Do we need characteristic $0$ in the next theorem?\color{blue} also fix the arxiv version.}

\begin{thm}[Kapranov]
If $k$ has characteristic $0$ or $\deg f<\operatorname{char}k$, then for a polynomial $f\in k\Puiseux[x_1,\ldots,x_n]$ one has
\[V_{\operatorname{Trop}}(``f")=\operatorname{Trop}(V(f)).\]
\end{thm}
%This also works in positive characteristic if $\deg f<\operatorname{char}k$.

\subsubsection{Dual subdivision}
We associate a refinement of the Newton polytope $\NP(f)$ called the \emph{dual subdivision} $\operatorname{DS}(f)$ of $f$.
%First, the \emph{Newton polytope} $\operatorname{NP}$ of $p$ is the convex hull
%\[
%\operatorname{NP}(p)=
%\operatorname{Conv}\left(\{
%(i_1,i_2)\mid a_{(i_1,i_2)}\neq-\infty
%\}\right)=\operatorname{Conv}(A),
%\]
The refinement is given by the projection to $\R^n$ of the boundary of the upper faces (with respect to the last coordinate) of the polyhedron
\[
\operatorname{Conv}\left(\{(I,a_I)\in\R^{n+1})\mid a_I\neq-\infty\}\right).
\]
There is a one-to-one correspondence of the elements

\begin{figure}[t]
     \begin{subfigure}[b]{0.45\textwidth}
         \centering
         \includegraphics[width=0.6\textwidth]{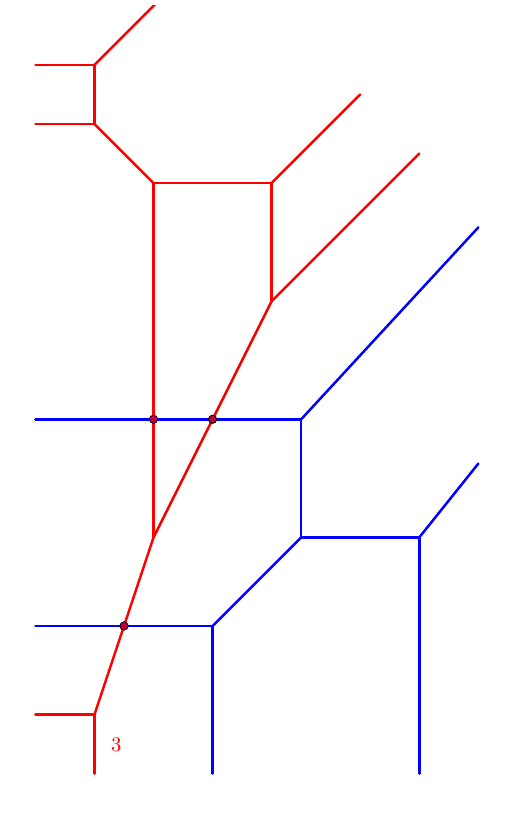}
         \caption{A cubic curve and a conic curve.}
     \end{subfigure}
     \hfill
     \begin{subfigure}[b]{0.45\textwidth}
         \centering
         \includegraphics[width=\textwidth]{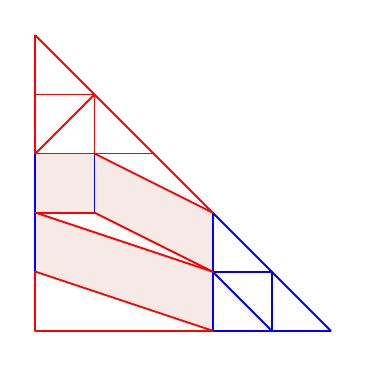}
         \caption{The dual subdivision of the union of the curves in (a).}
     \end{subfigure}
        \caption{Example of the dual subdivision.}
        \label{fig:intcurves}
\end{figure}

\begin{center}
\begin{tabular}{ |c|c| } 
 \hline
 $V_{Trop}(f)$ & $\operatorname{DS}(f)$ \\
  \hline
 vertex $v$ &                   connected component of $\operatorname{NP}(f)\setminus \operatorname{DS}(f)$\\ 
 $l$-dimensional face $e$ &         $n-l$-dimensional face $e'$\\ 
 connected component $K$ of $\R^n\setminus V_{\operatorname{Trop}}(f)$ & vertex $v_K$ \\ 
 \hline
\end{tabular}
\end{center}
Moreover, the corresponding dual faces $e$ and $e'$ are orthogonal, and inclusion of faces are inverted. 
Figure \ref{fig:dual subdivision} shows a tropical conic with its dual subdivision. Figure~\ref{fig:intcurves} shows the dual subdivision of a reducible quintic.
\begin{figure}[t]
     \begin{subfigure}[b]{0.45\textwidth}
         \centering
         \begin{tikzpicture}[scale=0.6]
\draw[red,thick] (-1,2)--(2,2)--(2,-1);
\draw[red,thick] (2,2)--(3,3)--(4,3)--(4,-1);
\draw[red,thick] (4,3)--(6.5,5.5);
\draw[red,thick] (3,3)--(3,4)--(5.5,6.5);
\draw[red,thick] (3,4)--(-1,4);
\end{tikzpicture}
         \caption{Tropical conic}
     \end{subfigure}
     \hfill
     \begin{subfigure}[b]{0.45\textwidth}
         \centering
         \begin{tikzpicture}
\filldraw[fill=gray!10!white,thick] (0,0)--(4,0)--(0,4)--cycle;
\draw[red,thick] (0,0)--(4,0)--(0,4)--cycle;
\draw[red,thick] (2,0)--(0,2)--(2,2)--cycle;
\end{tikzpicture} 
         \caption{Dual subdivision}
     \end{subfigure}
        \caption{A tropical conic with its dual subdivision.}
        \label{fig:dual subdivision}
\end{figure}
\subsubsection{Tropical intersections}
We say that $n$ tropical hypersurfaces $V_1,\ldots,V_n$ in $\R^n$ intersect \emph{tropically transversely} at $p\in V_1\cap\ldots \cap V_n$ if the point~$p$ is an isolated point belonging to the interior of a top dimensional face of $V_i$, for every $i=1,\dots, n$.
In particular, for $n$ tropical hypersurfaces $V_1,\ldots,V_n$ in $\R^n$ that intersect tropically transversely at~$p\in V_1\cap\ldots \cap V_n$, i.e., the point $p$ corresponds to a parallelepiped in the dual subdivision of the union $V_1\cup \ldots\cup V_n$. %Just as in the case for curves \eqref{eq:Mult and area parallelogram} 
\begin{df}
\label{df:nonenrichedintmult}
Assume $V_1,\ldots,V_n$ intersect tropically transversally at $p$.
The \emph{intersection multiplicity} of $V_1,\ldots,V_n$ at $p$ is the volume of the parallelepiped dual to $p$:
\[\operatorname{mult}_p(V_1,\dots,V_n)=\operatorname{Vol(\text{parallelepiped dual to }}p).\]
\end{df}
This definition is meaningful by the following result of Huber-Sturmfels \cite{HuberSturmfels}. 
Their argument also works in positive characteristic when the characteristic is big enough, that is all equations involved have degree less than the characteristic for example if the characteristic is bigger than the maximum of the diameters $\max\{\operatorname{diam}(NP(F_i))$ of all the Newton polytopes of the tropical polynomials $F_i$ defining the tropical hypersurfaces $V_i$.
\begin{lm}
\label{lm:tropintmult}
    Assume that $k$ is a field of $\operatorname{char}k=0$ or $\operatorname{char} k> \max\{\operatorname{diam}(\operatorname{NP}(F_i))\}$.
    Let $F_1,\ldots,F_n$ be~$n$ general Laurent polynomials in $k\Puiseux[x_1,\ldots,x_n]$. These define tropical polynomials $f_1,\ldots,f_n$.
    Let $V_i=V_{\operatorname{Trop}}(f_i)$ for $i=1,\ldots,n$ and $p\in V_1\cap\ldots\cap V_n$.
    Then the number of common zeros of the $F_i$ that tropicalize to $p$ equals $\operatorname{mult}_p(V_1,\ldots,V_n)$.
\end{lm}

Our main theorem \ref{thm:mainthm} is the quadaratic enrichment of this Lemma and in section \ref{section:applications} we will derive a quadractic enrichment of the following consequence of Lemma \ref{lm:tropintmult}.
Recall that the \emph{mixed volume} $\MV(\Delta_1,\ldots,\Delta_n)$ of $n$ polytopes in $\R^n$ is the coefficient of $\lambda_1\cdots \lambda_n$ in the polynomial $R(\lambda_1,\ldots,\lambda_n)$ given by 
\[\operatorname{Vol}(\lambda_1\Delta_1+\ldots+\lambda_n\Delta_n).\]

%{\color{blue} can we give a reference and a proof sketch of this theorem since it is kind of crucial for our result!}
\begin{thm}[Tropical B\'{e}zout and Bernstein-Kushnirenko Theorem]
\label{thm:tropical Bezout}
Let $V_1,\ldots,V_n$ be tropical hypersurfaces in $\R^n$ with Newton polytopes $\Delta_{1},\ldots,\Delta_{n}$, respectively. Then 
\[\sum_{p}\operatorname{mult}_p(V_1,\ldots,V_n)=\MV(\Delta_1,\ldots,\Delta_n).\]
In particular, if $\Delta_i=\Delta_{d_i}=\operatorname{Conv}\{0,d_ie_1,\ldots,d_ie_n\}$ %=d_i\cdot \Delta_1$ where $\Delta_1$ is the convex hull
%$$\Delta_1=\operatorname{Conv}\{0,e_1,\ldots,e_n\}\subset \R^n$$
for $i=1,\ldots,n$, then we get the tropical B\'ezout theorem
\[\sum_{p}\operatorname{mult}_p(V_1,\ldots,V_n)=d_1\cdots d_n.\]
%{\color{red} we already defined $\Delta_d$ when we define the Newton polytope}
\end{thm}

\subsection{Enriched tropical hypersurfaces and Viro Polynomials}

Viro's patchworking is a combinatorial construction yielding topological properties of real algebraic varieties. It is an algorithmic construction whose input is a subdivision of a polytope and a set of signs $\sigma(I)$ (either plus or minus) for every integer point $I$ in the dual subdivision of the polytope. A Viro polynomial associated to this data is a polynomial
\[
\displaystyle\sum_{I\in\operatorname{DS}(\Delta)\cap\Z^2}\sigma(I)x^I t^{\varphi(I)}
\]
where $\varphi$ is a convex piece-wise linear function inducing the subdivision and such that tropicalizing the polynomial yields back a defining polynomial for the tropical curve.
Based on this idea, we generalize this concept by replacing the signs $\sigma(I)$ with elements $\alpha_I\in k^\times/(k^\times)^2$ and call the following \emph{enriched Viro polynomial}

\begin{equation}
    \label{eq:enrichedViroPolynomial}
    \displaystyle\sum_{I\in\Delta\cap\Z^n} \alpha_Ix^It^{\varphi(I)}.
\end{equation}

Note that if $k=\R$, then this coincides with the original definition of a Viro polynomial since~$\R^\times/(\R^\times)^2=\{\pm1\}$.
Tropicalization gives back a tropical hypersurface $V$ which has the dual subdivision we started with. However, in the tropicalization process one loses information, namely the elements in $k^\times/(k^\times)^2$. We would like to remember these \emph{coefficients} by assigning them to the corresponding connected component in $\R^2\setminus V$, that is we assign the coefficient $\alpha_I$ of a monomial $x^It^{\varphi(I)}$ to the component where the $\nu(x^It^{\varphi(I)})=I\cdot \nu(x)-\varphi(I)$ attains the maximum. Equivalently, one can assign the coefficients $\alpha_I$ to the corresponding vertex $I$ in the dual subdivision.
This gives rise to the following definition.

\begin{figure}[t]
     \centering
     \begin{subfigure}[b]{0.45\textwidth}
         \centering
         \begin{tikzpicture}
\draw[blue,thick] (-3,0)--(0,0)--(0,-3);
\draw[blue,thick] (0,0)--(2,2);
\node[blue] at (2,-0.5) {$\alpha_{(1,0)}$};
\node[blue] at (-0.5,2) {$\alpha_{(0,1)}$};
\node[blue] at (-1,-1) {$\alpha_{(0,0)}$};
\node[gray,scale=0.7] at (0.5,-0.1) {(3,3)};
\end{tikzpicture}
         \caption{An enriched tropical line.}
     \end{subfigure}
     \hfill
     \begin{subfigure}[b]{0.45\textwidth}
         \centering
         \begin{tikzpicture}[scale=0.7]
\draw[red,thick] (-1,2)--(2,2)--(2,-1);
\draw[red,thick] (2,2)--(3,3)--(4,3)--(4,-1);
\draw[red,thick] (4,3)--(6.5,5.5);
\draw[red,thick] (3,3)--(3,4)--(5.5,6.5);
\draw[red,thick] (3,4)--(-1,4);
\node[red] at (3,0) {$\beta_{(1,0)}$};
\node[red] at (0,3) {$\beta_{(0,1)}$};
\node[red] at (0,0) {$\beta_{(0,0)}$};
\node[red] at (5,1) {$\beta_{(2,0)}$};
\node[red] at (1,5) {$\beta_{(0,2)}$};
\node[red] at (5,5) {$\beta_{(1,1)}$};
\node[gray,scale=0.7] at (2.5,1.9) {$(0,0)$};
\node[gray,scale=0.7] at (3.5,3.3) {$(2,2)$};
\node[gray,scale=0.7] at (4.5,2.9) {$(4,2)$};
\node[gray,scale=0.7] at (2.7,4.3) {$(2,4)$};
\end{tikzpicture}
         \caption{An enriched tropical conic.}
     \end{subfigure}
        \caption{Examples of enriched tropical curves}
        \label{figure:enrichedcurve}
\end{figure}

\begin{df}
An \emph{enriched tropical hypersurface} $\tV=(V,(\alpha_I))$ in $\R^n$ is a tropical hypersurface~$V$ in $\R^n$ together with a \emph{coefficient} $\alpha_I$ assigned to each connected component of $\R^n\setminus V$, or equivalently, to each vertex in the dual subdivision. We call such element $\alpha_I$ of $k^{\times}/(k^{\times})^2$ the \emph{coefficient} of the component/vertex of the dual subdivision. We write $V$ for the underlying non-enriched tropical hypersurface. 
\end{df}

To each enriched tropical hypersurface $\tV$ we can assign an \emph{(enriched) Viro polynomial} of the form \eqref{eq:enrichedViroPolynomial} such that tropicalizing and remembering the coefficients gives back $\tV$.

%{\color{blue} In this definition, we ask for the tropical multiplicity to be not divisible by the characteristic. 
%If we ask for the characteristic to be big enough, this condition is satisfied automatically so there is no need for it. But, is it equivalent? in such case I would prefer to delete the char assumption}

\begin{ex}
Figure \ref{figure:enrichedcurve} shows an enriched tropical line and an enriched tropical conic with enriched Viro polynomials
\[\alpha_{(0,0)}+\alpha_{(1,0)}xt^3+\alpha_{(0,1)}yt^3,\]
 and, respectively,
\[\beta_{(0,0)}+\beta_{(1,0)}x+\beta_{(0,1)}y+\beta_{(2,0)}x^2t^4+\beta_{(1,1)}xyt^2+\beta_{(0,2)}y^2t^4.\]
Tropicalizing the enriched Viro polynomial of the line yields the tropical polynomial
\[``0+(-3)x+(-3)y"=\max\{0,x-3,y-3\}\]
which has tropical vanishing locus the tropical line with $3$-valent vertex at $(3,3)$ displayed in Figure~\ref{figure:enrichedcurve}. The enrichment remembers the coefficients of the enriched Viro polynomials and assigns them to the connected components where the tropicalization of the corresponding monomial attains the maximum.

Similarly, one can see that the tropicalization of the enriched Viro polynomial of the conic yields the tropical conic in Figure \ref{figure:enrichedcurve} as well as the coefficients of the connected components.
\end{ex}

To prove a quadratically enriched tropical B\'ezout theorem, we also need to enrich the union of enriched tropical hypersurfaces $\tV_1\cup \ldots \cup \tV_n$. Note that to $f\in k\Puiseux[x_1,\ldots,x_n]$ we can associate an enriched Viro polynomial $f_{\circ}$ by replacing the coefficients of $f$ with its initials. Note that if $V_i$ is defined by a tropical polynomial $``f_i"$ for $i=1,\ldots,n$, then $``f_1\cdots f_n"$ is a tropical polynomial defining $V_1\cup \ldots \cup V_n$. The product of enriched Viro polynomials $f_1\cdots f_n$ is a polynomial in $k\Puiseux[x_1,\ldots,x_n]$ and to get an enriched Viro polynomial for $\tV_1\cup \ldots\cup \tV_n$ we take $(f_1\cdots f_n)_\circ$. The following lemma tells us how to geometrically determine the coefficients of $(f_1\cdots f_n)_\circ$ in terms of the coefficients of the tropical hypersurfaces $\tV_i$.
\begin{lm}
\label{lm:coefficients of union}
Let $\tV_1,\ldots,\tV_n$ be $n$ enriched tropical hypersurfaces in $\R^n$ with enriched Viro polynomials 
\[f_i(x)=\sum_{I^j\in A_j}\alpha_{I^j}x^{I^j}t^{\varphi_j(I^j)}\]
for $i=1,\ldots,n$. Then the coefficients of the enriched tropical hypersurface given by the union $\tV_1\cup \tV_2\cup \ldots\cup \tV_n$ %are the coefficients of of the tropical hypersurface associated to 
%the product $f_1\cdots f_n$ of enriched Viro polynomials and 
can be determined as follows.
Let $K$ be a vertex in the dual subdivision of $\tV_1\cup \tV_2\cup \ldots\cup \tV_n$. Let $J^i$ be the vertex in the dual subdivision of $\tV_i$ such that the connected component dual to $K$ in $\R^n \setminus (\tV_1\cup\ldots\tV_n)$ is a subset of the connected component in $\R^n\setminus \tV_i$ dual to $J^i$ for each $i=1,\ldots,n$. Then the coefficient of the vertex $K$ of $\tV_1\cup \tV_2\cup \ldots\cup \tV_n$ equals $\prod_{i=1}^n\alpha_{J^i}$.
\end{lm}

\begin{ex}
Figure \ref{fig:coeffprod} illustrates how to assign the coefficients to a union of enriched tropical curves.
    \begin{figure}[t]\begin{tabular}{ccc}
\begin{tikzpicture}
\clip (-2,-2) rectangle (2,2);
\draw[line width=2pt,color=red!10!] (-2,0)--(2,0);
\draw[line width=2pt,color=blue] (0,-2)--(0,2);
\node[color=blue] at (-1,0.3) {\Large$\alpha_I$};
\node[color=blue] at (1,0.3) {\Large$\alpha_J$};
\node[color=blue] at (0.4,-1.8) {$C_1$};
\end{tikzpicture}&
\begin{tikzpicture}
\clip (-2,-2) rectangle (2,2);
\draw[line width=2pt,color=blue!10!] (0,-2)--(0,2);
\draw[line width=2pt,color=red] (-2,0)--(2,0);
\node[color=red] at (0.4,1) {\Large$\beta_{I'}$};
\node[color=red] at (0.4,-1) {\Large$\beta_{J'}$};
\node[color=red] at (1.8,0.3) {$C_2$};
\end{tikzpicture}&
\begin{tikzpicture}
\clip (-2,-2) rectangle (2,2);
\draw[line width=2pt,color=blue] (0,-2)--(0,2);
\draw[line width=2pt,color=red] (-2,0)--(2,0);
\node[color=purple] at (-1,1) {\Large$\alpha_I\beta_{I'}$};
\node[color=purple] at (-1,-1) {\Large$\alpha_{I}\beta_{J'}$};
\node[color=purple] at (1,1) {\Large$\alpha_{J}\beta_{I'}$};
\node[color=purple] at (1,-1) {\Large$\alpha_{J}\beta_{J'}$};
\node[color=blue] at (0.4,-1.8) {$C_1$};
\node[color=red] at (1.8,0.3) {$C_2$};
\draw[color=purple,fill=purple] (0,0) circle (1.4pt);
\end{tikzpicture}\\
\end{tabular}
\caption{Coefficients of the union $\tC_1\cup \tC_2$ around an intersection point.\label{fig:coeffprod}}
\end{figure}
\end{ex}

%{\color{blue} are we sure off all the indices, inequalities here?}
\begin{proof}[Proof of Lemma \ref{lm:coefficients of union}]
Let $p=(p_1,\ldots,p_n)$ be a point in the interior of the connected component in the complement $\R^n\setminus (\tV_1\cup\ldots\cup\tV_n)$ that is dual to $K$.
Then $p$ is in the interior of the connected component dual to $J^i$ in $\R^n\setminus \tV_i$ for $i=1,\ldots,n$.
That means that 
\[\sum_{j=1}^nJ^i_j\cdot p_j-\varphi_i(J^i)>\sum_{j=1}^nI^i_j\cdot p_j-\varphi_i(I^i)\]
for any $I^i\neq J^i$ and all $i=1,\ldots,n$. 
Hence,
\[\sum_{i=1}^n\left(\sum_{j=1}^nJ^i_j\cdot p_j\right)-\sum_{i=1}^n\varphi_i(J^i)>\sum_{i=1}^n\left(\sum_{j=1}^nI^i_j\cdot p_j\right)-\sum_{i=1}^n\varphi_i(I^i)\]
for any $(I^1,\ldots,I^n)\neq (J^1,\ldots,J^n)$. 
In particular for any $(I^1,\ldots,I^n)\neq (J^1,\ldots,J^n)$ such that $\sum_{i=1}^nI^i= \sum_{i=1}^nJ^i=K$ we get that
\[\sum_{i=1}^n\varphi_i(J^i)<\sum_{i=1}^n\varphi_i(I^i)\]
and hence the monomial with exponent $K$ in $f_1\cdots f_n$ equals
\[\left(\prod_{i=1}^n \alpha_{J^i}t^{\sum_{i=1}^n\varphi_i(J^i)}+\text{higher order terms in }t \right)\cdot x^K.\]
%Hence, the monomial with exponent $K$ in $f_1\cdot f_2$ equals
%\[\left(\alpha_I\beta_Jt^{\varphi(I)+\psi(J)}+\text{higher order terms in } t\right)\cdot x^{k_1}y^{k_2}\]
%and thus, the coefficient of $\tC_1\cup \tC_2$ at the vertex $K$ is $\alpha_I\beta_J$.
%Since in $k\Puiseux^\times/(k\Puiseux^\times)^2$ we have 
%\[\left(\prod_{i=1}^n \alpha_{J^i}t^{\sum_{i=1}^n\varphi_i(J^i)}+\text{higher order terms in }t \right)=\prod_{i=1}^n \alpha_{J^i}\]
and thus, the coefficient of $K$ equals $\prod_{i=1}^n\alpha_{J^i}$.
\end{proof}

%%%%%%%%%%%%%%%%%%%%%%%%%%%%%%%%%%%%%%%%%%%%%%%%%
\subsection{Combinatorics of tropical curves}
We redefine the concept of a tropical curve in $\R^n$ from a combinatorial point of view. We will need this in the proof of Proposition \ref{prop:number of parallelepipeds}. The following definition coincides with the definition of a tropical curve in $\R^2$ defined algebraically as before and it is known in the literature as an \emph{embedded abstract tropical curve}.

\begin{df}
A \emph{tropical curve} $C$ is a finite weighted graph $(V,E,\omega)$ embedded in~$\R^n$, where ~$E$ %=E^{\circ}\cup E^{\infty}$ 
is the disjoint union of non-directed edges $E^{\circ}\subset \{e\subset V\mid \Card(e)=2\}$ and univalent edges $E^{\infty}\subset V$, such that every edge $e\in E^{\circ}$ embeds into a segment of the graph of an integer line, i.e. given by $\bar{a}_e\cdot t+\bar{b}_e$ with $\bar{a}_e\in\Z^n\setminus\{\bar{0}\},\bar{b}_e\in\Q^n$, every edge $l\in E^{\infty}$ embeds into a ray of an integer line, and every vertex $v\in V$ satisfies
the \emph{balancing condition}
\[
\sum_{e\in E, v\in e} \omega(e)\cdot \mathrm{u}_e=0
\]
where $\mathrm{u}_e=\displaystyle\frac{\pm 1}{\gcd(\bar{a}_e)}\bar{a}_e$ oriented outwards from $v$, and $\omega:E\lra\Z$ is a non-negative function.
We call $\bar{a}_e$ a \emph{direction} vector of $e$ and $\mathrm{u}_e$ a \emph{primitive} vector of $e$ at $v$.
When drawing a tropical curve we write the weights not equal to $1$ next to the edges.
\end{df}
\begin{ex}
Figure \ref{figure:balancingCondition} shows a 3-valent vertex $v$ with its three primitive vectors in blue. The one edge labeled $2$ has weight $2$ while the other edges have weight $1$. The balancing condition is satisfied since
\[1\cdot \begin{pmatrix}-1\\0\end{pmatrix}+1\cdot \begin{pmatrix}1\\2\end{pmatrix}+2\cdot \begin{pmatrix}0\\-1\end{pmatrix}=\begin{pmatrix}0\\0\end{pmatrix}.\]
\end{ex}
\begin{figure}[t]
\begin{center}
\begin{tikzpicture}
\draw[step=1cm,gray!20!white,very thin] (-1.9,-1.9) grid (1.9,3.9);
\draw[thick] (-2,0)--(0,0)--(2,4);
\draw[thick] (0,0)--(0,-2);
\node[] at (0.3,-1) {$2$};
\draw[->,blue, thick](0,0) -- (-1,0);
\draw[->,blue, thick](0,0) -- (0,-1);
\draw[->,blue, thick](0,0) -- (1,2);
\node at (0.3,-0.1) {$v$};
\end{tikzpicture}
\caption{Balancing condition}
\label{figure:balancingCondition}
\end{center}
\end{figure}

% Here, the polygon $\Delta_1+\Delta_2$ is the  Minkowski sum of the polygons. 

% B\'ezout's theorem for tropical curves is the following special instance.
% Assume that the degrees~$\Delta_1=\Delta_{d_1}=\operatorname{Conv}(\{(0,0),(d_1,0),(0,d_1)\})$ and $\Delta_2=\Delta_{d_2}=\operatorname{Conv}(\{(0,0),(d_2,0),(0,d_2)\})$, i.e, the tropical curves $C_1$ and $C_2$ correspond to curves in the projective plane. Since
% \[\Delta_{d_1}+\Delta_{d_2}=\operatorname{Conv}(\{(0,0),(d_1+d_2,0),(0,d_1+d_2)\})=\Delta_{d_1+d_2},\]
% then we have that in this case
% \[
% C_1\cdot C_2= \operatorname{Area}(\Delta_{d_1+d_2})-\operatorname{Area}(\Delta_{d_1})-\operatorname{Area}(\Delta_{d_2})=\frac{(d_1+d_2)^2}{2}-\frac{d_1}{2}-\frac{d_2^2}{2}=d_1\cdot d_2.
% \]

%%%%%%%%%%%%%%%%%%%%%%%%%%%%%%%%%%%%%%%%%%%%%%%%%

%%%%%%%%%%%%%%%%%%%%%%%%%%%%%%%%%%%%%%%%%%%%%%%%%

%% file: 5higherdimensions.tex
\section{Enriched tropical intersection multiplicities}
\label{section:Higher dimensions}

%We do an analog of the preceding section in higher dimensions.
In this section we define the \emph{enriched intersection multiplicity} $\mult_p(\tV_1,\dots,\tV_n)$ of $n$ tropical hypersurfaces in $\R^n$ at an intersection point $p$ in terms of their defining polynomials. We prove that it can be computed by a purely combinatorial formula.
 %{\color{red} do we need to define what the Newton polytopes of $f_i$ are}

\subsection{Notation and conventions}
We start by establishing notation for this section.
The \emph{Newton polytope} of a Laurent polynomial $f=\sum a_Ix^I\in k\Puiseux[x_1,\ldots,x_n]$ is the convex hull $\operatorname{Conv}\{I: a_I\neq 0\}\subset \R^n$. Note that this agrees with the Newton polytope of the associated tropical polynomial.
In the whole section we work over a field $k$ of characteristic $0$ or characteristic bigger than the diameter $\operatorname{diam}(NP(f_i))$ of the Newton polytopes of the polynomials $f_i\in k\Puiseux[x_1,\ldots,x_n]$.
%This notation does not agree with the notation for curves, but both are natural notations for the cases they address.
We use $\tV$ for enriched tropical hypersurfaces. We write $V$ for the underlying (non-enriched) tropical hypersuface.
Let $\tV$ be an enriched tropical hypersurface in $\R^n$ and let
\[\sum_{I\in A} \alpha_{I}x^It^{\phi(I)}\]
be an enriched Viro polynomial for $\tV$. 
Here, the sum runs over a finite set $A$ of integer $n$-tuples $I=(I_1,\ldots,I_n)$ in $\Z^n$, the map $\phi:A\longrightarrow \Q$ is the restriction of a rational convex function to $A$, and $x^I=x_1^{I_1}\cdots x_n^{I_n}$. The coefficients $\alpha_I$ are elements of $k^{\times}$.
Now assume we have $n$ tropical hypersurfaces $\tV_1,\ldots,\tV_n$ in $\R^n$ with enriched Viro polynomials
\[f_i=\sum_{I^{i}\in A_i}\alpha_{I^{i}}x^{I^i}t^{\phi_i(I^i)}\]
for $i=1,\ldots,n$ and assume that $\tV_1,\ldots,\tV_n$ intersect tropically transversally at $p$.
Then for each~$i\in \{1,\ldots,n\}$, the point $p$ lies on a top dimensional face of $\tV_i$ separating two connected components of $\R^n$ (see Figure \ref{fig:coeffprod}). Assume these components are the components where $I_1^ix_1+\ldots+I_n^ix_n-\phi_i(I^i)$ and
$J_1^ix_1+\ldots+J_n^ix_n-\phi_i(J^i)$ attain the maximum for some $I^i,J^i\in A_i$.
Then we say that for $\alpha_i=\alpha_{I^i}$ and $\beta_j=\beta_{J^j}$
\[f_i^{\operatorname{local}}\coloneqq\alpha_ix^{I^i}+\beta_ix^{J^i}\]
is the \emph{local binomial equation} of $\tV_i$ at $p$.
Put $\Delta^i\coloneqq I^i-J^i$ (the order of $I,J$ matters) and
let 
\[M=(\Delta^1,\ldots,\Delta^n)\]
be the matrix with columns the $\Delta^i$ and
\[m=\vert \det M \vert\]
be the absolute value of its determinant.
%and write $\det (\Delta^i)$ for the determinant of the matrix with rows $\Delta^1,\ldots,\Delta^n$. {\color{red} rows vs columns}
%

For an intersection point~$p$ of $\tV_1,\ldots,\tV_n$, we write $P$ for the parallelepiped in the dual subdivision of $\tV_1\cup\ldots\cup \tV_n$ dual to $p$.
The parallelepiped $P$ is given by
\[
P=\operatorname{Conv}\left\{ \sum_{i=1}^n K^i \middle| K^i=I^i \textrm{ or } K^i=J^i \right\}= \operatorname{Conv}\left\{ v_A=v_\circ-\sum_{i\in A}\Delta^i\middle|A\subset \{1,\ldots,n\} \right\}
\]
with $v_\circ=\sum_{i=1}^n I^i$.
%
%For a vertex $v=\sum_{i=1}^nI^i$ for some $I^i\in A_i$ for $i=1,\ldots,n$ in the dual subdivision   of $\tV_1\cup\ldots\cup \tV_n$ we write $\alpha_v\coloneqq \prod_{i=1}^n\alpha_{I^i}$ for its coefficient in .
%
In particular, the non-enriched tropical intersection multiplicity at~$p$ (see Definition \ref{df:nonenrichedintmult}) equals 
\[\operatorname{mult}_p(V_1,\ldots,V_n)\overset{\ref{df:nonenrichedintmult}}{=}\operatorname{Vol}(P)=\vert \det M\vert=m.\]

%For $z\in \widebar{k}\Puiseux^n$ write $z_\circ\coloneqq(z_1,\ldots,z_n)= \Ini z\in \bar{k}^n$ for the initial of $z$. 
%
Furthermore, we let $e=(1,\ldots,1)\in \Z^n$ and $e_i=(0,\ldots,0,1,0,\ldots,0)$ with $1$ in the $i$th position. 
%We set the vector $\hat{e}=(1,\ldots,1)\in \Z^{n-1}$ and write $\hat{x}=(x_1,\ldots,x_{n-1})$.
%For $I\in \Z^n$, we put $\hat{I}\coloneqq (I_1,\ldots,I_{n-1})\in \Z^{n-1}$.
%{\color{blue} shall we delete this last lines? I re-introduce them when we use them anyways}

\subsection{Enriched tropical intersection multiplicity}

We want to define the enriched intersection multiplicity to agree with the local index, as defined in Definition \ref{df:local index}, at the zero of the section of $V\coloneqq\glob(d_1)\oplus\glob(d_2)\oplus\ldots\oplus \glob(d_n)\longrightarrow \mathbb{P}^n_{k\Puiseux}$ defined by the $n$ hypersurfaces in $\mathbb{P}^n_{k\Puiseux}$. Assume our tropical hypersurfaces have enriched Viro polynomials~$f_1,\ldots,f_n$ in $ k\Puiseux[x_1,\ldots,x_n]$.
We have seen that the local index in the Poincar\'e-Hopf theorem for the vector bundle $V$ equals the local $\A^1$-degree of $(f_1,\ldots,f_n)\colon \A^n_{k\Puiseux}\longrightarrow \A^n_{k\Puiseux}$
for which we have an explicit formula, namely the trace of the determinant of the Jacobian evaluated at the zero (see \eqref{eq:nonrational}). Note that to $f\in k\Puiseux[x_1,\ldots,x_n]$ we can associate an enriched Viro polynomial $f_{\circ}$ by replacing the coefficients of $f$ by its initials. We say that $f$ is a \emph{defining polynomial} for the enriched tropical hypersurface with enriched Viro polynomial $f_\circ$.
This motivates the following definition.

For the following definition of the enriched tropical intersection multiplicity we will assume that the underlying tropical hypersurfaces intersect tropically transversally at an intersection point $p$.
Lemma \ref{lm:tropintmult} also holds for tropical non-transverse intersections. 
However, in the enriched setting, the local contribution given by a small perturbation is not invariant and hence we restraint our consideration to the transversal case.

%{\color{blue} I want a big characteristic assumption and a tropically transverse intersection at $p$ and we need to explain why we need these assumptions!}
\begin{df}
\label{def:enrichedMult}
Let $\tV_1,\ldots,\tV_n$ be $n$ tropical hypersurfaces in $\R^n$ with defining polynomials $f_1,\ldots,f_n\in k\Puiseux[x_1,\ldots,x_n]$. Let $p$ be an intersection point $p\in V_1\cap\ldots\cap V_n$ at which the $V_i$ interesect tropically transversally. Assume that there are exactly $s$ closed points $z^{(1)},\ldots,z^{(s)}$ in $\A^n_{k\Puiseux}$ such that $-\operatorname{val}(z^{(i)})=p$.
%let $z$ be a closed point in $\A^n_{k\Puiseux}$ such that $-\operatorname{val}(z)=p$ and such that $z$ is a zero of $f_1=\ldots =f_n=0$. {\color{blue} I don't like the definition of $z$ is should somehow be the product of all these $z$'s}
We define the \emph{enriched intersection multiplicity} of $\tV_1,\ldots,\tV_n$ at $p$ to be 
\begin{equation}
\label{eq:defof multp}
\mult_p(\tV_1,\ldots,\tV_n)\coloneqq \sum_{i=1}^s\Tr_{\kappa(z^{(i)})/k\Puiseux}\qqinv{\det \operatorname{Jac}(f_1,\ldots,f_n)(z^{(i)})}\in \GW(k\Puiseux),
\end{equation}
where %the sum ranges over zeros $z$ of~$f_1,\ldots,f_n$ with minus the valuation equal to $p$ and 
$\kappa(z^{(i)})$ is the residue field of $z^{(i)}$.
\end{df}
%\begin{rmk}
By the assumption on the characteristic of $k$, we have that $\kappa(z^{(i)})= L_i\Puiseux$ for some finite separable field extension $L_i$ of $k$ for $i=1,\ldots,s$. In particular, we can collect all $\kappa(z^{(i)})$ in a finite \'etale algebra
\[E_t\coloneqq \prod_{i=1}^s\kappa(z^{(i)})=L_1\Puiseux\times\ldots \times L_q\Puiseux.\]
%where the product is over all $z$ with $-\operatorname{val}(z)=p$.
Since $\Tr_{E_t/k\Puiseux}=\sum_i \Tr_{L_i\Puiseux/k\Puiseux}$ we have
\begin{equation}
\label{eq:enrichedmultwithfiniteetalealgebra}
\mult_p(\tV_1,\ldots,\tV_n)=\Tr_{E_t/k\Puiseux}(\qinv{\det \operatorname{Jac}(f_1,\ldots,f_n)(z)})\in \GW(k\Puiseux)\end{equation}
where $z=(z^{(1)},\ldots,z^{(s)})\in E_t$.

%$E_t\cong F_1\Puiseux\times\ldots \times F_r\Puiseux$ where the $F_i$ are finite separable field extensions of $k$ and $F_i\Puiseux$ are the residue fields of the common zeros of the $f_i$. Since $\Tr_{E_t/k\Puiseux}=\sum_i \Tr_{F_i\Puiseux/k}$ we have that \eqref{eq:defof multp} agrees with \eqref{eq:enrichedIMDef}.
%\end{rmk}
%Recall from Example \ref{ex:GW of Puiseux} that $\GW(k\Puiseux)\cong \GW(k)$ and that the isomorphism sends a generator $\qinv{\sum_{i=i_0}^\infty a_it^{\frac{i}{n}}}\in \GW(k\Puiseux)$ to its non-zero coefficient of lowest degree $\qinv{ a_{i_0}} \in \GW(k)$. Hence, it suffices to find the non-zero coefficient of lowest degree of $\det \Jac(f_1,\ldots,f_n)(z)$. 
%Recall that we write $z_\circ=(z_1,\ldots,z_n)$ for the closed point in $\A^n_k$ defined by the initials of~$z$. Further, let $E$ be the coordinate ring of all $z_\circ$'s such that $z$ is a zero of $f_1,\ldots,f_n$ and such that~$-\operatorname{val}(z)=p=(p_1,\ldots,p_n)$.

Let
\begin{equation}\label{eq:locbinhyp}
    f_i^{\operatorname{local}}= \alpha_iz^{I^i}+\beta_iz^{J^i}
\end{equation}
be the local binomial equation of $\tV_i$ at $p$, for $i=1,\ldots,n$.
Set 
\[E=\frac{k[x_1^{\pm1},\ldots,x_n^{\pm1}]}{(\alpha_ix^{\Delta^i}+\beta_i)_{i=1,\ldots,n}}\cong L_1\times\ldots\times L_s.\]
 Using the isomorphism $\GW(L_i\Puiseux)\cong \GW(L_i)$ from Example \ref{ex:GW of Puiseux} we identify \eqref{eq:enrichedmultwithfiniteetalealgebra}
with 
\begin{equation}
\label{eq:enrichedmultwithfiniteetalealgebrawithoutPuiseux}
\mult_p(\tV_1,\ldots,\tV_n)=\Tr_{E/k}\qqinv{\det \operatorname{Jac}(f_1^{\operatorname{local}},\ldots,f_n^{\operatorname{local}})(z_\circ)}\in\GW(k)\end{equation}
where $z_\circ=(z_1,\ldots,z_n)\in E^n$ defined by the initials of~$z$.

%Write $z_i=(z_\circ)_it^{-p_i}+\text{higher order terms}$ for $i=1,\ldots,n$. Then 
%\[f_i(z)=f_i(z_1,\ldots,z_n)=\alpha_iz_\circ^It^{\phi_i(I^i)-\sum_{j=1}^np_jI_j^i}+\beta_iz_\circ^It^{\phi_i(J^i\sum_{j=1}^np_jI_j^i)}+\text{higher order terms}\]
%We know that 
%\[\phi_i(I^i)-\sum_{j=1}^np_jI_j^i=\phi_i(J^i)-\sum_{j=1}^np_jJ_j^i.\]

%{\color{red} for tropically transverse intersections $\alpha_i,\beta_i\in k^\times$?}
%The following Theorem is the $n$-dimensional analog of Theorem \ref{thm:enriched intersection multiplicity} for curves.%and recall that to find $z_\circ$ we need to solve for $f_i^{\operatorname{local}}=0$ for all $i=1,\ldots,n$.
%{\color{blue} for big characteristic or characteristic 0 $E_t=E\Puiseux$}

\begin{prop}
\label{prop:enriched intersection multiplicity dim n}
With the notation from Definition \ref{def:enrichedMult},
if $\tV_1,\ldots,\tV_n$ intersect tropically transversely at~$p$, then
\[\qinv{\det\operatorname{Jac}(f_1^{\operatorname{local}},\ldots,f_n^{\operatorname{local}})(z)}=
\qinv{\det M \cdot \prod_{i=1}^n\alpha_{i} \cdot z_{\circ}^{\sum_{i=1}^nI^i-e}}\in \GW(E),\]
where $e=(1,1,\ldots,1)\in \Z^n$.
In particular, the enriched intersection multiplicity at $p$ equals
\[
\mult_p(\tV_1,\ldots,\tV_n)=
\Tra[E/k]{
\qinv{ 
\det M\cdot \prod_{i=1}^n\alpha_{i} z_{\circ}^{\sum_{i=1}^nI^i-e} }}\in\GW(k).
\] 
%where $E$ is the coordinate ring of all such $z_\circ$.
\end{prop}

\begin{proof}

%By Example \ref{ex:GW of Puiseux} it suffices to compute
%$\det \frac{\partial f_i^{\operatorname{local}}}{\partial x_j}(z_\circ)$. Using Equations~\eqref{eq:locbinhyp} 
We calculate that 

\begin{align*}
    \det \frac{\partial f_i^{\operatorname{local}}}{\partial x_j}(z_\circ)_{i,j}&= \det(I^i_j\alpha_iz_\circ^{I^i-e_j}+J^i_j\beta_iz_\circ^{J^i-e_j})_{i,j}\\
     &=\det(I^i_j\alpha_iz_\circ^{I^i-e_j}-J^i_j\alpha_iz_{\circ}^{I^i-e_j})_{i,j}\\
    &=\prod_{i=1}^n\alpha_i\cdot z_\circ^{\sum_{i=1}^nI^i-e}\cdot \det (I^i_j-J^i_j)_{i,j}.
\end{align*}
\end{proof}

The following Lemma shows that the enriched intersection multiplicity as calculated in Proposition~\ref{prop:enriched intersection multiplicity dim n} is independent of the order of the exponent vectors $I^j$ and $J^j$, for $j=1,\ldots,n$.

\begin{lm}
\label{lm: symmetry dim n}
The determinant of the Jacobian
\begin{equation}\label{eq:detjac}
    \det M \cdot \prod_{i=1}^n\alpha_{i} \cdot  z_{\circ}^{\sum_{i=1}^nI^i-e}
\end{equation}
is invariant under the exchange of the roles of $I^j$ and $J^j$ and simulaneously the roles of $\alpha_j$ and $\beta_j$ for any $j\in \{1,\ldots,n\}$.
\end{lm}
\begin{proof}
Recall that we have 
\[f_i^{\operatorname{local}}(z_\circ)=\alpha_{i}z_\circ^{I^i}+\beta_iz_\circ^{J^i}=0\]
for $i=1,\ldots,n$, or equivalently,
\begin{equation}
\label{eq:replaceIbyJ}
-\frac{\beta_i}{\alpha_{i}}z_{\circ}^{J^i-I^i}=1.
\end{equation}

If we exchange the roles of $I^j$ and $J^j$, we replace $\Delta^j=J^j-I^j$ by $-\Delta^j=I^j-J^j$ and then~$\det M \cdot \prod_{i=1}^n\alpha_{i} \cdot  z_{\circ}^{\sum_{i=1}^nI^n-e}$ becomes 
\begin{align*}
    &(-\det M )\cdot \prod_{\substack{i=1\\i\neq j}}^n\alpha_{i}\cdot \beta_j \cdot  z_{\circ}^{\sum_{i=1}^nI^i-e+J^j-I^j}\\
    =&\left(-\frac{\beta_j}{\alpha_{j}}z_\circ^{J^j-I^j}\right)\cdot\det M \cdot \prod_{i=1}^n\alpha_{i}  \cdot  z_{\circ}^{\sum_{i=1}^nI^i-e}\\
    \overset{\eqref{eq:replaceIbyJ}}{=}&\det M \cdot \prod_{i=1}^n\alpha_{i} \cdot  z_{\circ}^{\sum_{i=1}^nI^i-e}.
\end{align*}
\end{proof}

\subsection{A combinatorial formula for $\mult_p(\tV_1,\ldots,\tV_n)$}

%Recall that an intersection of $n$ tropical hypersurfaces corresponds to a parallelepiped in the dual subdivision of $V_1\cup\ldots\cup V_n$ and the volume of this parallelepiped equals the classical tropical intersection multiplicity by Lemma \ref{lm:tropintmult}. Theorem \ref{thm:enriched intersection multiplicity dim n} assigns an enriched intersection multiplicity $\mult_p(\tV_1,\ldots,\tV_n)$ to each intersection point $p$.
We identify the formula for the enriched intersection multiplicity from Proposition~\ref{prop:enriched intersection multiplicity dim n} with an element of $\GW(k)$ which can be read of the dual subdivision of~$\tV_1\cup\ldots\cup \tV_n$. We will see that the intersection multiplicity we computed in Proposition \ref{prop:enriched intersection multiplicity dim n} is determined by coefficients at the ``odd vertices" in the dual subdivision of $\tV_1\cup\ldots\cup\tV_n$ in the sense of the following definition.
\begin{df}
We call a lattice point $v=(v_1,v_2,\ldots,v_n)\in\Z^n$ \emph{odd}, if its class in $(\Z/2\Z)^n$ equals $(1,1,\ldots,1)$. 
\end{df}
 %in the dual subdivision of $\tV_1\cup\ldots\cup\tV_n$ in arbitrary dimension.

Let $p$ be an intersection point of enriched tropical hypersurfaces $\tV_1,\ldots,\tV_n$ and let $P$ be the parallelepiped in the dual subdivision of $\tV_1\cup\ldots\cup\tV_n$ dual to $p$.
%Furthermore, we assign a sign to the vertices $v$: 
Recall that the local binomials equation of $\tV_i$ at $p$ are of the form 
\[f_i^{\operatorname{local}}=\alpha_{i}x^{I^i}+\beta_ix^{J^i}\]
with $\alpha_i,\beta_i\in k^\times$ and 
\[\Delta^i=I^i-J^i\]
for $i=1,\ldots,n$.
Let $v$ be a corner vertex of $P$. Then $v$ can be uniquely expressed as 
\[v= v_\circ -\sum_{i\in A}\Delta^i\eqqcolon v_A\]%\sum_{i=1}^n I^i-\delta_{K^i,J^i}\Delta^i,\]
for some subset $A\subset \{1,\ldots,n\}$, where $v_\circ =\sum_{i=1}^n I^i$. 
The coefficient of $v_A$ equals $$\alpha_{v_A}=\prod_{i\not \in A}\alpha_i \cdot \prod_{i\in A}\beta_i$$ by Lemma \ref{lm:coefficients of union}.
%$\delta_{K^i,J^i}$ is the Kronecker delta. 
Furthremore, we define the \emph{sign} of the vertex $v_A$ with respect to the parallelepiped~$P$ as  
\begin{equation}
\label{eq:signdimn}
\epsilon_P(v)\coloneqq (-1)^{\# A}\cdot \operatorname{sign}(\det M).
\end{equation}
and say that
\[\gamma_{v_A}\coloneqq \epsilon_P(v)\cdot \alpha_{v_A}\]
is the \emph{signed coefficient} of $v_A$.
\begin{ex}
\label{ex:epsilons}
Figure \ref{figure:EnrichedIntersection} shows an intersection of two tropical curves $C_1$ and $C_2$ and the dual parallelogram. %On the right is an example of the assignment of signs. 
%\begin{figure}[t]
%\begin{tabular}{ccc}
%\includegraphics{Figures/sign0.pdf}&
%\includegraphics{Figures/sign1.pdf}&
%\includegraphics{Figures/sign2.pdf}
%\end{tabular}
%\caption{The dual parallelogram associated to an intersection point and a local picture of signed vertices. {\color{red}Do we still need this figure?} \label{fig:sign}} 
%\end{figure}
\begin{figure}
     \centering
     \begin{subfigure}[b]{0.45\textwidth}
         \centering
         \begin{tikzpicture}[scale=0.7]
\draw[red,very thick,->] (-3,-2)--(3,2);
\draw[blue,very thick,->] (0,-3)--(0,3);
\filldraw[black] (0,0) circle (2pt) node[anchor=east]{$p$};
\node at (-1,1.5) {$\alpha_1\beta_2$};
\node at (1,2) {$\beta_1\beta_2$};
\node at (-1,-2) {$\alpha_1\alpha_2$};
\node at (1,-1.5) {$\beta_1\alpha_2$};
\node[blue] at (0.5,-2.5) {$C_1$};
\node[red] at (1.5,0.5) {$C_2$};
\end{tikzpicture}
\caption{Intersection at $p$}
     \end{subfigure}
     \hfill
     \begin{subfigure}[b]{0.45\textwidth}
         \centering
         \begin{tikzpicture}
\fill[gray!10!white,thick] (-2,0)--(0,0)--(2,-3)--(0,-3)--cycle;
\draw[blue,thick] (-2,0)--(0,0);
\draw[blue,thick] (2,-3)--(0,-3);
\draw[red,thick] (0,0)--(2,-3);
\draw[red,thick] (0,-3)--(-2,0);
\node at (0,-1.5) {$P$};
\node at (0.5,0.2) {$v_{\{1,2\}}$};
\node at (-2.4,0.2) {$v_{\{2\}}$};
\node at (-0.4,-3.2) {$v_\circ=v_{\emptyset}$};
\node at (2.5,-3.2) {$v_{\{1\}}$};
\node[blue] at (1,-3.3) {$\Delta I^1$};
\node[red] at (-1.5,-1.5) {$\Delta I^2$};
\end{tikzpicture} 
\caption{Dual parallelogram $P$}
     \end{subfigure}
        \caption{An enriched tropical intersection together with the dual parallelogram.}
         \label{figure:EnrichedIntersection}
\end{figure}
%{\color{blue} add more details to the example and maybe mention the rule for curves, we might want to add $I$, $I'$, $J$, $J'$ to the picture and explain in detail how to assign the signs for curves}
The determinant $\det M$ records the order in which the curves intersect in the following sense: Choosing the roles of $I^i$ and $J^i$ orients the edges of $C_i$ such that $I^i$ is on the left and $J^i$ is on the right for $i=1,2$. We can always arrange the ordering of $I^i$ and $J^i$ so that $\det M=+1$. Swapping $C_1$ and $C_2$ changes the sign of $\det M$ for all intersections of $C_1$ and $C_2$. 
For the four vertices $v_\circ=v_\emptyset=I^1+I^2$, $v_{\{2\}}=I^1+J^2$, $v_{\{1\}}=J^1+I^2$ and $v_{\{1,2\}}=J^1+J^2$ we get
\[\epsilon_P(v_\circ)=+1, \; \epsilon_P(v_{\{1\}})=-1,\; \epsilon_P(v_{\{2\}})=-1,\; \epsilon_P(v_{\{1,2\}})=+1 \]
given that $\sign (\det M)=+1$.
%$(-1)\cdot(-1)$, $(-1)\cdot (+1)$, $(+1)\cdot(-1)$ and $(+1)\cdot (+1)$ respectively and it can be determined in the following way.
This sign can be determined in the following way.
 When you walk around the vertex $v$ inside $P$ anticlockwise and the edge you start at is dual to an edge of $C_1$ then $\epsilon_P(v)=+1$, if it is dual to an edge of $C_2$ then $\epsilon_P(v)=-1$.
 The signed coefficients of the vertices equal 
 \[\gamma_{v_\circ}=\alpha_1\alpha_2,\; \gamma_{v_{\{2\}}}=-\alpha_1\beta_2,\;  \gamma_{v_{\{1\}}}=-\beta_1\alpha_2,\; \gamma_{v_{\{1,2\}}}=\beta_1\beta_2.\]
\end{ex}

Geometrically, the sign of the vertex $v$ with respect to the parallelepiped~$P$ is the sign of the determinant of the edges of the parallelepiped~$P$ adjacent to $v$ oriented outwards from~$v$. Namely, if for every $i=1,\dots,n$ we put $\epsilon_i$ as the sign $\pm1$ such that $v+\epsilon_i\Delta^i\in P$, then $\epsilon_P(v)=\operatorname{sign}(\det \epsilon_i\Delta^i)= \prod_{i=1}^n\epsilon_i\cdot \sign(\det M)$ where $ \det\epsilon_i\Delta^i$ is the determinant of the matrix with columns $\epsilon_i\Delta^i$.  
In particular, we have that
\begin{equation}
\label{eq:epsilonsgeometrically}
\epsilon_P(v_\circ)\cdot m=\det M\end{equation}
where $m=\vert \det M\vert $.
%and the sign, as defined here. %, agrees with the definition of the sign \eqref{eq:sign} for curves.

\begin{rmk}
\label{rmk:signs dim n}
Note that the sign of the vertex $\sum_{i=1}^nI^i$ is the  opposite of the sign of $\sum_{i=1}^nI^i+J^j-I^j$ for any $j\in \{1,\ldots,n\}$.
For example in the case of curves, we have that the sign of the vertex $v_\circ=I^1+I^2$ is the same as the sign of $v_{\{1,2\}}=J^1+J^2$ and opposite of the sign of $v_{\{2\}}=I^1+J^2$ and $v_{\{1\}}=J^1+I^2$. %(compare with Remark \ref{rmk:signs}). 
\end{rmk}

%Let $q\in \{0,1,\ldots, 2^n\}$ be the number of odd corner vertices $v_1,\ldots,v_{q}$ of $P$. For $l=1,\ldots, q$ let $\alpha_{v_l}$ be the coefficient of $v_l$. 
\begin{thm}[Main theorem]
\label{thm:mainthm}
Assume that $k$ is a field of characteristic $0$ or characteristic bigger than the diameter of the Newton polytopes of the enriched Viro polynomials $f_i$ of $\tV_i$ for $i=1,\ldots,n$. Further, assume $\operatorname{char}k\neq 2$.
Let $p$ be an intersection point of enriched tropical hypersurfaces~$\tV_1,\ldots,\tV_n$ that intersect tropically transversally at $p$. Let $P$ be the parallelepiped in the dual subdivision of $\tV_1,\ldots,\tV_n$ corresponding to $p$ and let $v_1,\ldots,v_q$ be the odd corner vertices of $P$. 
Assume the non-enriched tropical intersection multiplicity $\operatorname{mult}_p(V_1,\ldots,V_n)$ equals $m$, then
\[\mult_p(\tV_1,\ldots,\tV_n)=\sum_{l=1}^q\qinv{ \epsilon_P(v_l)\alpha_{v_l}}+\frac{m-q}{2}\;h\in \GW(k),\]
where
$\alpha_{v_l}$ is the coefficient of the odd vertex $v_l$ in the dual subdivision of $\tV_1\cup \ldots\cup \tV_n$, for $l=1,\ldots,q$.
\end{thm}
Before proving the theorem we apply it to an example in dimension $2$.

\begin{ex}
\label{ex:enrichedMult}
Figure \ref{fig:innerdual} shows the dual subdivision of the union of the tropical cubic and the tropical conic from Figure \ref{fig:intcurves}. The parallelograms corresponding to the intersections are highlighted. Let~$\alpha_{(3,1)}$ and $\alpha_{(1,3)}$ be the coefficients of the two odd vertices $(3,1)$ and $(1,3)$ of the dual subdivision. Then the intersection corresponding to the upper left parallelogram has enriched intersection multiplicity of rank $1$ since the area of the upper left parallelogram is $1$. Furthermore, it has one odd corner vertex, namely $(1,3)$ and $\epsilon_P(v)$ at this vertex is~$+1$ (see Example \ref{ex:epsilons} for a rule for determining $\epsilon_P(v)$ when $n=2)$. So its enriched intersection multiplicity is
$\qinv{\alpha_{(1,3)}}$.
The enriched intersection multiplicity of the intersection corresponding to the upper right parallelogram
has rank $2$. The upper right parallelogram has two odd corner vertices $(1,3)$ and $(3,1)$. One computes that both signs $\epsilon_P(v)$ at these vertices are $-1$ using Example~\ref{ex:epsilons} and thus this intersection has enriched intersection multiplicity 
$\qinv{ -\alpha_{(1,3)}} +\qinv{-\alpha_{(3,1)}}$.
The remaining intersection point is dual to the lower left parallelogram which has area $3$, and hence the rank of its enriched intersection multiplicity must be $3$. The lower left parallelogram has only one odd corner vertex $(3,1)$, where the sign $\epsilon_P(v)$ equals $+1$. Note that there is an odd point in the interior of this parallelogram but this does not contribute to the enriched intersection multiplicity, only the corner vertices do. To get something of rank $3$, we have to add a hyperbolic form $h$ to $\qinv{\alpha_{(3,1)}}$ to get the enriched intersection multiplicity $h+\qinv{\alpha_{(3,1)}}$ of this last intersection point. 

\begin{figure}[t]
    \centering
    \includegraphics{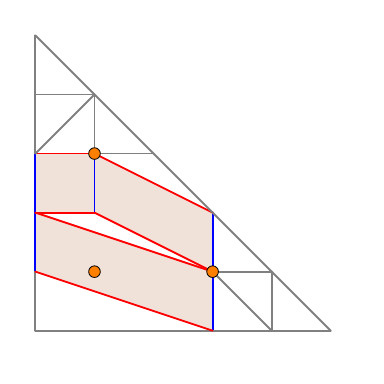}
    \caption{Odd vertices and their adjacent parallelograms in the dual subdivision of the curves in Figure~\ref{fig:intcurves}.}
    \label{fig:innerdual}
\end{figure}
\end{ex}

In order to prove Theorem~\ref{thm:mainthm} we use the Smith normal form of the matrix of exponents in the defining equations of the algebra $E$ to construct an isomorphism to an algebra $E'$ with a diagonal presentation.
We then compute the enriched multiplicity with the presentation of~$E'$ by using Proposition~\ref{prop:traces} and construct an isomorphism in Proposition~\ref{prop:bijection} to relate the obtained computation to the combinatorial data given by the defining equations of $E$.\medbreak

By assumption we have $\vert\det M \vert =m\neq 0$. Hence, there exist matrices $S,T\in \Z^{n\times n}$ both invertible with inverses in $\Z^{n\times n}$ such that
\begin{equation}
\label{eq:smithnormalform}
S M T
% =\begin{pmatrix}
%     m_1& 0 & \ldots &\ldots & 0\\
%     0&m_2& 0& \ldots & 0\\
%     0 & 0 & \ldots &\ldots &0\\
%     0&\ldots & \ldots & \ldots & m_n
% \end{pmatrix}
=\operatorname{diag}(m_1,m_2,\ldots,m_n),\end{equation}
with unique positive integers $m_1,\ldots, m_n$ satisfying $m_i\vert m_{i+1}$ for $i=1,\ldots,n-1$. That is, the matrix $M$ is in \emph{Smith normal form} and the $m_i$ are the \emph{elementary divisors}.
Recall that
\[E=\frac{k[x_1^{\pm1},\ldots,x_n^{\pm1}]}{(\alpha_ix^{I^i}+\beta_ix^{J^i})_{i=1,\ldots,n}}.\]
The $i$th defining equation  of $E$
$\alpha_ix^{I^i}+\beta_ix^{J^i}=0$
is equivalent to
\begin{equation}
\label{eq:lambdai}
    x^{\Delta^i}=-\frac{\beta_{i}}{\alpha_{i}}\eqqcolon \lambda_i
\end{equation}
for $i=1,\ldots,n$.
Since $T$ in \eqref{eq:smithnormalform} is invertible, the ideal generated by the equations in \eqref{eq:lambdai} coincides with the ideal generated by
\begin{equation}
\label{eq:newequations}
    \prod_j x_j^{(M T)_{ji}}=\prod_{j=1}^n \lambda_j^{T_{ji}}\eqqcolon \mu_i, \quad i=1,2,\dots,n.
\end{equation}
%So we get a vector $\lambda =(\lambda_i)_{i=1,\ldots,n}\in k^n$.
%Let $\lambda^T$ be the vector in $k^n$ with $i$th entry 
%\[(\lambda^T)_i=\prod_{j=1}^n \lambda_j^{T_{ji}}.\]
Furthermore, the variable change given by
\[y_i=\prod_{j=1}^n x_j^{(S^{-1})_{ji}}\]
yields an isomorphism of finite \'etale $k$-algebras between $E$ and the finite \'etale $k$-algebra 
\begin{equation}
\label{eq:E'}
E'\coloneqq \frac{k[y_1,\ldots,y_n]}{(y^{m_i}_i-\mu_i)_{i=1,\ldots,n}}\end{equation}
where the $m_i$ are the elementary divisors of $M$ and $\mu_i\in k^\times$ as defined in \eqref{eq:newequations}.
Therefore, for~$m=\vert \det M \vert$,  $\gamma_\circ=\sign(\det M)\cdot \prod_{i=1}^n \alpha_i$ and $v_\circ=\sum_{i=1}^n I^i$, we have that
\begin{align*}
\mult_p(\tV_1,\ldots,\tV_n)&\overset{\ref{prop:enriched intersection multiplicity dim n}}{=}\Tr_{E/k}\qqinv{\det M\cdot \prod_{i=1}^n\alpha_i\cdot x^{\sum_{i=1}^nI^i-e}}\\
&=\Tr_{E/k}\qqinv{m \cdot \gamma_\circ \cdot  x^{v_\circ-e}}\\&=\Tr_{E'/k}\qqinv{m\cdot  \gamma_\circ \cdot y^{
{S(
v_\circ-e)}}}\end{align*}
and the latter can be computed easily using Proposition \ref{prop:traces} as shown in the following proposition for the case when $v_\circ$ is odd.
Note that if $v_\circ$ is odd, then $v_\circ-e$ has only even entries and thus also $S(v_\circ-e)$ has only even entries. So $\qinv{m\cdot \gamma_\circ \cdot y^{S(v_\circ-e)}}=\qinv{m\cdot \gamma_\circ}$ in $\W(E')$.
\begin{prop}\label{prop:traceEp}
Let $E'$ be the finite \'etale $k$-algebra defined in Equation~\eqref{eq:E'}, then
\begin{equation}
    \Tr_{E'/k }\qqinv{m_1\cdots m_n \cdot \gamma_\circ}=
    \sum_{A'} \qinv{\gamma_\circ\cdot\prod_{i\in A'}\mu_i}\in \W(k),
\end{equation}
where the sum runs over all sets $A'\subset\{1,2,\dots,n\}$ such that $m_i$ is even for every $i\in A'$.
\end{prop}

\begin{proof}
    Let 
\[L_i=\frac{k[y_1,\ldots,y_i]}{(y_j^{m_j}-\mu_j)_{j=1,\ldots,i}}.\] 
Then $E'=L_n$, $k=L_0$ and $L_{i+1}%=\frac{L_i[y_{i+1}]}{(y_{i+1}^{d_{i+1}}-(\lambda^T)_{i+1})}
=L_i[y_{i+1}]/(y_{i+1}^{m_{i+1}}-\mu_{i+1})$ for $i=0,\ldots,n-1$. In particular,
\[\Tr_{E'/k }\qqinv{m_1\cdots m_n \cdot \gamma_\circ}=\Tr_{L_1/L_0}\circ\ldots\circ \Tr_{L_n/L_{n-1}}\qqinv{m_1\cdots m_n\cdot \gamma_\circ}.\]
Proposition \ref{prop:traces} implies
\[
\Tr_{L_n/L_{n-1}}\mkern-6mu\qinv{m_1\cdots m_n\cdot \gamma_\circ}=
\begin{cases}
    \qinv{m_1\cdots m_{n-1}\cdot \gamma_\circ} &\text{if $m_i$ is odd},\\
    \qinv{m_1\cdots m_{n-1}\cdot \gamma_\circ}+\qinv{m_1\cdots m_{n-1}\cdot \gamma_{v_\circ}\mu_n} &\text{if $m_i$ is even.}
\end{cases}
\]
Thus, the statement follows from iterating this relation.
% Iterating this, we get
% \[\Tr_{E/k}\qqinv{m_1\cdots m_n \cdot\gamma_\circ}= \sum_{v'\in \Lambda'^{\text{odd}}}\qinv{\gamma_{v'}}\overset{\ref{prop:bijection}}{=}\sum_{v\in \Lambda_E^{\text{odd}}}\qinv{\gamma_v}\overset{\eqref{eq:gammav}}{=}\sum_{v\in \Lambda_E^{\text{odd}}}\qinv{\epsilon_P(v)\alpha_v}\]
% which agrees with what we wanted to show \eqref{eq:multinW(k)}.
\end{proof}

Let $P$ be the parallelepiped dual to $p$ and let $\Lambda_E\subset \Z^n$ be the set of corner vertices of $P$.
For  $v_A=v_\circ-\sum_{i\in A}\Delta^i\in \Lambda_E$ 
for some subset $A\subset \{1,\ldots,n\}$,
let
\[\gamma_{v_A}\coloneqq \epsilon_{P}(v_A)\cdot \prod_{i \not\in A}\alpha_i\cdot \prod_{i\in A}\beta_i =\sign (\det M) \cdot \prod_{i\not \in A}\alpha_i\cdot \prod_{i\in A}(-\beta_i)= \gamma_\circ\cdot \prod_{i\in A}\lambda_i\]
in $k^\times/(k^\times)^2$ be its signed coefficient. 
Here, $\lambda_i=-\frac{\beta_i}{\alpha_i}$ for $i=1,\ldots,n$ as in \eqref{eq:lambdai} and $\gamma_\circ=\sign(\det M)\cdot\prod_{i=1}^n\alpha_i$. %We call $\gamma_{v_A}$ the \emph{signed coefficient} of $v_A$. 
To $E'$ we also associate a parellelepiped, namely the one with corner vertices 
\[\Lambda_{E'}\coloneqq\{v_A'=Sv_\circ-\sum_{i\in A}m_i\cdot e_i\vert A\subset \{1,\ldots,n\}\}\]
where $e_i$ is the $i$th standard basis vector,
and we define the signed coefficients of these corner vertices to be 
\[\gamma_{v_A'}=\gamma_\circ\cdot \prod_{i\in A}\mu_i\in k^\times/(k^\times)^2.\]
To prove our main theorem we will first show that there is a bijection $\phi\colon \Lambda_E\rightarrow \Lambda_{E'}$ such that $\gamma_{v_A}=\gamma_{\phi(v_A)}$ in $k^\times/(k^\times)^2$ for all $A\subset \{1,\ldots,n\}$. %which restrict to a bijection on the odd elements of $\Lambda_E$ and $\Lambda_{E'}$.  

Recall that $T$ %and $S$ 
in \eqref{eq:smithnormalform} is the product of column operations on $M$ where one is allowed to
 swap columns, multiply a column by $-1$, and add a column to another column.
%\begin{itemize}
    % \item swap columns,
    % \item multiply a column by $-1$,
    % \item add a column to another column.
%\end{itemize}
We construct a bijection $\phi$ as described above whenever we perform one of the operations listed above, or any set of row operations.
To show this, assume you start wit a matrix $M=(\Delta^1,\ldots,\Delta^n)$ with columns $\Delta^i$, and with $\lambda_i\in k^\times$ for $i=1,\ldots,n$.
Further, let $M'=(\Delta'^1,\ldots,\Delta'^n)\in \Z^{n\times n}$ be the matrix with columns $\Delta'^i$ for $i=1,\ldots,n$ that one obtains after having applied one of the operations listed above, that is $M'=M\cdot T$ or $M'=S\cdot M$ where $T$ is one of the column operations and $S$ is a set of row operations. For some fixed $v_\circ\in \Z^n$ let 
\[\Lambda=\{v_A=v_\circ-\sum_{i\in A}\Delta^i\vert A\subset \{1,\ldots,n\}\}\]
and 
\[\Lambda'=\{v_A'=S v_\circ-\sum_{i\in A}\Delta'^i\vert A\subset \{1,\ldots,n\}\}.\]
%and let $\Lambda^{\text{odd}}\subset \Lambda$ and $\Lambda'^{\text{odd}}\subset \Lambda'$ be the subsets of odd elements.
Furthermore,
define $\lambda'_i=\lambda_i$ in case of a row operation and define $\lambda'_i\coloneqq\prod_{j=1}^n\lambda_j^{T_{ji}}$ in case of a column operation. 
Finally for $v_A\in \Lambda$ let
 \[\gamma_{v_A}\coloneqq\gamma_\circ \cdot \prod_{i\in A}\lambda_i\in k^\times/(k^\times)^2\]
$v'_A\in \Lambda'$ let
 \[\gamma_{v_A'}\coloneqq \gamma_\circ \cdot \prod_{i\in A}\lambda_i'\in k^\times/(k^\times)^2\]
for some fixed $\gamma_\circ\in k^\times$.

\begin{prop}
\label{prop:bijection}
%Let $M=(\Delta^1,\ldots,\Delta^n)$ and let $\Lambda $
%Let $\Lambda=\Lambda_E$ be the set of corner vertices of the original paralellepiped. 
    There is a bijection 
    \[\phi\colon \Lambda \rightarrow \Lambda'\]
    such that $\gamma_{v_A}=\gamma_{\phi(v_A)}$ in $k^\times/(k^\times)^2$ for all $A\subset \{1,\ldots,n\}$. Furthermore, 
    %in case $v_\circ$ is odd $\phi$ restricts to a bijection \[\phi^{\text{odd}}\colon \Lambda^{\text{odd}}\rightarrow \Lambda'^{\text{odd}}.\]
     %In case $v_\circ$ is not odd, we still get that $\#\Lambda^{\text{odd}}=\#\Lambda'^{\text{odd}}$. 
     \begin{equation}
     \label{eq:evenphi}
     v_A-v_\circ\equiv0\mod2 \iff \phi(v_{A})-Sv_\circ\equiv0\mod2\end{equation}
\end{prop}
\begin{proof}
\begin{figure}[ht]
    \centering
    \begin{tikzpicture}[scale=0.8]
    % Nodes
    \coordinate (A) at (0,0) {};
    \coordinate (B) at (3,0) {};
    \coordinate (C) at (2,2) {};
    \coordinate (D) at (5,2) {};
    
    \coordinate (E) at (7,0) {};
    \coordinate (F) at (10,0) {};
    \coordinate (G) at (12,2) {};
    \coordinate (H) at (15,2) {};
    
    % Parallelograms
    \draw[thick] (A)--(B)-- (D) --(C) -- (A);
    \draw[thick] (E) -- (F) -- (H) -- (G) -- (E);
    
    % Arrows
    %\draw[->, thick, blue] (A) to[out=150,in=210] (C);
    %\draw[->, thick, blue] (A) to[out=210,in=150] (B);
    %\draw[->, thick, blue] (C) to[out=30,in=330] (D);
    %draw[->, thick, blue] (B) to[out=30,in=330] (D);
    
    %\draw[->, thick, blue] (E) to[out=150,in=210] (G);
    %\draw[->, thick, blue] (E) to[out=210,in=150] (F);
    %\draw[->, thick, blue] (G) to[out=30,in=330] (H);
    %\draw[->, thick, blue] (F) to[out=30,in=330] (H);
    
    \draw[->, thick, gray, shorten >=4pt, shorten <=4pt] (A) to[out=-20,in=-160] (E);
    \draw[->, thick, gray, shorten >=4pt, shorten <=4pt] (B) to[out=-20,in=-160] (F);
    \draw[->, thick, gray, shorten >=4pt, shorten <=4pt] (C) to[out=20,in=160] (H);
    \draw[->, thick, gray, shorten >=4pt, shorten <=4pt] (D) to[out=20,in=160] (G);

    % Labels
    %\node[below left] at (A) {$v_\circ$};
    %\node[below ] at (E) {$v_\circ$};
    \node[above] at (1.5,0) {$\Delta^i$};
    \node[above] at (8.5,0) {$\Delta^i$};
    %\node[right] at (F) {$i$};
    %\node[right] at (B) {$i$};
    \node[above] at (1,1.2) {$\Delta^j$};
    \node[above] at (9.2,1.2) {$\Delta^i + \Delta^j$};
\end{tikzpicture}
    \caption{Adding the $i$th column to the $j$th column.}
    \label{fig:columntrafo}
\end{figure}
    We construct this bijection for every allowed row or column operation.
    %We use the notation
    %\[v_A=v_\circ-\sum_{i\in A}\Delta^i\in \Lambda \]
    %and
    %\[v_A'=v_\circ-\sum_{i\in A}\Delta'^i \in \Lambda'\]
    %for $A\subset \{1,\ldots,n\}$.
    We start with the column operations. %We will see that for a column operation $\phi$ will always restrict to a bijection on the odd elements.
   If the operation is swapping the $i$th column with the $j$th column in $M$ to obtain $M'$: Then one can define $\phi$ as follows. 
        \begin{itemize}
            \item When $i,j\not \in  A$ or $i,j\in A$, we define $\phi(v_A)=v_A'$, since 
            \[v_A=v_\circ-\sum_{l\in A}\Delta^l= v_\circ-\sum_{l\in A}\Delta'^l=v_A' \text{ and }\gamma_{v_A}=\gamma_\circ \cdot \prod_{l\in A}\lambda_l= \gamma_\circ \cdot \prod_{l\in A}\lambda_l' =\gamma_{v_A'}.\]
            %Similarly when $i,j \in  A$.
            
            \item When $i\not\in A$ and $j\in A$, observe that \[v'_{(A\setminus\{j\})\cup\{i\}}= v_\circ-\sum_{l\in A}\Delta'^l+\Delta'^j-\Delta'^i=v_\circ-\sum_{l\in A}\Delta^l=v_A\] and  \[\gamma_{v'_{(A\setminus\{j\})\cup \{i\}}}= \gamma_\circ\cdot \prod_{l\in (A\setminus\{j\})\cup \{i\}}\lambda_l'= \gamma_\circ\cdot \prod_{l\in A}\lambda_l  =\gamma_{v_A}.\]
            We set $\phi(v_A)=v'_{(A\setminus\{j\})\cup \{i\}}$. 
             Alike, when $i\in A$ and $j\not\in A$, let $\phi(v_A)=v'_{(A\setminus\{i\})\cup\{j\}}$.
        \end{itemize}
        %Clearly, \eqref{eq:evenphi} holds.
        %{\color{blue} replace $\gamma$ by $\sign(\det M)\gamma$}
    If the operation is multipling the $i$th column by $-1$: We define $\phi$ by $\phi(v_A)=v_A'$. 
    We have that $\gamma_{v_A}=\gamma_{\phi(v_A)}$ in $k^\times/(k^\times)^2$, since  
    $\gamma_{v_A}=\gamma_{\phi(v_A)}$ if $i\not \in A$, or $\gamma_{v_A}=\gamma_{v_A'}\cdot \lambda_i^2$ if $i\in A$.\medbreak
    
    Lastly, if the operation is adding the $i$th column to the $j$th column: We define $\phi$ as follows (see Figure \ref{fig:columntrafo} for an illustration). 
        \begin{itemize}
            \item When $j\not \in  A$, we define $\phi(v_A)=v_A'$ since $v_A=v'_A$ and $\gamma_{v_A}=\gamma_{v_A'}$. In the example in Figure \ref{fig:columntrafo}, these are the two bottom vertices.
            \item When $j\in A$ and $i\not \in A$, we define $\phi(v_A)=v'_{A\cup \{i\}}$ since $v_A=v'_{A\cup \{i\}}+2\Delta^i$ and
            \[\gamma_{v'_{A\cup \{i\}}}=\gamma_\circ\cdot \prod_{l\in A}\lambda_l \cdot \lambda_i^2= \gamma_{v_A}\in k^\times/(k^\times)^2.\]
            In the example in Figure \ref{fig:columntrafo}, the vertex $v_A$ is the upper left vertex in the left parallelogram and $v'_{A\cup \{i\}}$ is the upper right one in the right parallelogram.
            
            \item When $j\in A$ and $i\in A$, we define $\phi(v_A)=v'_{A\setminus \{i\}}$ since $v_A=v'_{A\setminus\{i\}}$  and
            \[\gamma_{v'_{A\setminus  \{i\}}}=\gamma_\circ\cdot \prod_{l\in A}\lambda_l= \gamma_{v_A}.\]
             In the example in Figure \ref{fig:columntrafo}, the vertex $v_A$ is the upper right vertex in the left parallelogram and $v'_{A\setminus \{i\}}$ is the upper left one in the right parallelogram.
        \end{itemize}
        One checks easily that $\phi$ defined as in all column operations above satisfies the relation \eqref{eq:evenphi}.
    It remains to look at the row operations. Note that for row operations we always have $\gamma_{v_A}=\gamma_{v'_A}$ for $A\subset \{1,\ldots,n\}$. We set $\phi(v_A)=v_A'$.
    Since $v_A-v_\circ=\sum_{i\in A}\Delta^i=M\cdot(\sum_{i\in A}e_i)$ and 
    $v'_A-Sv_\circ=\sum_{i\in A}\Delta'^i=SM\cdot(\sum_{i\in A}e_i)$, the relation
    \[v_A-v_\circ\equiv0\mod2 \iff v'_{A}-Sv_\circ\equiv0\mod2\]
    holds since the class of $S$ modulo two is invertible.
    % Clearly, for the first two row operations, namely multiplying a row with $-1$ and swapping two rows, this preserves odd points. When $v_\circ$ is odd, adding a row to another also preserves odd points: In this case $v_A$ is odd if and only if $\Delta^i$ has only even entries for all $i\in A$ which is equivalent to $\Delta'^i$ having only even entries for all $i\in A$ that is if and only if $v_A'$ is even.Furthermore, $v_A$ is odd if and only if $v'_A$ is odd {\color{red} maybe say something here for the third row operation, basically when a column is even, that it remains even after adding a row to another..., \color{blue} this might not be true when $v_\circ$ is not odd, but we can fix this by either only allowing special column operations, namely the ones we actually use or making a case distinction between the cases of $v_\circ$ odd or not}
    % For the first two row operations it is easy to check that the following $\phi$ does the job.
    % \begin{itemize}
    %    \item Swapping two rows: Let $\phi(v_A)=v_A'$.
    %    \item Multiplying a row with $-1$: Let $\phi(v_A)=v_A'$.
    % \end{itemize}
\end{proof}

We will prove Theorem \ref{thm:mainthm} by applying Proposition \ref{prop:bijection} several times, namely until we get a finite \'etale $k$-algebra of the form ${k[y_1,\ldots,y_n]}/{(y_i^{m_i}-\mu_i)_{i=1,\ldots,n}}$ where the $m_i$ are the elementary divisors of $M$, as in Equation~\eqref{eq:E'}, and using Proposition~\ref{prop:traceEp} and the relation in Equation~\eqref{eq:evenphi} to show that the enriched intersection multiplicity  $\mult_p(\tV_1,\ldots,\tV_n)$ is given by the sum over the odd vertices $\Lambda_E^\text{odd}\subset\Lambda_E$ of the quadratic forms given by the signed coefficients. %that in each step the sum $\sum_{v\in \Lambda_E^{\text{odd}}}\qinv{\gamma_v}\in \W(k)$ where $\Lambda_E^{\text{odd}}\subset \Lambda_E$ is the subset of odd points, stays the same. This will allow us to assume that $E=\frac{k[x_1,\ldots,x_n]}{(x_i^{m_i}-\lambda_i)_{i=1,\ldots,n}}$ and use Proposition \ref{prop:traces} to compute that \[\mult_p(\tV_1,\ldots,\tV_n)=\Tr_{E/k}\qinv{m_1\cdots m_n\cdot \gamma_\circ \cdot x^{v_\circ-e}}= \sum_{v\in \Lambda_E^{\text{odd}}}\qinv{\gamma_v}\in \W(k).\]

\begin{ex}
    Here is an example of how to get from $M=\begin{pmatrix}
        3 & 2\\ 5 & 4
    \end{pmatrix}$ to $M'=\begin{pmatrix}1& 0\\ 0&2
    \end{pmatrix}$ using only the allowed row and column operations. To shorten the example, we sometimes performed several operations in one step.
The first column of the following table consists of the matrices $M$, the second column consists of the the equations defining the $k$-algebra and the third column consists of the $\gamma_{v_A}$ in $k^\times/(k^\times)^2$ for $A\subset \{1,2\}$.
One gets from the first matrix to the second by multiplying with $T_1=\begin{pmatrix}
    1& 0\\-1&1
\end{pmatrix}$ from the right (column operation), from the second to the third by multiplying with $T_2=\begin{pmatrix}
    1&-2\\0&1
\end{pmatrix}$ from the left (column operation) and from the third to the fourth by multiplying with $S=\begin{pmatrix}
    1&0\\-1&1
\end{pmatrix}$ from the right (row operation).
\begin{center}
\begin{tabular}{c|c|c}
        $\begin{pmatrix}
            3& 2\\ 5&4
        \end{pmatrix}$ &  $\begin{array}{c}
              x_1^3x_2^5=\lambda_1 \\
              x_1^2x_2^4=\lambda_2
        \end{array}$ &$\begin{array}{c}
\gamma_\emptyset=\gamma_\circ \\
\gamma_{\{1\}}=\gamma_\circ\lambda_1\\
\gamma_{\{2\}}=\gamma_\circ\lambda_2\\
\gamma_{\{1,2\}}=\gamma_\circ\lambda_1\lambda_2
\end{array}$\\
\hline
$\begin{pmatrix}
            1& 2\\ 1&4
        \end{pmatrix}$ &  $\begin{array}{c}
              x_1x_2=\lambda_1\lambda_2^{-1} \\
              x_1^2x_2^4=\lambda_2
        \end{array}$ &$\begin{array}{c}
\gamma_\emptyset=\gamma_\circ \\
\gamma_{\{1\}}=\gamma_\circ\lambda_1\lambda_2\\
\gamma_{\{2\}}=\gamma_\circ\lambda_2\\
\gamma_{\{1,2\}}=\gamma_\circ\lambda_1
\end{array}$\\
\hline
$\begin{pmatrix}
            1& 0\\ 1&2
        \end{pmatrix}$ &  $\begin{array}{c}
              x_1x_2=\lambda_1\lambda_2^{-1} \\
              x_2^2=\lambda_1^{-2}\lambda_2^{3}
        \end{array}$ &$\begin{array}{c}
\gamma_\emptyset=\gamma_\circ \\
\gamma_{\{1\}}=\gamma_\circ\lambda_1\lambda_2\\
\gamma_{\{2\}}=\gamma_\circ\lambda_2\\
\gamma_{\{1,2\}}=\gamma_\circ\lambda_1
\end{array}$\\
\hline
$\begin{pmatrix}
            1& 0\\ 0&2
        \end{pmatrix}$ &  $\begin{array}{c}
              x_1=\lambda_1\lambda_2^{-1} \\
              x_2^2=\lambda_1^{-2}\lambda_2^{3}
        \end{array}$ &$\begin{array}{c}
\gamma_\emptyset=\gamma_\circ \\
\gamma_{\{1\}}=\gamma_\circ\lambda_1\lambda_2\\
\gamma_{\{2\}}=\gamma_\circ\lambda_2\\
\gamma_{\{1,2\}}=\gamma_\circ\lambda_1
\end{array}$\\
\end{tabular}
\end{center}
In each step, the set of odd vertices $\Lambda^\text{odd}$ is the same, either empty, or one of the sets~$\{v_\emptyset,v_{\{2\}}\}$ or $\{v_{\{1\}},v_{\{1,2\}}\}$.
Remark that in each of the three cases $\sum_{\text{$v_A$ odd}}\qinv{\gamma_{v_A}}$ is the same for every row in the table above.

\end{ex}

\begin{proof}[Proof of Theorem \ref{thm:mainthm}]

Recall from \eqref{eq:E'} that 
\[E\cong E'=\frac{k[y_1,\ldots,y_n]}{(y_i^{m_i}-\mu_i)_{i=1,\ldots,n}}\] and the $m_i$ are the elementary divisors of $M$.

As before $\Lambda_E=\{v_\circ -\sum_{i\in A}\Delta^i\vert A\subset \{1,\ldots,n\}\}$ and $\Lambda_E^{\text{odd}}$ is its subset of odd elements. Further, let $\Lambda_{E'}=\{Sv_\circ -\sum_{i\in A}m_i\cdot e_i\vert A\subset \{1,\ldots,n\}\}$.
%As before set $v_\circ=\sum_{i=1}^n I^i$, $\gamma_\circ=\sign(\det M)\cdot \prod_{i=1}^n$ and 
For $v_A=v_\circ-\sum_{i\in A} \Delta^i\in \Lambda_E$ set
\begin{equation}
\label{eq:gammav}
\gamma_{v_A}=\epsilon_P(v_A)\cdot \alpha_{v_A}=\gamma_\circ \cdot\prod_{i\in A}\lambda_i\end{equation}
and for $v'_A=Sv_\circ -\sum_{i\in A}m_i\cdot e_i$ let
\[v_A'=\gamma_\circ \cdot  \prod_{i\in A}\mu_i. \]
After applying Proposition \ref{prop:bijection} finitely many times, we get that a bijection
$\phi\colon \Lambda_E\rightarrow \Lambda_{E'}$
such that
 $\gamma_{v_A}=\gamma_{\phi(v_A)}$ in $k^\times/(k^\times)^2$ for all $A\subset \{1,\ldots,n\}$, and such that $v_A-v_\circ\equiv0\mod2$ if and only if $\phi(v_{A})-Sv_\circ\equiv0\mod2$.

If $\Lambda_E^{\text{odd}}$ is not empty, by Lemma~\ref{lm: symmetry dim n} we can assume that $v_\circ$ is odd without loss of generality.
In this case,
    \[
        \Tr_{E/k}\qqinv{m \cdot \gamma_\circ \cdot x^{v_\circ-e}}
        =\Tr_{E'/k}\qqinv{m \cdot \gamma_\circ \cdot y^{S(v_\circ-e)}}
        = \Tr_{E'/k }\qqinv{m_1\cdots m_n \cdot \gamma_\circ}
        %\overset{\ref{prop:traces}}{=}& \sum_{v'\in \Lambda'^{\text{odd}}}\qinv{\gamma_{v'}}\overset{\ref{lm:bijection}}{=}\sum_{v'\in \Lambda_E^{\text{odd}}}\qinv{\gamma_v}.
    \]
and by Proposition~\ref{prop:traceEp}
\[
    \Tr_{E'/k }\qqinv{m_1\cdots m_n \cdot \gamma_\circ}=
    \sum_{A'} \qinv{\gamma_\circ\cdot\prod_{i\in A'}\mu_i},
\]
where the sum runs over all sets $A'\subset\{1,2,\dots,n\}$ such that $m_i$ is even for every $i\in A'$, i.e. the sets such that $\phi(v_{A})-Sv_\circ\equiv0\mod2$. 
By Proposition \ref{prop:bijection}, this sum runs over all sets $A$ such that $v_A-v_\circ\equiv0\mod2$, or equivalently, over all odd vertices $v\in\Lambda_E^\text{odd}$.
Since $\gamma_{v_A}=\gamma_{\phi(v_A)}$ in $k^\times/(k^\times)^2$ for all $A\subset \{1,\ldots,n\}$, then
\[\Tr_{E/k}\qqinv{m \cdot \gamma_\circ \cdot x^{v_\circ-e}}=\sum_{v\in\Lambda_E^\text{odd}}\qinv{\gamma_v}\overset{\eqref{eq:gammav}}{=}\sum_{v\in\Lambda_E^\text{odd}}\qinv{\epsilon_P(v)\alpha_v},\]
and thus the statement follows.

Lastly, if $\Lambda_E^{\text{odd}}=\emptyset$, then there is a $j$ such that $(v_0)_j\equiv0\mod 2$ and $\Delta^i_j\equiv0\mod2$ for every $i$.
For this $j$ let 
\[L=\frac{k[x_1,\ldots, x_{j-1},y,x_{j+1},\ldots,x_n]}{(x^{\Delta'^i}-\lambda_i)_{i=1,\dots,n}}\]
where $\Delta'^i_j=\Delta^i_j/2$ and $\Delta'^i_l=\Delta^i_l$ for $l\neq j$.
Then $E=L[x_j]/(x_j^2-y)$ and thus $\dim_L E=2$. Hence,
\begin{align*}
        \Tr_{E/k}\qqinv{m \cdot \gamma_\circ \cdot x^{v_\circ-e}}
        =\Tr_{L/k}\circ\Tr_{E/L}\qqinv{m \cdot \gamma_\circ \cdot x^{v_\circ-e}}
        \overset{\ref{prop:traces}}{=}\Tr_{L/k}(0)=0
    \end{align*}
in $\W(k)$, by Proposition~\ref{prop:traces} since $x_j\in E\setminus L$ is not a square in $E$. This agrees with the formula in the statement of the theorem since $P$ has no odd corner vertices.
    
\end{proof}

%% file: 6applications.tex
\section{Enriched Tropical B\'ezout and Bernstein-Kushnirenko \allowbreak theorems}
\label{section:applications}

We use properties of toric varieties to give applications of the computation we obtained in the precedent sections. For more details on toric varieties we refer to \cite{Fulton} and \cite{gkz}. 
 
\subsection{A tropical proof of B\'ezout's theorem enriched in $\GW(k)$}

With the combinatorial formulas in Theorem \ref{thm:mainthm} for the enriched intersection multiplicity, we can quadratically enrich the proof of the tropical B\'{e}zout theorem \ref{thm:tropical Bezout}. The resulting count agrees with McKean's nontropical B\'{e}zout's theorem \ref{thm:McKean} in the relatively orientable case. 
In the non-relatively orientable case, we do not get an invariant result for the sum of enriched intersection multiplicities at the intersection points as expected. However, our methods tell us all possible counts for this sum.

The proof of the enriched tropical B\'{e}zout theorem is an easy corollary of the following Proposition.

\begin{prop}
\label{prop:number of parallelepipeds}
Let $V_1,\ldots,V_n$ be $n$ tropical hypersurfaces in $\mathbb{R}^n$ with associated Newton polytopes \allowbreak $\Delta_1,\ldots,\Delta_n$, respectively. Let $v$ be a lattice point in the interior of the Minkowski sum $\Delta_1+\ldots+\Delta_n$ and let
 \[P_v\coloneqq\{\text{$P$ in $\operatorname{DS}\left(\bigcup_{i=1}^n V_i\right)$ dual to some $p\in V_1\cap\ldots\cap V_n$, s.t. $v$ is a corner vertex of $P$}\}.\]
If the hypersurfaces $V_1,\ldots,V_n$ intersect tropically transversely, then
\begin{enumerate}
\item The cardinality $\operatorname{Card}(P_v)$ of $P_v$ is even.
\item There are equally many parallelepipeds $P$ in $P_v$ such that the sign $\epsilon_P(v)$ (as defined in \eqref{eq:signdimn}) in $P$ is positive as there are with negative sign
\[\operatorname{Card}(\{P\in P_v\mid \epsilon_P(v)=+1\})=\operatorname{Card}(\{P\in P_v\mid \epsilon_P(v)=-1\}).\]
\end{enumerate}

\end{prop}

\begin{figure}[t]
\begin{subfigure}[b]{0.45\textwidth}
\begin{tikzpicture}[line cap=round,line join=round,scale=0.7]
\clip(-1,-1) rectangle (7,8);
\fill[line width=2pt,color=purple,fill=purple,fill opacity=0.1] (1,5) -- (2,6) -- (4,6) -- (6,4) -- (6,2) -- (5,1) -- (3,1) -- (1,3) -- cycle;
\draw [line width=2pt,color=red] (-1,2)-- (1,2);
\draw [line width=2pt,color=red] (1,2)-- (2,1);
\draw [line width=2pt,color=red] (2,1)-- (2,-1);
\draw [line width=2pt,color=red] (1,2)-- (1,5);
\draw [line width=2pt,color=red] (1,5)-- (-1,5);
\draw [line width=2pt,color=red] (1,5)-- (3,7);
\draw [line width=2pt,color=red] (3,7)-- (3,8);
\draw [line width=2pt,color=blue] (4,0)-- (4,-1);
\draw [line width=2pt,color=red] (3,7)-- (7,7);
\draw [line width=2pt,color=red] (2,1)-- (7,1);
\draw [line width=2pt,color=blue] (0,4)-- (4,0);
\draw [line width=2pt,color=blue] (4,0)-- (6,2);
\draw [line width=2pt,color=blue] (6,2)-- (7,2);
\draw [line width=2pt,color=blue] (6,2)-- (6,4);
\draw [line width=2pt,color=blue] (6,4)-- (4,6);
\draw [line width=2pt,color=blue] (4,6)-- (0,6);
\draw [line width=2pt,color=blue] (0,6)-- (0,4);
\draw [line width=2pt,color=blue] (0,4)-- (-1,4);
\draw [line width=2pt,color=blue] (0,6)-- (-1,7);
\draw [line width=2pt,color=blue] (4,6)-- (4,8);
\draw [line width=2pt,color=blue] (6,4)-- (7,4);
\draw [color=purple,fill=purple] (1,3) circle (2.5pt);
\draw node at (0.7,2.8)  {\large$p_1$};
\draw [color=purple,fill=purple] (3,1) circle (2.5pt);
\draw node at (3,0.6)  {\large$q_1$};
\draw [color=purple,fill=purple] (2,6) circle (2.5pt);
\draw node at (2,6.4)  {\large$p_2$};
\draw [color=purple,fill=purple] (5,1) circle (2.5pt);
\draw node at (5,0.6)  {\large$q_2$};
\draw node at (4.2,-0.6)  {\large$2$};
\draw [color=blue] node at (6.6,4.4)  {\large$C$};
\draw [color=red] node at (-0.6,1.6)  {\large$V_n$};
\draw [color=purple] node at (3.5,3.5)  {\Large$R_v$};
\end{tikzpicture}
\end{subfigure}\quad
\begin{subfigure}[b]{0.45\textwidth}
\begin{tikzpicture}[line cap=round,line join=round,scale=0.7]
\clip(-1,-1) rectangle (7,8);
\fill[line width=2pt,color=purple,fill=purple,fill opacity=0.01] (1,5) -- (2,6) -- (4,6) -- (6,4) -- (6,2) -- (5,1) -- (3,1) -- (1,3) -- cycle;
\draw [line width=2pt,opacity=0.1,color=red] (-1,2)-- (1,2);
\draw [line width=2pt,opacity=0.1,color=red] (1,2)-- (2,1);
\draw [line width=2pt,opacity=0.1,color=red] (2,1)-- (2,-1);
\draw [line width=2pt,opacity=0.1,color=red] (1,2)-- (1,5);
\draw [line width=2pt,opacity=0.1,color=red] (1,5)-- (-1,5);
\draw [line width=2pt,opacity=0.1,color=red] (1,5)-- (3,7);
\draw [line width=2pt,opacity=0.1,color=red] (3,7)-- (3,8);
\draw [line width=2pt,opacity=0.1,color=blue] (4,0)-- (4,-1);
\draw [line width=2pt,opacity=0.1,color=red] (3,7)-- (7,7);
\draw [line width=2pt,opacity=0.1,color=red] (2,1)-- (7,1);
\draw [line width=2pt,opacity=0.1,color=blue] (0,4)-- (4,0);
\draw [line width=2pt,opacity=0.1,color=blue] (4,0)-- (6,2);
\draw [line width=2pt,opacity=0.1,color=blue] (6,2)-- (7,2);
\draw [line width=2pt,opacity=0.1,color=blue] (6,2)-- (6,4);
\draw [line width=2pt,opacity=0.1,color=blue] (6,4)-- (4,6);
\draw [line width=2pt,opacity=0.1,color=blue] (4,6)-- (0,6);
\draw [line width=2pt,opacity=0.1,color=blue] (0,6)-- (0,4);
\draw [line width=2pt,opacity=0.1,color=blue] (0,4)-- (-1,4);
\draw [line width=2pt,opacity=0.1,color=blue] (0,6)-- (-1,7);
\draw [line width=2pt,opacity=0.1,color=blue] (4,6)-- (4,8);
\draw [line width=2pt,opacity=0.1,color=blue] (6,4)-- (7,4);
\draw [line width=2pt,color=purple] (1,3)-- (3,1);
\draw[color=purple] node at (2.25,2.25)  {\large$\gamma_1$};
\draw [line width=2pt,color=purple] (2,6)-- (4,6) -- (6,4) -- (6,2) -- (5,1);
\draw[color=purple] node at (5.25,5.25)  {\large$\gamma_2$};
\draw [color=purple,fill=purple] (1,3) circle (2.5pt);
\draw node at (0.7,2.8)  {\large$p_1$};
\draw [color=purple,fill=purple] (3,1) circle (2.5pt);
\draw node at (3,0.6)  {\large$q_1$};
\draw [color=purple,fill=purple] (2,6) circle (2.5pt);
\draw node at (2,6.4)  {\large$p_2$};
\draw [color=purple,fill=purple] (5,1) circle (2.5pt);
\draw node at (5,0.6)  {\large$q_2$};
\draw node at (4.2,-0.6)  {\large$2$};
\draw [color=blue,opacity=0.3] node at (6.6,4.4)  {\large$C$};
\draw [color=red,opacity=0.3] node at (-0.6,1.6)  {\large$V_n$};
\draw [color=purple,opacity=0.3] node at (3.5,3.5)  {\Large$R_v$};
\end{tikzpicture}
\end{subfigure}
    \caption{Partition of intersection points into pairs next to the region $R_v$ associated to a vertex.}
    \label{fig:regions}
\end{figure}

\begin{proof}

%{\color{red} Let's illustrate this with an example for curves.}
Due to the transversality hypothesis, the hypersurfaces $V_1,\ldots,V_{n-1}$ intersect along a tropical curve $C\subset\R^n$. This curve intersects $V_n$ tropically transversely. 
Let us denote by $R_{v}\subset\R^n$ the connected component of $\R^n\setminus V_1\cup \ldots\cup V_n$ where the monomial of the exponent $v$ is maximal (see Figure \ref{fig:regions}). %in the product of defining equations of the hypersurfaces $V_1,\ldots,V_n$. 
Since $v$ is an inner lattice point of the dual polytope, the region $R_v$ is a bounded polytope.
If $P_v$ is empty, our assertion follows. Otherwise, let $p\in C\cap V_n$ be an intersection point such that its dual polytope $P\in P_v$. 
Let $\gamma$ be the connected component of~$C\cap \partial R_v$ containing the point $p$. We claim that $\gamma$ is a piecewise linear path. Namely, the set~$\gamma$ is formed by two segments of edges of $C$ which contain an intersection point with $V_n$, together with possible bounded edges (for example $\gamma_2$ in Figure \ref{fig:regions}), or only one segment containing two intersection points with $V_n$ (for example $\gamma_1$ in Figure \ref{fig:regions}). In particular, exactly one of the two intersection points of $\gamma$ with $V_n$ is $p$. If a vertex $w$ of $C$ is in $\gamma$, its valency in $\gamma$ is $2$, corresponding to the edges in $C$ adjacent to the region $R_v$. Since $R_v$ is a bounded polytope, the curve $\gamma$ is compact, having an endpoint~$q$ that is in $V_n$. Indeed, the point $q\neq p$ and cannot be a vertex of $C$, %or its valency in $\gamma$ would be $1$, 
hence it is an inner point of an edge of $C$. Since there are no changes in the monomials where the maximum is achieved in the interior of an edge, this change is produced by the hypersurface $V_n$. Therefore, the paths $\gamma$ establish a partition of the parallelepipeds in $P_v$ into pairs $P_v=\bigsqcup \{P_\gamma,Q_\gamma\}$ with $P_\gamma$ and $Q_\gamma$ the dual parallelepipeds of the endpoints $p_\gamma$ and $q_\gamma$ of $\gamma$ for each such path $\gamma$ (see Figure \ref{fig:regions}).
%a pairing between the intersection points of $C\cap V_n$ adjacent to $R_v$.

We transfer the frame in $p$ through $\gamma$ to show that the polytopes corresponding to the endpoints of $\gamma$ have opposite sign at the vertex $v$ (see Figure \ref{fig:paths}).
For that, let us start by recalling that the sign $\epsilon_P(v)$ is the sign of the determinant $\left(\Delta^i\right)_{i=1}^n$, where every $\Delta^i$ has been oriented in such a way that~${v+\Delta^i\in P}$ (see \eqref{eq:epsilonsgeometrically}). This oriented vector is the normal vector of the facet of $V_i$ pointing outwards to the region $R_v$.
Let us define $w_p\coloneqq\bigwedge_{i=1}^{n-1}\Delta^i$ (the vector of alternating minors of the $(n-1)\times n$ matrix $(\pm\Delta^i)_{i=1}^{n-1}$).
The vector $w_p$ is a direction vector of the edge of $C$ containing $p$, albeit not a primitive one. Moreover, the sign $\epsilon_P(v)=\sign(\det\left(\Delta^i\right)_{i=1}^n)=\sign(w_p\cdot \Delta^n)$ can be computed as the sign of the inner product of $w_p$ with the normal vector of the facet of $V_n$ containing $p$, oriented outwards the region~$R_v$. 
We can define $w_{\gamma(t)}$ for every point of $\gamma$ that is an inner point of an edge of $C$. 
If we orient $\gamma$ as a path starting at $p$, at every point $\gamma(t)$, the orientations of $w_{\gamma(t)}$ and $\gamma$ either coinside for all $t$ or are opposite for all $t$, since the relative position of the normal vectors of the facets of $V^i$ at $\gamma(t)$ does not change in the boundary of $R_v$. This implies that exactly one of the vectors $w_p$ at $p\in V_n$ or $w_q$ at $q\in V_n$ is oriented towards the region $R_v$ while the other one is not. Hence, the endpoints of $\gamma$ have opposite signs.
\end{proof}

\begin{figure}[t]
\begin{subfigure}[b]{0.30\textwidth}
\begin{tikzpicture}[line cap=round,line join=round,scale=0.5]
\clip(-1,-1) rectangle (7,8);
\fill[line width=2pt,color=purple,fill=purple,fill opacity=0.01] (1,5) -- (2,6) -- (4,6) -- (6,4) -- (6,2) -- (5,1) -- (3,1) -- (1,3) -- cycle;
\draw [line width=2pt,opacity=0.1,color=red] (-1,2)-- (1,2);
\draw [line width=2pt,opacity=0.1,color=red] (1,2)-- (2,1);
\draw [line width=2pt,opacity=0.1,color=red] (2,1)-- (2,-1);
\draw [line width=2pt,opacity=0.1,color=red] (1,2)-- (1,5);
\draw [line width=2pt,opacity=0.1,color=red] (1,5)-- (-1,5);
\draw [line width=2pt,opacity=0.1,color=red] (1,5)-- (3,7);
\draw [line width=2pt,opacity=0.1,color=red] (3,7)-- (3,8);
\draw [line width=2pt,opacity=0.1,color=red] (3,7)-- (7,7);
\draw [line width=2pt,opacity=0.1,color=red] (2,1)-- (7,1);
% \draw [line width=2pt,opacity=0.1,color=blue] (4,0)-- (4,-1);
% \draw [line width=2pt,opacity=0.1,color=blue] (0,4)-- (4,0);
% \draw [line width=2pt,opacity=0.1,color=blue] (4,0)-- (6,2);
% \draw [line width=2pt,opacity=0.1,color=blue] (6,2)-- (7,2);
% \draw [line width=2pt,opacity=0.1,color=blue] (6,2)-- (6,4);
% \draw [line width=2pt,opacity=0.1,color=blue] (6,4)-- (4,6);
% \draw [line width=2pt,opacity=0.1,color=blue] (4,6)-- (0,6);
% \draw [line width=2pt,opacity=0.1,color=blue] (0,6)-- (0,4);
% \draw [line width=2pt,opacity=0.1,color=blue] (0,4)-- (-1,4);
% \draw [line width=2pt,opacity=0.1,color=blue] (0,6)-- (-1,7);
% \draw [line width=2pt,opacity=0.1,color=blue] (4,6)-- (4,8);
% \draw [line width=2pt,opacity=0.1,color=blue] (6,4)-- (7,4);
%\draw [color=blue,opacity=0.3] node at (6.6,4.4)  {\large$C$};
%\draw node at (4.2,-0.6)  {\large$2$};
\draw [line width=2pt,color=purple,opacity=0.3] (1,3)-- (3,1);
\draw [color=purple,opacity=0.3] node at (2.25,2.25)  {\large$\gamma_1$};
\draw [line width=2pt,color=purple,opacity=0.3] (2,6)-- (4,6) -- (6,4) -- (6,2) -- (5,1);
\draw [color=purple,opacity=0.3] node at (5.25,5.25)  {\large$\gamma_2$};
\draw [color=red,opacity=0.3] node at (-0.6,1.6)  {\large$V_n$};
\draw [color=purple,opacity=0.3] node at (3.5,3.5)  {\Large$R_v$};
\draw[->] (1,3) -- (-0.25,3);
\draw node at (0.25,3.6)  {$\Delta^n(p_1)$};
\draw[->] (2,6) -- (0.75,7.25);
\draw node at (2.25,7.5)  {$\Delta^n(p_2)$};
\draw[->] (3,1) -- (3,-0.25);
\draw node at (1.5,-0.25)  {$\Delta^n(q_1)$};
\draw[->] (5,1) -- (5,-0.25);
\draw [color=purple,fill=purple] (1,3) circle (2.5pt);
\draw [color=purple,fill=purple] (3,1) circle (2.5pt);
\draw [color=purple,fill=purple] (2,6) circle (2.5pt);
\draw [color=purple,fill=purple] (5,1) circle (2.5pt);
\end{tikzpicture}
\begin{tikzpicture}[overlay,scale=0.5]
\draw node at (-0.75,0.75)  {$\Delta^n(q_2)$};    
\end{tikzpicture}
\end{subfigure}\quad
\begin{subfigure}[b]{0.30\textwidth}
\begin{tikzpicture}[line cap=round,line join=round,scale=0.5]
\clip(-1,-1) rectangle (7,8);
\fill[line width=2pt,color=purple,fill=purple,fill opacity=0.01] (1,5) -- (2,6) -- (4,6) -- (6,4) -- (6,2) -- (5,1) -- (3,1) -- (1,3) -- cycle;
% \draw [line width=2pt,opacity=0.1,color=red] (-1,2)-- (1,2);
% \draw [line width=2pt,opacity=0.1,color=red] (1,2)-- (2,1);
% \draw [line width=2pt,opacity=0.1,color=red] (2,1)-- (2,-1);
% \draw [line width=2pt,opacity=0.1,color=red] (1,2)-- (1,5);
% \draw [line width=2pt,opacity=0.1,color=red] (1,5)-- (-1,5);
% \draw [line width=2pt,opacity=0.1,color=red] (1,5)-- (3,7);
% \draw [line width=2pt,opacity=0.1,color=red] (3,7)-- (3,8);
% \draw [line width=2pt,opacity=0.1,color=red] (3,7)-- (7,7);
% \draw [line width=2pt,opacity=0.1,color=red] (2,1)-- (7,1);
\draw [line width=2pt,opacity=0.1,color=blue] (4,0)-- (4,-1);
\draw [line width=2pt,opacity=0.1,color=blue] (0,4)-- (4,0);
\draw [line width=2pt,opacity=0.1,color=blue] (4,0)-- (6,2);
\draw [line width=2pt,opacity=0.1,color=blue] (6,2)-- (7,2);
\draw [line width=2pt,opacity=0.1,color=blue] (6,2)-- (6,4);
\draw [line width=2pt,opacity=0.1,color=blue] (6,4)-- (4,6);
\draw [line width=2pt,opacity=0.1,color=blue] (4,6)-- (0,6);
\draw [line width=2pt,opacity=0.1,color=blue] (0,6)-- (0,4);
\draw [line width=2pt,opacity=0.1,color=blue] (0,4)-- (-1,4);
\draw [line width=2pt,opacity=0.1,color=blue] (0,6)-- (-1,7);
\draw [line width=2pt,opacity=0.1,color=blue] (4,6)-- (4,8);
\draw [line width=2pt,opacity=0.1,color=blue] (6,4)-- (7,4);
\draw [line width=2pt,color=purple,opacity=0.3] (1,3)-- (3,1);
\draw[color=purple] node at (2.35,2.35)  {\large$\gamma_1$};
\draw [line width=2pt,color=purple,opacity=0.3] (2,6)-- (4,6) -- (6,4) -- (6,2) -- (5,1);
\draw[color=purple] node at (5.35,5.35)  {\large$\gamma_2$};
\draw [color=blue,opacity=0.3] node at (4.3,-0.6)  {$2$};
\draw [color=blue,opacity=0.3] node at (6.6,4.4)  {\large$C$};
%\draw [color=red,opacity=0.3] node at (-0.6,1.6)  {\large$V_n$};
\draw [color=purple,opacity=0.3] node at (3.5,3.5)  {\Large$R_v$};
\draw[->] (1,3) -- (-0.25,4.25);
\draw node at (0.75,4)  {$w_{p_1}$};
\draw[->] (2,6) -- (3.25,6);
\draw node at (3,6.75)  {$w_{p_2}$};
\draw[->] (3,1) -- (1.75,2.25);
\draw node at (1.5,1)  {$w_{q_1}$};
\draw[->] (5,1) -- (3.75,-0.25);
\draw node at (5,-0.25)  {$w_{q_2}$};
\draw [color=purple,fill=purple] (1,3) circle (2.5pt);
\draw [color=purple,fill=purple] (3,1) circle (2.5pt);
\draw [color=purple,fill=purple] (2,6) circle (2.5pt);
\draw [color=purple,fill=purple] (5,1) circle (2.5pt);
\end{tikzpicture}
\end{subfigure}
\begin{subfigure}[b]{0.30\textwidth}
\begin{tikzpicture}[line cap=round,line join=round,scale=0.5]
\clip(-1,-1) rectangle (7,8);
\fill[line width=2pt,color=purple,fill=purple,fill opacity=0.01] (1,5) -- (2,6) -- (4,6) -- (6,4) -- (6,2) -- (5,1) -- (3,1) -- (1,3) -- cycle;
\draw [line width=2pt,opacity=0.1,color=red] (-1,2)-- (1,2);
\draw [line width=2pt,opacity=0.1,color=red] (1,2)-- (2,1);
\draw [line width=2pt,opacity=0.1,color=red] (2,1)-- (2,-1);
\draw [line width=2pt,opacity=0.1,color=red] (1,2)-- (1,5);
\draw [line width=2pt,opacity=0.1,color=red] (1,5)-- (-1,5);
\draw [line width=2pt,opacity=0.1,color=red] (1,5)-- (3,7);
\draw [line width=2pt,opacity=0.1,color=red] (3,7)-- (3,8);
\draw [line width=2pt,opacity=0.1,color=blue] (4,0)-- (4,-1);
\draw [line width=2pt,opacity=0.1,color=red] (3,7)-- (7,7);
\draw [line width=2pt,opacity=0.1,color=red] (2,1)-- (7,1);
\draw [line width=2pt,opacity=0.1,color=blue] (0,4)-- (4,0);
\draw [line width=2pt,opacity=0.1,color=blue] (4,0)-- (6,2);
\draw [line width=2pt,opacity=0.1,color=blue] (6,2)-- (7,2);
\draw [line width=2pt,opacity=0.1,color=blue] (6,2)-- (6,4);
\draw [line width=2pt,opacity=0.1,color=blue] (6,4)-- (4,6);
\draw [line width=2pt,opacity=0.1,color=blue] (4,6)-- (0,6);
\draw [line width=2pt,opacity=0.1,color=blue] (0,6)-- (0,4);
\draw [line width=2pt,opacity=0.1,color=blue] (0,4)-- (-1,4);
\draw [line width=2pt,opacity=0.1,color=blue] (0,6)-- (-1,7);
\draw [line width=2pt,opacity=0.1,color=blue] (4,6)-- (4,8);
\draw [line width=2pt,opacity=0.1,color=blue] (6,4)-- (7,4);
\draw [line width=2pt,color=purple,opacity=0.1] (1,3)-- (3,1);
\draw [color=purple,opacity=0.3] node at (2.35,2.35)  {\large$\gamma_1$};
\draw [line width=2pt,color=purple,opacity=0.1] (2,6)-- (4,6) -- (6,4) -- (6,2) -- (5,1);
\draw [color=purple,opacity=0.3] node at (5.35,5.35)  {\large$\gamma_2$};
\draw [color=blue,opacity=0.3] node at (4.3,-0.6)  {$2$};
\draw [color=blue,opacity=0.3] node at (6.6,4.4)  {\large$C$};
\draw [color=red,opacity=0.3] node at (-0.6,1.6)  {\large$V_n$};
\draw [color=purple,opacity=0.3] node at (3.5,3.5)  {\Large$R_v$};
\draw[->] (1,3) -- (-0.25,4.25);
%\draw node at (0.75,4)  {$w_{p_1}$};
\draw[->] (2,6) -- (3.25,6);
%\draw node at (3,6.75)  {$w_{p_2}$};
\draw[->] (3,1) -- (1.75,2.25);
%\draw node at (1.5,1)  {$w_{q_1}$};
\draw[->] (5,1) -- (3.75,-0.25);
%\draw node at (5,-0.25)  {$w_{q_2}$};
\draw[->] (1,3) -- (-0.25,3);
%\draw node at (0.25,3.6)  {$\Delta^n(p_1)$};
\draw[->] (2,6) -- (0.75,7.25);
%\draw node at (2.25,7.5)  {$\Delta^n(p_2)$};
\draw[->] (3,1) -- (3,-0.25);
%\draw node at (1.5,-0.25)  {$\Delta^n(q_1)$};
\draw[->] (5,1) -- (5,-0.25);
\draw [color=purple,fill=purple] (1,3) circle (2.5pt);
%\draw node at (0.7,2.8)  {\large$p_1$};
\draw [color=purple,fill=purple] (3,1) circle (2.5pt);
%\draw node at (3,0.6)  {\large$q_1$};
\draw [color=purple,fill=purple] (2,6) circle (2.5pt);
%\draw node at (2,6.4)  {\large$p_2$};
\draw [color=purple,fill=purple] (5,1) circle (2.5pt);
%\draw node at (5,0.6)  {\large$q_2$};
\draw [color=white,fill=white] (2.25,0.75) circle (13pt);
\draw [color=red] node at (2.25,0.75) {$(-)$};
\draw [color=white,fill=white] (1.25,3.75) circle (13pt);
\draw [color=blue] node at (1.25,3.75) {$(+)$};
\draw [color=white,fill=white] (2.25,6.75) circle (13pt);
\draw [color=red] node at (2.25,6.75) {$(-)$};
\draw [color=white,fill=white] (5,1.75) circle (13pt);
\draw [color=blue] node at (5,1.75) {$(+)$};
\end{tikzpicture}
\end{subfigure}

    \caption{Transfer of frames defining the local sign of vertex relative to an intersection point.}
    \label{fig:paths}
\end{figure}

Let $\Delta_d$ be the Newton polygon of a general degree $d$ polynomial in $n$ variables.
Note that if $V_1,\dots, V_n$ are tropical hypersurfaces with Newton polytopes $\Delta_{d_1},\dots,\Delta_{d_n}$, then the union of hypersurfaces $V_1\cup\ldots\cup V_n$ has Newton polytope~$\Delta_{d_1+\dots+d_n}$.

\begin{coro}[Enriched tropical B\'ezout]
\label{coro:EnrichedTropicalBezout}
Assume that $\operatorname{char}k\neq 2$.
Let $\tV_1,\ldots,\tV_n$ be enriched tropical hypersurfaces in $\mathbb{R}^n$ with Newton polytopes $\Delta_{d_1},\ldots, \Delta_{d_n}$. Assume that $\tV_1,\ldots,\tV_n$ intersect tropically transversally.
If  $\sum_{i=1}^nd_i\equiv n+1 \mod 2$, then 
\[\sum_{p\in \tV_1\cap\ldots\cap \tV_n}\widetilde{\operatorname{mult}}_p(\tV_1,\ldots,\tV_n)=\frac{d_1\cdots d_n}{2}\;h\in\GW(k).\]
\end{coro}
\begin{ex}
We continue Example \ref{ex:enrichedMult} as an example of the statement in Corollary \ref{coro:EnrichedTropicalBezout}, that is we compute sum of the enriched tropical intersection multiplicities of the interesection of the two enriched tropical curves in Figure \ref{fig:innerdual}.
%there are three odd points in the interior of the Newton polygon. %The points $(3,1)$ and $(1,3)$ are both vertices of two parallelograms corresponding to intersection points, while $(1,1)$ is not a vertex of a parallelogram. 
Again let $\alpha_{(1,1)}$, $\alpha_{(1,3)}$ and $\alpha_{(3,1)}$ be the coefficients at the odd vertices $(1,1)$, $(1,3)$ and $(3,1)$.
Recall from Example \ref{ex:enrichedMult} that the enriched intersection multiplicities at the three intersection points are $\qinv{\alpha_{(1,3)}}$, $\qinv{-\alpha_{(1,3)}}+\qinv{-\alpha_{(3,1)}}$ and $\qinv{\alpha_{3,1}}+h$. 
Summing up the intersection multiplicities and using the identity $\qinv{a}+\qinv{-a}=h$ for $a\in k^\times$ from Remark \ref{rmk:hyperbolic} we get
\[\qinv{\alpha_{(1,3)}}+\qinv{-\alpha_{(1,3)}}+\qinv{-\alpha_{(3,1)}}+\qinv{\alpha_{3,1}}+h=3h\] which coincides with McKean's enriched (non-tropical) B\'ezout Theorem \ref{thm:McKean}.
\end{ex}

\begin{proof}[Proof of Corollary \ref{coro:EnrichedTropicalBezout}]

Note that in the relatively orientable case, that is, when $\sum_{i=1}^n d_i$
 is congruent to $n+1\mod 2$, all odd points in the dual subdivision of $\tV_1\cup \ldots\cup \tV_n$ lie in the interior of $\Delta_{d_1+\ldots+d_n}=\Delta_1+\ldots+\Delta_n$ and none on the boundary.

We know from the classical tropical non-enriched B\'ezout theorem \ref{thm:tropical Bezout} that the sum of $\operatorname{mult}_p(V_1,\ldots,V_n)$ over all the intersections of the $V_i$ is equal to $d_1\cdots d_n$. Since $$\operatorname{rank}\mult_p(\tV_1,\ldots,\tV_n)=\operatorname{mult}_p(V_1,\ldots,V_n)$$ we get taht the rank of 
\[\sum_{\text{intersections } p}\widetilde{\operatorname{mult}}_p(\tV_1,\ldots,\tV_n)\]
equals $d_1\cdots d_n$.

For $v$ an odd point in the interior of $\Delta_{d_1+\ldots+d_n}$ let $N(v)$ be the number of parallelepipeds in the dual subdivision of $\tV_1\cup\ldots\cup \tV_n$ which correspond to an intersection of the $V_i$. Let $\alpha_v$ be the coefficient of $v$ in $\operatorname{DS(\tV_1\cup\ldots \cup \tV_n)}$.
By Theorem \ref{thm:mainthm} and Proposition \ref{prop:number of parallelepipeds} we have
\begin{align*}
    \sum_{p}\mult_p(\tV_1,\ldots,\tV_n)=\sum_{v}\left(\frac{N(v)}{2}\qinv{ \alpha_v}+\frac{N(v)}{2}\qinv{ -\alpha_v}\right) \text{ in }\operatorname{W}(k),
\end{align*}
where the first sum runs over the intersection points of the $V_i$ and the second sum runs over the odd vertices in the interior of $\Delta_{d_1+\ldots+d_n}$. 
Since $\qinv{ a}+\qinv{ -a} =0$ in $\operatorname{W}(k)$ by (ii) in Remark \ref{rmk:hyperbolic}, we get that $\sum_{v}\left(\frac{n(v)}{2}\qinv{ \alpha_v}+\frac{n(v)}{2}\qinv{ -\alpha_v}\right)=0$ in $\operatorname{W}(k)$ and thus
$\sum_{p}\mult_p(\tV_1,\ldots,\tV_n)$ equals a multiple of $h$ in~$\GW(k)$. 

%We know from the (non-enriched) tropical B\'ezout theorem \ref{thm:tropical Bezout} that the sum of volumns of the parallelepipeds in the dual subdivision of $V_1\cup\ldots\cup V_n$ corresponding to intersections of the $V_i$ is equal to $d_1\cdots d_n$. In other words 
%\[\sum_{\text{intersections } p}\widetilde{\operatorname{mult}}_p(\tV_1,\ldots,\tV_n)\]
%equals $d_1\cdots d_n$.

%Now note that $\operatorname{rank}\widetilde{\operatorname{mult}}_p(\tV_1,\ldots,\tV_n)=\operatorname{mult}_p(V_1,\ldots,V_n)$.
%By Example \ref{ex:GW(C)} the classical tropical B\'ezout theorem tells us the rank of the sum $\sum_{p}\mult_p(\tV_1,\ldots,\tV_n)$.

Finally, recall that an element of $\GW(k)$ is determined by its rank and its image in the Witt group $\W(k)$ (see Remark \ref{rmk:GWandW}). Thus we get
\[\sum_{p}\mult_p(\tV_1,\ldots,\tV_n)=\frac{d_1\cdots d_n}{2}\; h\text{ in }\GW(k).\]

\end{proof}

\begin{rmk}
\label{rmk:BezoutProof}
This gives a new proof of McKean's non-tropical quadratically enriched B\'ezout's theorem enriched \eqref{eq:McKean} with the following argument. 
Recall that one proves this theorem by computing the $\A^1$-Euler number $n^{\A^1}(V_k)$ of the vector bundle
\[V_k= \mathcal{O}_{\mathbb{P}^n_k}(d_1)\oplus\ldots\oplus\mathcal{O}_{\mathbb{P}^n_k}(d_n)\longrightarrow \mathbb{P}_k^n.\]
Since $\mathbb{P}^n$ is smooth and proper over $\Z$, the problem can actually defined over $\Z$ %{\color{blue} double check with Stephen} 
and there is a well defined answer $n^{\A^1}(V_\Z)\in \GW(\Z)$.
By our main theorem \ref{thm:mainthm} we have that 
\[n^{\A^1}(V_{k\Puiseux})=\sum_{p\in \tV_1\cap\ldots\cap\tV_n}\mult_p(\tV_1,\ldots,\tV_n)=\frac{d_1\cdots d_n}{2}\;h\in \GW(k\Puiseux)\]
whenever $k$ is a field of characteristic $0$ or characteristic bigger that the maximum of the $d_i$'s.
The natural map $\GW(\Z)\rightarrow \GW(\R)\cong \GW(\R\Puiseux)$ is an isomorphism that maps $n^{\A^1}(V_\Z)$ to $n^{\A^1}(V_{\R\Puiseux})$. Hence, we get that $n^{\A^1}(V_\Z)$ is equal to $\frac{d_1\cdots d_n}{2}h\in\GW(\Z)$. There is a natural map $\GW(\Z)\rightarrow \GW(k)$ which maps $n^{\A^1}(V_\Z)\rightarrow n^{\A^1}(V_k)$ which proves the theorem over an arbitrary field $k$.

%and let $V_{k\Puiseux}$ be its base change to the field $k\Puiseux$ of Puiseux series. Recall that McKean showed the enriched B\'ezout theorem by computing the $\A^1$-Euler number of $V$.
%Furthermore, recall that the $\A^1$-Euler number of $V_{k\Puiseux}$ equals the sum of local indices at the zeros of a general section of~$V_{k\Puiseux}$.
%A general section of~$V_{k\Puiseux}$ is defined by $n$ general polynomials over $k\Puiseux$ which give rise to enriched tropical hypersurfaces $\tV_1,\ldots,\tV_n$.
%We defined the enriched intersection multiplicities of the corresponding enriched tropical hypersurface to be this local index. Thus it follows directly from Corollary \ref{coro:EnrichedTropicalBezout} that 
%\[n^{\A^1}(V_{k\Puiseux})=\sum_{p\in \tV_1\cap\ldots\cap\tV_n}\mult_p(\tV_1,\ldots,\tV_n)=\frac{d_1\cdots d_n}{2}\;h\in \GW(k\Puiseux).\]
%Euler classes commute with base change. Hence, the image of the Euler class $n^{\A^1}(V)$ under the map $\GW(k)\longrightarrow \GW(V_{k\Puiseux})$ induced by $k\longrightarrow k\Puiseux$ is
%$n^{\A^1}( V_{k\Puiseux})$. Since this map is an isomorphism and it sends a generator $\qinv{ a} $ of $\GW(k)$ to $\qinv{ a}\in \GW(k\Puiseux)$, we conclude that  
%\[n^{\A^1}(V)=\frac{d_1\cdots d_n}{2}\;h\in \GW(k).\]
\end{rmk}

\subsubsection{Non-relatively orientable case}

In the non-relatively orientable case, that is, when $\sum_{i=1}^nd_i\not \equiv n+1\mod 2$, we do not get an invariant count. This can also be seen in our proof for the enriched tropical B\'{e}zout theorem: In case $\sum_{i=1}^nd_i\not \equiv n+1\mod 2$, not all odd points are in the interior of  $\Delta_{d_1+\ldots+d_n}$, but some are on the boundary. For these points on the boundary, we cannot apply Proposition \ref{prop:number of parallelepipeds}.

\begin{ex}
Figure \ref{figure:non-orientable} shows the intersection of two tropical conics and the dual subdivision of the union of the conics. There are two odd points on the boundary of the Newton polygon of the union. We enrich the two tropical conics by assigning coefficients $\alpha_{(i,j)}$ to a $(i,j)\in \Z^2$. Then the sum over the enriched intersection multiplicities equals 
\[h+\qinv{ \alpha_{(3,1)},-\alpha_{(1,3)}}\]
where $\alpha_{(3,1)}$ is the coefficient of the vertex $(3,1)$ and $\alpha_{(1,3)}$ is the coefficient of the vertex $(1,3)$ in the dual subdivision.
Hence, the sum depends on the choice of coefficients of the enriched tropical conics, but there is always a hyperbolic summand.

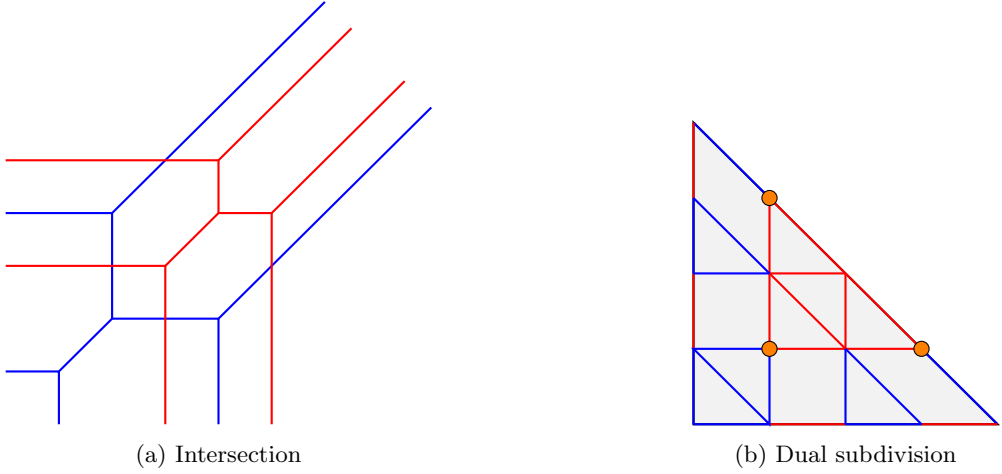
\begin{figure}[t]
     \centering
     \begin{subfigure}[b]{0.45\textwidth}
         \centering
         \begin{tikzpicture}[scale=0.60]
\draw[blue,thick] (-1,0)--(0,0)--(0,-1);
\draw[blue,thick] (0,0)--(1,1)--(3,1)--(3,-1);
\draw[blue,thick] (1,1)--(1,3)--(-1,3);
\draw[blue,thick] (3,1)--(7,5);
\draw[blue,thick] (1,3)--(5,7);
\draw[red,thick] (-1,2)--(2,2)--(2,-1);
\draw[red,thick] (2,2)--(3,3)--(4,3)--(4,-1);
\draw[red,thick] (4,3)--(6.5,5.5);
\draw[red,thick] (3,3)--(3,4)--(5.5,6.5);
\draw[red,thick] (3,4)--(-1,4);
\end{tikzpicture}
\caption{Intersection}
     \end{subfigure}
     \quad
     \begin{subfigure}[b]{0.45\textwidth}
         \centering
         \begin{tikzpicture}[scale=1.1]
\filldraw[fill=gray!10!white,thick] (0,0)--(4,0)--(0,4)--cycle;
\draw[blue, thick] (0,4) -- (1,3);
\draw[red, thick] (1,3)--(3,1);
\draw[blue, thick] (4,0)--(3,1);
\draw[red, thick] (1,1)--(3,1);
\draw[red, thick] (1,1)--(1,3);
\draw[red, thick] (1,2)--(2,2);
\draw[red, thick] (2,1)--(2,2);
\draw[red, thick] (2,1)--(1,2);
\draw[red, thick] (4,0)--(3,0);
\draw[red, thick] (2,0)--(1,0);
\draw[red, thick] (0,4)--(0,3);
\draw[red, thick] (0,2)--(0,1);
\draw[blue, thick] (3,0)--(2,0);
\draw[blue, thick] (1,0)--(0,0);
\draw[blue, thick] (0,3)--(0,2);
\draw[blue, thick] (0,1)--(0,0);
\draw[blue, thick] (3,0)--(2,1);
\draw[blue, thick] (0,3)--(1,2);
\draw[blue, thick] (2,0)--(2,1);
\draw[blue, thick] (0,2)--(1,2);
\draw[blue, thick] (1,0)--(0,1)--(1,1)--cycle;
\fill[orange, draw=black] (1,1) circle (0.1 cm);
\fill[orange, draw= black] (3,1) circle (0.1 cm);
\fill[orange, draw=black] (1,3) circle (0.1 cm);
\end{tikzpicture} 
\caption{Dual subdivision}
     \end{subfigure}
        \caption{Intersection of two tropical conics}
        \label{figure:non-orientable}
\end{figure}
\end{ex}

Let $\tV_1,\ldots,\tV_n$ be enriched tropical hypersurfaces in $\R^n$, with Newton polytopes \linebreak$\Delta_{d_1},\ldots,\Delta_{d_n}$ such that $\sum d_i\not \equiv n+1\mod 2$. 
As suggested in the example, we can find a lower bound for the number of hyperbolic summands in $\sum_p\mult_p(\tV_1,\ldots,\tV_n)$.
For $n$ odd let $m\in \Z$ such that the sum~$d_1+\ldots+d_n=2m$ and for $n$ even let $m\in \Z$ such that the sum $d_1+\ldots+d_n=2m+1$.
Set 
\[N(m)\coloneqq \text{the number of odd points on }\Delta_{d_1+\ldots+d_n}.\]
The following table computes $N(m)$ for $n\le 6$.
\label{section:non orientable}
\[\begin{array}{cc}
d_1=2m+1     & 1 \\
d_1+d_2=2m   & m\\
d_1+d_2 +d_3=2m+1 &\displaystyle \frac{m(m+1)}{2}\\
d_1+\dots +d_4=2m & \displaystyle\frac{m(m+1)(m-1)}{3!}\\
d_1+\dots +d_5=2m+1 & \displaystyle\frac{m(m+1)(m-1)(m+2)}{4!}\\
d_1+\dots +d_6=2m &\displaystyle\frac{m(m+1)(m-1)(m+2)(m-2)}{5!}\\
\end{array}\]

Since the only non-hyperbolic contribution to $\sum_{p\in \tV_1\cap\ldots\cap \tV_n}\mult_p(\tV_1,\ldots,\tV_n)$ comes from the odd points on the boundary, we get the following Corollary.

\begin{coro}[Enriched tropical B\'{e}zout in the non-relatively orientable case]
\label{cor:non-orientable}
Let \linebreak$\tV_1,\ldots,\tV_n$ be enriched tropical hypersurfaces in $\mathbb{R}^n$ with Newton polygons $\Delta_{d_1},\ldots, \Delta_{d_n}$ such that $\sum_{i=1}^nd_i\not\equiv n+1 \mod 2$.
Then 
\[\sum_{\text{intersections }p} \widetilde{\operatorname{mult}}_p(\tV_1,\ldots,\tV_n)=\frac{d_1\cdots d_n-r}{2}\;h+\qinv{ a_1,\ldots,a_{r}}\in\GW(k)\]
where $r$ has to be smaller or equal the number $N(m)$ of odd points on $\partial\Delta_{d_1+\ldots+d_n}$.
\end{coro}

\begin{rmk}
In the non-relatively orientable case McKean shows that one can orient the vector bundle
\[V=\mathcal{O}_{\mathbb{P}^n}(d_1)\oplus\ldots\oplus\mathcal{O}_{\mathbb{P}^n}(d_n)\longrightarrow \mathbb{P}^n\]
relative to a divisor at infinity and compute the $\A^1$-Euler number of $V$ relative to this divisor in the sense of Larson-Vogt \cite{LV}. 
Corollary \ref{cor:non-orientable} gives us a lower bound for the number of hyperbolic forms in this $\A^1$-Euler number by the same argument as in Remark \ref{rmk:BezoutProof}.
\end{rmk}

The lower bound on the number of hyperbolic summands is in Corollary \ref{cor:non-orientable} is not necessarily strict.
For enriched tropical curves we find a better, strict bound.

\begin{coro}[Enriched tropical B\'ezout for curves in the non-relatively orientable case]
Let~$\tC_1$ and $\tC_2$ be two enriched tropical curves of degree $\Delta_{d_1}$ and $\Delta_{d_2}$, respectively, with $d_1+d_2\equiv 0\mod 2$ that intersect tropically transversely. Let $d\coloneqq \min\{d_1,d_2\}$. Then
\[\sum_{p\in \tC_1\cap \tC_2}\mult_p(\tC_1,\tC_2)=\frac{d_1\cdot d_2-d}{2}\;h+\qinv{ a_1,\ldots,a_d}\in \GW(k)\]
for some $a_1,\ldots,a_d\in k^{\times}/(k^\times)^2$.
\end{coro}

\begin{proof}
In case, $d_1+d_2\equiv0\mod 2$, there are $\frac{d_1+d_2}{2}$ odd points on the boundary of $\Delta_{d_1}+\Delta_{d_2}$, all lying on the hypotenuse. To get a non-hyperbolic summand, one of the two edges adjacent to an odd vertex on the hypotenuse has to belong to $\tC_1$ and the other one has to belong to $\tC_2$. 
This can happen at most $d=\min(d_1,d_2)$ times since only $d_1$ segments on the hypotenuse of $\Delta_{d_1}+\Delta_{d_2}$ correspond to edges of $\tC_1$ and $d_2$ correspond to edges of $\tC_2$.
\end{proof}

%We compare these bounds with the inequalities of Eisenbud-Levine-Tessier (short ELT) \cite{EL}: Let $f:\A^n\longrightarrow \A^n$ with an isolated zero at zero. Then Eisenbud and Levine use results of Tessier to find an upper bound for the dimension of the anisotropic subspace $D$ of the local $\A^1$-degree of $f$ at $0$
%\[\dim_kD\le (\operatorname{rank}\deg_x^{\A^1}f)^{1-\frac{1}{n}}.\]

%{\color{red} this is a stricter upper bound than our result :-(} {\color{yellow} not anymore!}{\color{red} I really think this is realated, but I don't know how...}

\subsection{Bernstein-Kushnirenko theorem}
The results above do not restrict to hypersurfaces in $\mathbb{P}^n$. The tools of tropical geometry can be applied to toric varieties, where the action of the torus yields a combinatorial approach to their study. 

\begin{ex}
\label{ex:p1xp1}
Let $C_1$ and $C_2$ be two curves in $\mathbb{P}^1 \times \mathbb{P}^1$ defined by $f_1$ and $f_2$ of bidegree $(d_1,d_2)$ and $(e_1,e_2)$, respectively. Then $f_1$ and $f_2$ define a section of 
\[V\coloneqq \glob(d_1,d_2)\oplus \glob(e_1,e_2)\longrightarrow\mathbb{P}^1\times \mathbb{P}^1 \]
where $\glob(a,b)=\pi_1^*\glob_{\mathbb{P}^1}(a)\otimes \pi_2^*\glob_{\mathbb{P}^1}(b)$ and $\pi_i:\mathbb{P}^1\times \mathbb{P}^1\longrightarrow\mathbb{P}^1$ is the $i$th projection for $i=1,2$.
The vector bundle $V$ is relatively orientable if and only if 
$\det V\otimes \omega_{\mathbb{P}^1\times \mathbb{P}^1}$, that is, $\glob(d_1+e_1-2,d_2+e_2-2)$ is a square. This is the case if and only if both $d_1+e_1$ and $d_2+e_2$ are even.
The Newton polygons $\NP(f_1)$ and $\NP(f_2)$ are rectangles with corners $(0,0),(d_1,0),(0,d_2),(d_1,d_2)$, respectively, with corners $(0,0),(e_1,0),(0,e_2),(e_1,e_2)$. The Minkowski sum $\NP(f_1)+\NP(f_2)$ is the rectangle with corners $(0,0),(d_1+e_1,0),(0,d_2+e_2)$, $(d_1+e_1,d_2+e_2)$. This rectangle $\NP(f_1)+\NP(f_2)$ has no odd vertices on the boundary if and only if both $d_1+e_1$ and $d_2+e_2$ are even, that is exactly when the vector bundle $V$ is relatively orientable.
Let $\tC_1$ and $\tC_2$ be two enriched tropical curves with Newton polytopes equal to~$\NP(f_1)$ and $\NP(f_2)$ and assume that both $d_1+e_1$ and $d_2+e_2$ are even.
Then Proposition \ref{prop:number of parallelepipeds} implies that 
\begin{align*}\sum_{p\in C_1\cap C_2}\mult_p(\tC_1,\tC_2)=\frac{d_1e_2+d_2e_1}{2}\;h\in \GW(k).
\end{align*}
Equivalently, we get that the $\A^1$-Euler number of the vector bundle $V$ equals
\begin{align*}n^{\A^1}(V)&=\sum_{p\in C_1\cap C_2}\operatorname{ind}_p(f_1,f_2)\\&=\frac{\operatorname{Area}(\NP(f_1)+\NP(f_2))-\operatorname{Area}(\NP(f_1))-\operatorname{Area}(\NP(f_2))}{2}\;h\\&=\frac{d_1e_2+d_2e_1}{2}\;h\in \GW(k)
\end{align*}
which yields an enriched count of intersection points of two curves in $\mathbb{P}^1\times \mathbb{P}^1$.
\end{ex}
%%%%%%%%%%%%%%%%%%%%%%%%%%%%%%%
\begin{ex}
\label{ex:Hirzebruch}
Let $C_1$ and $C_2$ be two curves in the Hirzebruch surface $\Sigma_n$ defined by $f_1$ and $f_2$ of bidegree $(a_1,b_1)$ and $(a_2,b_2)$, respectively, where $a_i=C_i\cdot F, b_i=C_i\cdot E+na_i$, for $F$ a generic fiber and $E$ the exceptional divisor of $\Sigma_n$. Then $f_1$ and $f_2$ define a section of 
\[V\coloneqq \glob(a_1 E+b_1 F)\oplus \glob(a_2 E+b_2 F) \longrightarrow\Sigma_n. \]
The vector bundle $V$ is relatively orientable if and only if 
\[\det V\otimes \omega_{\Sigma_n}=\glob\left((a_1+a_2-2)E+(b_1+b_2-(n+2)F)\right)\] is a square, which is the case if and only if $a_1+a_2$ is even and $b_1+b_2\equiv n\mod2$.
The Newton polygons $\NP(f_1)$ and $\NP(f_2)$ are trapezia with corners $(0,0),(a_1 n+b_1,0),(b_1,a_1),(0,a_1)$ and~$(0,0),(a_2 n+b_2,0),(b_2,a_2),(0,a_2)$, respectively. 
The Min\-kowski sum $\NP(f_1)+\NP(f_2)$ is the trapezium with corners $(0,0),((a_1+a_2) n+b_1+b_2,0),(b_1+b_2,a_1+a_2),(0,a_1+a_2)$. This trapezium~$\NP(f_1)+\NP(f_2)$ has no odd vertices on the boundary if and only if $a_1+a_2$ is even and~$b_1+b_2\equiv n\mod2$, that is exactly when $V$ is relatively orientable.
Let $\tC_1$ and $\tC_2$ be two tropical curves with Newton polytopes equal to $\NP(f_1)$ and $\NP(f_2)$ and assume that $a_1+a_2$ is even and $b_1+b_2\equiv n\mod2$.
Then Proposition \ref{prop:number of parallelepipeds} implies that 
\begin{align*}\sum_{p\in \tC_1\cap \tC_2}\mult_p(\tC_1,\tC_2)=\frac{a_1a_2n+a_1b_2+a_2b_1}{2}\;h\in \GW(k).
\end{align*}
Equivalently, as in the previous example, this coincides with $n^{\A^1}(V)$ and yields an enriched count of intersection points of two curves in $\Sigma_n$.
\end{ex}

The examples above as well as B\'ezout's theorem are special cases of a quadratic enrichment of the Bernstein-Kushnirenko theorem.
We recall the classical statement of this theorem.
Let $A$ be a finite subset of $\mathbb{Z}^n$ and let 
\[L_A\coloneqq \left\{f\middle| f(x)=\sum_{I\in A}c_Ix^I=\sum_{I\in A}c_Ix_1^{I_1}\ldots x_n^{I_n}\text{, }c_I\in k\right\}\] 
be the space of Laurent polynomials whose exponents are in $A$. 
Let $\Delta_A$ be the convex hull of the points in $A$.
The classical Bernstein-Kushnirenko theorem says.
\begin{thm}[Bernstein-Kushnirenko theorem]
\label{thm:classicalBK}
For $n$ finite subsets $A_1,\ldots,A_n$ of $\mathbb{Z}^n$ and for a generic system of equations 
\[f_1(x)=\ldots=f_n(x)=0\]
where $f_i\in L_{A_i}$, the number of solutions in $(\C\setminus\{0\})^n$ equals the mixed volume
\[\operatorname{MVol}(\Delta_{A_1},\ldots,\Delta_{A_n}).\]
\end{thm}

Before we state the quadratically enriched version of this theorem, we define the following condition that can be seen as the combinatorial analogue of relative orientability.
\begin{df}
We say that the tuple $(A_1,\ldots,A_n)$ is \emph{combinatorially oriented} if the Minkowski sum $\Delta_{A_1}+\ldots+\Delta_{A_n}$ has no odd points on the boundary.
\end{df}

\begin{ex}[B\'ezout theorem]
\label{ex:bezout combi oriented}
Let $A_i=\Delta_{d_i}\cap \Z^n$ for some positive integer $d_i$, i.e., $\Delta_{A_i}=\Delta_{d_i}$, for~$i=1,\ldots,n$. 
Then $(A_1,\ldots,A_n)$ is combinatorially oriented if and if~$\Delta_{d_1}+\ldots+\Delta_{d_n}$ has no odd boundary points which is exactly the case if $d_1+\ldots+d_n\equiv n+1\mod 2$, i.e., exactly when the vector bundle  
\[\mathcal{O}_{\mathbb{P}^n}(d_1)\oplus \ldots \oplus \mathcal{O}_{\mathbb{P}^n}(d_n)\longrightarrow \mathbb{P}^n\]
is relatively orientable.
\end{ex}

In all examples \ref{ex:bezout combi oriented}, \ref{ex:p1xp1} and \ref{ex:Hirzebruch} above the condition of being combinatorially oriented coincides with the condition for the corresponding vector bundle to be relatively orientable which motivates the following conjecture. 

\begin{conj}
Let $X$ be a smooth toric variety of dimension $n$ and let $f_1,\ldots,f_n$ be regular functions on $X$ such that $f_i$ is a non-trivial section of a line bundle $\mathcal{L}_i\longrightarrow X$ such that the system~$f_1=\dots=f_n=0$ has a non-empty solution set formed of isolated zeros. Furthermore, let~$A_i\coloneqq \{\text{exponents of }f_i\}\subset \Z^n$. Then $(A_1,\ldots,A_n)$ is combinatorially oriented if and only if the vector bundle $\mathcal{L}_1\oplus \ldots\oplus \mathcal{L}_n\longrightarrow X$ is relatively orientable.
\end{conj}

Let us say that a variety has the \emph{combinatorial orientability property} if it this conjecture holds for every sum of $\dim X$ line bundles satisfying the hypothesis. We prove in the following theorem that the class of varieties that satisfy this conjecture is closed under products. In particular, this property holds on products of projective spaces and Hirzebruch surfaces.

\begin{thm}
If $X_1$ and $X_2$ are smooth toric varieties that satisfy the combinatorial orientability property, then $X_1\times X_2$ also satisfies the combinatorial orientability property.
\end{thm}

\begin{proof}
The product $X\coloneqq X_1\times X_2$ has a toric structure given by the product of the toric structures on each component. Since $X$ is smooth, we have that
\[\operatorname{Pic}(X)\simeq H^2(X,\Z)=H^2(X_1,\Z)\oplus H^2(X_2,\Z).\]
by the K\"unneth formula and the fact that $H^1(X_i,\Z)=0$ due to $X_i$ being a smooth toric variety, $i=1,2$.
Hence, through this isomorphism, every line bundle $\cL$ is determined by its bidegree $\bar{d}=(d_1,d_2)$, where each degree class $d_i\in H^2(X_i,\Z), i=1,2$. Therefore, for every line bundle $\cL$ over $X$ of degree $\bar{d}$, there are line bundles $\cL^1$ and $\cL^2$ over $X_1$ and $X_2$ of degree $d_1$ and $d_2$, respectively, such that
\[\cL=p_1^*\cL^1\otimes p_2^*\cL^2,\] where $p_i:X\lra X_i, i=1,2$ is the component projection. 
Given this decomposition, the line bundle $\cL$ is a square if and only if each $\cL^i$, $i=1,2$, is a square.
Moreover, the polytope $\Delta$ associated to $\cL$ in $\Lambda\otimes\R$ is the product of the polytopes $\Delta_1$ and $\Delta_2$ associated to $\cL_1$ and $\cL_2$ in~$\Lambda_1\otimes\R$ and $\Lambda_2\otimes\R$,  where $\Lambda_1$ and $\Lambda_2$ are the lattices associated to $X_1$ and $X_2$, respectively, and the lattice $\Lambda\coloneqq\Lambda_1\times\Lambda_2$ is the one associated to $X$. The polytope~$\Delta$ has boundary 
\[\partial(\Delta)=\partial(\Delta^1\times\Delta^2)=(\partial(\Delta^1)\times\Delta^2)\cup(\Delta^1\times\partial(\Delta^2)),\]
and so, the odd lattice points are
$\partial(\Delta)^{\mathrm{odd}}=(\partial(\Delta^1)^{\mathrm{odd}}\times(\Delta^2)^{\mathrm{odd}})\cup((\Delta^1)^{\mathrm{odd}}\times\partial(\Delta^2)^{\mathrm{odd}}).$
Now, let~$V$ be the vector bundle $\mathcal{L}_1\oplus \ldots\oplus \mathcal{L}_n\longrightarrow X$, where $n=\dim X$, $L_i$ is a line bundle over~$X$ and~$A_i$ is the exponent set of a generic section $f_i$ of $\mathcal{L}_i,i=1,\dots,n$.
Put $\Delta_i=\operatorname{Conv}(A_i)$ the convex hull of the set $A_i$ in $\Lambda\otimes\R.$ Assume that the system $\left\{f_i=0\right\}_{i=1,\dots,n}$ has an isolated zero. In this case we have that $(\sum_{i=1}^n\Delta_i)^{\mathrm{odd}}\neq\emptyset$, otherwise there would be an $i_0$ for which $\Delta_{i_0}=\{\mathrm{pt}\}$ or there would be a vector subspace $H\subset\Lambda\otimes\R$ of lower dimension, containing all $\Delta_i\subset H$, which contradicts the fact that the system has only isolated zeros. This implies that $(\sum_{i=1}^n\Delta_i^1)^{\mathrm{odd}}\neq\emptyset$ and~$(\sum_{i=1}^n\Delta_i^2)^{\mathrm{odd}}\neq\emptyset$.
Lastly, since the Minkowski sum commutes with products, the odd lattice points in the boundary of the Minkowski sum satisfy
\[
\partial(\sum_{i=1}^n\Delta_i)^{\mathrm{odd}}=(\partial(\sum_{i=1}^n\Delta_i^1)^{\mathrm{odd}}\times(\sum_{i=1}^n\Delta_i^2)^{\mathrm{odd}})\cup((\sum_{i=1}^n\Delta_{i=1}^1)^{\mathrm{odd}}\times\partial(\sum_{i=1}^n\Delta_i^2)^{\mathrm{odd}}).
\]
These facts imply our statement. Namely, if the $n$-tuple $(A_1,\ldots,A_n)$ is combinatorially oriented, then $\partial(\sum_{i=1}^n\Delta_i)^{\mathrm{odd}}=\emptyset$
by definition. Since $(\sum_{i=1}^n\Delta_i^1)^{\mathrm{odd}}\neq\emptyset, (\sum_{i=1}^n\Delta_i^2)^{\mathrm{odd}}\neq\emptyset$ in this case, we have that
the $n$-tuple $(A_1,\ldots,A_n)$ is combinatorially oriented if and only if both of the sets~$\partial(\sum_{i=1}^n\Delta_i^1)^{\mathrm{odd}}$ and $\partial(\sum_{i=1}^n\Delta_i^2)^{\mathrm{odd}}$ are empty. Since $X_1$ and $X_2$ satisfy the combinatorial orientability property, the sets~$\partial(\sum_{i=1}^n\Delta_i^1)^{\mathrm{odd}}$ and $\partial(\sum_{i=1}^n\Delta_i^2)^{\mathrm{odd}}$ are empty if and only if the vector bundles given by the direct sum of the components of each of the line bundles~${V^1\coloneqq\mathcal{L}_1^1\oplus \ldots\oplus \mathcal{L}_n^1\longrightarrow X_1}$ and
$V^2\coloneqq\mathcal{L}_1^2\oplus \ldots\oplus \mathcal{L}_n^2\longrightarrow X_2$ are relatively orientable. Finally, the vector bundles $V_1$ and $V_2$ are relatively orientable if and only the vector bundle $V$ is relatively orientable since $\det V=p_1^*\det V^1\otimes p_2^*\det V^2$ and $\omega_X=p_1^*\omega_{X_1}\otimes p_2^*\omega_{X_2}$, so 
\[\det V\otimes \omega_X=p_1^*(\det V^1\otimes\omega_{X_1})\otimes p_2^*(\det V^2\otimes\omega_{X_2})\] 
is a square if and only if $\det V^1\otimes\omega_{X_1}$ and $\det V^2\otimes\omega_{X_2}$ are squares.
\end{proof}

\begin{ex}
\label{ex:p1n}
Let $V_1,V_2,\dots,V_n$ be hypersurfaces of $(\mathbb{P}^1)^n$ defined by $f_i$, of multidegree $(d_1^i,d_2^i,\dots,d_n^i)$, for~$i=1,2,\dots, n$. Then $(f_1,f_2,\dots,f_n)$ defines a section of 
\[V\coloneqq\bigoplus_{i=1}^n \glob(d_1^i,d_2^i,\dots,d_n^i)\longrightarrow(\mathbb{P}^1)^n\]
where $\glob(d_1^i,d_2^i,\dots,d_n^i)\coloneqq\bigotimes_{j=1}^n \pi_j^*\glob_{\mathbb{P}^1}(d_j^i)$ and $\pi_j:(\mathbb{P}^1)^n\longrightarrow\mathbb{P}^1$ is the $j$th projection.
The vector bundle $V$ is relatively orientable if and only if 
\[\det V\otimes \omega_{(\mathbb{P}^1)^n}=\glob\left(\sum_{i=1}^n d_1^i-2,\sum_{i=1}^n d_2^i-2,\dots,\sum_{i=1}^n d_n^i-2\right)\] is a square, which is the case if and only if $\sum_{i=1}^n d_j^i-2$ even for every $j=1,2,\dots,n$.
For every~$i=1,2,\dots, n$, the Newton polygon $\NP(f_i)$ is the parallelepiped with a corner in $\bar{0}$ and side edges~$d_j^i\mathrm{e}_j$, where $\{\mathrm{e}_j\}_{j=1}^n$ is the standard basis. The Minkowski sum $\sum_{i=1}^n\NP(f_i)$ is the parallelepiped with a corner in $\bar{0}$ and side edges $\sum_{i=1}^n d_j^i\mathrm{e}_j$. This parallelepiped $\sum_{i=1}^n\NP(f_i)$ has no odd vertices on the boundary if and only if every  $\sum_{i=1}^n d_j^i, j=1,2,\dots,n$  is even, that is exactly when $V$ is relatively orientable.
Let $\tV_i,i=1,2,\dots,n$ be enriched tropical hypersurfaces with Newton polytope equal to~$\NP(f_i)$ and assume that $\sum_{i=1}^n d_j^i, j=1,2,\dots,n$ is even. Then, Proposition~\ref{prop:number of parallelepipeds} implies that 
\begin{align*}\sum_{p\in \bigcap_{i=1}^n V_i}\mult_p(\tV_1,\tV_2,\dots, \tV_n)=\frac{1}{2}\left(\sum_{\sigma\in\mathcal{S}_n}d_{\sigma(1)}^1d_{\sigma(2)}^2\cdots d_{\sigma(n)}^n\right)h=n^{\A^1}(V)\in \GW(k).
\end{align*}

Equivalently, this coincides with $n^{\A^1}(V)$ and yields an enriched count of intersection points of~$n$ curves in $(\mathbb{P}^1)^n$.

\end{ex}

For $(A_1,\ldots,A_n)$ combinatorially oriented, we get that the enriched count of zeros of the system of equations $f_1=\ldots=f_n=0$ is independent of the coefficients of $f_1,\ldots,f_n$.
%{\color{red} for the next two theorems I added a characteristic assumption, if the toric varieties are actually defined over $\Z$ we do not need this...}
\begin{thm}[Enriched Bernstein-Kushnirenko theorem]
\label{thm:Enriched Bernstein-Kushnirenko}
Assume $k$ is a field with $\operatorname{char}k=0$ or $\operatorname{char}k>\max_i \{\operatorname{diam}(\Delta_{A_i})\}$.
For a combinatorially oriented $n$-tuple of indexing sets $(A_1,\ldots,A_n)$ and for a generic system of equations
\[f_1(x)=\ldots=f_n(x)=0\]
where $f_i\in L_{A_i}$, the enriched count of solutions $z$ in $\operatorname{Spec} k[x_1^{\pm1},\ldots,x_n^{\pm 1}]$ equals 
\[\sum_z\Tr_{\kappa(z)/k}\qinv{ \det \operatorname{Jac}(f_1,\ldots,f_n)(z)})=\frac{\operatorname{MVol}(\Delta_{A_1},\ldots,\Delta_{A_n})}{2}\;h\in \GW(k).\]
\end{thm}

We can also say something about the non-orientable case. Just like in the case of B\'ezout, then the enriched count depends on the choice of coefficients of $f_1,\ldots,f_n$. However, we still get a lower bound for the number of hyperbolic summands depending on the number of odd points on the boundary of $\Delta_{A_1}+\ldots+\Delta_{A_n}$.

\begin{thm}
\label{thm:BK non-orientable}
Assume $k$ is a field with $\operatorname{char}k=0$ or $\operatorname{char}k>\max_i \{\operatorname{diam}(\Delta_{A_i})\}$.
Let $(A_1,\ldots, A_n)$ be a sequence of finite subsets of $\Z^n$. 
Let 
\[N\coloneqq \operatorname{Card}(\partial(\Delta_{A_1}+\ldots+\Delta_{A_n})\cap\Lambda^{\text{odd}})\] be the number of odd points on the boundary of $\Delta_{A_1}+\ldots+\Delta_{A_n}$.
Then for a generic system of equations 
\[f_1(x)=f_2(x)=\ldots=f_n(x)=0\]
where $f_i\in L_{A_i}$, we get that 
the enriched count of solutions $z$ in $\operatorname{Spec}k[x_1^{\pm1},\ldots,x_n^{\pm 1}]$
equals
\begin{equation}
\label{eq:BernsteinKushnirenko}
\sum_z\Tr_{\kappa(z)/k}\qinv{\det \operatorname{Jac}(f_1,\ldots,f_n)(z)}=\frac{\operatorname{MVol}(\Delta_{A_1},\ldots,\Delta_{A_n})-r}{2}\;h+\qinv{ a_1,\ldots,a_r}\end{equation}
in $\GW(k)$, for some $r\le N$ and $a_1,\ldots,a_r\in k^\times$.
\end{thm}
\begin{proof}[Proof of Theorem \ref{thm:Enriched Bernstein-Kushnirenko} and Theorem \ref{thm:BK non-orientable}]

%et's assume that first that our base field $k$ is actually the field of Puiseux series $k\Puiseux$ over $k$ and that $\operatorname{char}k=0$ or $\operatorname{char} k>\max_i\{\operatorname{diam}(\Delta_{A_i})\}$.
The argument is the same as in Corollary \ref{coro:EnrichedTropicalBezout}. We know by Theorem \ref{thm:tropical Bezout} that $\sum_z\Tr_{\kappa(z)/k\Puiseux}\qinv{\det \operatorname{Jac}(f_1,\ldots,f_n)(z)}$ is an element of $\GW(k)$ of rank $\operatorname{MVol}(\Delta_{A_1},\ldots,\Delta_{A_n})$. If there are no odd points on the boundary of the Minkowski sum $\Delta_{A_1}+\ldots+\Delta_{A_n}$ then by the same argument as in Corollary $\ref{coro:EnrichedTropicalBezout}$ we can have at most $N$ summands in $\sum_z\Tr_{\kappa(z)/k\Puiseux}\qinv{\det \operatorname{Jac}(f_1,\ldots,f_n)(z)}$ that are not multiples of $h$.

%{\color{blue} ask Stephen if this is defined over $\Z$, need that the toric variety is smooth and proper over $\Z$}
To derive the theorem for $k$, we use the natural isomorphism $\GW(k)\cong \GW(k\Puiseux)$ from Example \ref{ex:GW of Puiseux}.%{\color{red} Should we say more here?}
\end{proof}

We provide examples where \eqref{eq:BernsteinKushnirenko} is invariant, i.e., does not depend of the coefficients of the $f_i$ and where it depends on the coefficients.

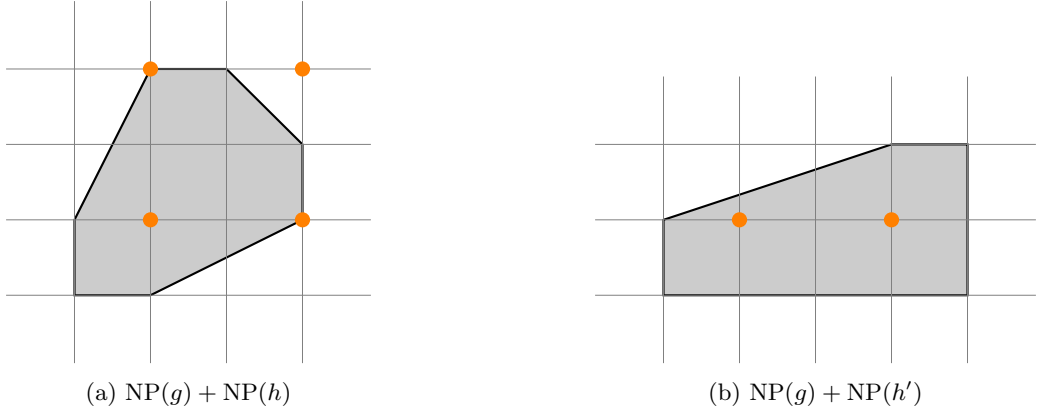
\begin{figure}[t]
     \centering
     \begin{subfigure}[b]{0.45\textwidth}
         \centering
         \begin{tikzpicture}
\filldraw[fill=gray!40!white,thick] (0,0)--(1,0)--(3,1)--(3,2)--(2,3)--(1,3)--(0,1)--cycle;
\draw[step=1cm,gray,very thin] (-0.9,-0.9) grid (3.9,3.9);
\fill[orange] (1,1) circle (0.1 cm);
\fill[orange] (3,1) circle (0.1 cm);
\fill[orange] (1,3) circle (0.1 cm);
\fill[orange] (3,3) circle (0.1 cm);
\end{tikzpicture}
\caption{$\NP(g)+\NP(h)$}
     \end{subfigure}
     \hfill
     \begin{subfigure}[b]{0.45\textwidth}
         \centering
         \begin{tikzpicture}
\filldraw[fill=gray!40!white,thick] (0,0)--(4,0)--(4,2)--(3,2)--(0,1)--cycle;
\draw[step=1cm,gray,very thin] (-0.9,-0.9) grid (4.9,2.9);
\fill[orange] (1,1) circle (0.1 cm);
\fill[orange] (3,1) circle (0.1 cm);
\end{tikzpicture}
\caption{$\NP(g)+\NP(h')$}
     \end{subfigure}
        \caption{Newton polygons}
         \label{figure:MinkowskiSum}
\end{figure}

\begin{ex}
The following is the leading example in \cite{Sturmfels}.
Let 
\[g(x,y)=a_1+a_2x+a_3xy+a_4y\]
and 
\[h(x,y)=b_1+b_2x^2y+b_3xy^2.\]
By Theorem \ref{thm:classicalBK} the number of solutions in $\operatorname{Spec}\C[x^{\pm1},y^{\pm 1}]$ to $g(x,y)=h(x,y)=0$ is equals the mixed volume $\MV(\NP(g),\NP(h))$ of the Newton polygons $\NP(g)$ and $\NP(h)$ of $g$ and $h$ for a generic choice of coefficients $a_1,a_2,a_3,a_4,b_1,b_2,b_3$.
The enriched count of solutions to $g(x,y)=0$ and~$h(x,y)=0$ depends on the choice of coefficients since there are $2$ odd points on the boundary of $\NP(g)+\NP(h)$ as shown the first picture in Figure \ref{figure:MinkowskiSum}. For example if $k=\R$ and $a_1=a_2=a_3=a_3=a_4=b_1=b_2=b_3=1$ we get that the enriched count of zeros equals $2h$. However, if we set $b_3=-1$ (all other coefficients are still equal to $1$), then the enriched count of zeros equals $h+\qinv{1,1}$.

If we replace $h$ by 
$h'(x,y)=b_1+b_2x^3+b_3x^3y$, the Minkowski sum $\NP(g)+\NP(h')$ of the Newton polygons $\NP(g)$ and $\NP(h')$ of $g$ and $h'$ has no odd points on the boundary as one can see in the second picture of Figure \ref{figure:MinkowskiSum}. The mixed volume of $\NP(g)$ and $\NP(h')$ equals $4$, hence the enriched count of solutions to $g(x,y)=h'(x,y)=0$ in $\operatorname{Spec}(k[x^{\pm 1},y^{\pm1}])$ equals $2h$ and this is independent of the choice of coefficients.

\end{ex}